\documentclass{IEEEtran4PSCC}

\usepackage{cite}
\usepackage[cmex10]{amsmath}
\usepackage{algorithmic}
\usepackage{url}
\usepackage{subcaption} 
\usepackage{supertabular}
\usepackage{eurosym}
\usepackage{pifont}
\usepackage{amssymb}
\usepackage{xcolor}

\ifCLASSINFOpdf
\usepackage[pdftex]{graphicx}
\else
\usepackage[dvips]{graphicx}
\fi
\makeatletter
\let\old@ps@headings\ps@headings
\let\old@ps@IEEEtitlepagestyle\ps@IEEEtitlepagestyle
\def\psccfooter#1{%
    \def\ps@headings{%
        \old@ps@headings%
        \def\@oddfoot{\strut\hfill#1\hfill\strut}%
        \def\@evenfoot{\strut\hfill#1\hfill\strut}%
    }%
    \def\ps@IEEEtitlepagestyle{%
        \old@ps@IEEEtitlepagestyle%
        \def\@oddfoot{\strut\hfill#1\hfill\strut}%
        \def\@evenfoot{\strut\hfill#1\hfill\strut}%
    }%
    \ps@headings%
}
\makeatother

\usepackage[draft=False]{hyperref} % put False when you've finished all your writing
\hypersetup{
     colorlinks=true,
     linkcolor=blue,
     filecolor=blue,
     citecolor=black,      
     urlcolor=cyan,
     }

\begin{document}

\title{Stochastic optimization for unit commitment applied to the security of supply: extended version}

%% To specify the authors when (number of affiliations <= 2)
\author{
\IEEEauthorblockN{Jonathan Dumas\\}
\IEEEauthorblockA{RTE Research \& Development, Paris, France\\
jonathan.dumas@rte-france.com}
}

%% To specify the authors when (number of affiliations > 2)
% \author{\IEEEauthorblockN{Author n.1\IEEEauthorrefmark{1},
% Author n.2\IEEEauthorrefmark{2},
% Author n.3\IEEEauthorrefmark{3}, 
% Author n.4\IEEEauthorrefmark{3} and
% Author n.5\IEEEauthorrefmark{4}}
% \IEEEauthorblockA{\IEEEauthorrefmark{1} Department Name of Organization A\\
% Name of the organization A,
% Address A\\ Emails if wanted}
% \IEEEauthorblockA{\IEEEauthorrefmark{2} Department Name of Organization B\\
% Name of the organization B,
% Address B\\ Emails if wanted}
% \IEEEauthorblockA{\IEEEauthorrefmark{3} Department Name of Organization C\\
% Name of the organization C,
% Address C\\ Emails if wanted}
% \IEEEauthorblockA{\IEEEauthorrefmark{4}Department Name of Organization D\\
% Name of the organization D,
% Address D\\ Emails if wanted}
% }

% make the title area
\maketitle

\begin{abstract}
Transmission system operators employ reserves to deal with unexpected variations of demand and generation to guarantee the security of supply. The French transmission system operator RTE dynamically sizes the required margins using a probabilistic approach relying on continuous forecasts of the main drivers of the uncertainties of the system imbalance and a 1\% risk threshold. However, this criterion does not specify which means to activate upward/downward and when to face a deficit of available margins versus the required margins. Thus, this work presents a strategy using a probabilistic unit commitment with a stochastic optimization-based approach, including the fixed and variable costs of units and the costs of lost load and production. The abstract problem is formulated with a multi-stage stochastic program and approximated with a heuristic called \textit{two-stage stochastic model predictive control}. It solves a sequence of two-stage stochastic programs to conduct the central dispatch. An implementation is conducted by solving an approximated version with a single two-stage stochastic program.
This method is tested on a real case study comprising nuclear and fossil-based units with French electrical consumption and renewable production.
\end{abstract}


% Use this to place sponsorships
%\thanksto{\noindent Submitted to the 22nd Power Systems Computation Conference (PSCC 2022).}

\section{Introduction}

\subsection{Context and motivations}

% Context
Transmission System Operators (TSOs) procure ancillary services from balancing service providers (BSPs) to ensure operational security \cite{code-eb-A23}. However, the increasing integration of renewable energy sources has raised the uncertainties TSOs face in the daily operation of power grids. TSOs manage remaining imbalances in the system by employing contracted and non-contracted power reserves supplied by BSPs. In particular, TSOs can activate up or downward operating reserves which are defined, under current European legislation \cite{eu-guideline}, by the System Operation Guidelines: Frequency Containment Reserve (FCR), automatic Frequency Restoration Reserve (aFRR), manual Frequency Restoration Reserve (mFRR), and Replacement Reserve (RR).

European TSOs use various methods to secure sufficient capacity to balance the system. They can be categorized into a fully \textit{contracted strategy}, and a \textit{dynamic margin monitoring} approach, presented in \cite{dumas2023dynamic}.
The operational approach adopted by RTE for dynamically sizing the required margins in the context of the dynamic margin monitoring strategy uses a short-term security threshold. 
Currently, France encloses two security of supply criteria. 
% Critere long terme
First, a reliability target as defined in Regulation (EU) 2019/943 Article 25 \cite{code-eb-A25}, named "\textit{the three-hour failure criterion}" set to three hours, including two hours of loss of load per year in terms of expectancy. It is used to size and design the power system in prospective studies, also called adequacy reports\footnote{\url{https://www.rte-france.com/en/analyses-trends-and-perspectives/projected-supply-estimates}}.
% Critere court terme
Second, a short-term criterion, named "\textit{the one \% criterion}". RTE currently uses it to size the required margins \cite{dumas2023dynamic} and aFRR to cover at least 99 \% of the French imbalances following Regulation (UE) 2017/1485 Article 157 \cite{code-eb-A157}.

%%%%
% Limitations of the short-term criterion
% 1. why 1 % ?
% 2. actual approach does not specify which unit to activate and when
An RTE internal technical report \cite{bienvenu2019} pointed out several research directions to question the value and the concept of the short-term criterion. Indeed, the one \% value was set 50 years ago in a context where significant system imbalances occurred mainly in winter during the consumption peaks, with an essential need for upward margins. This one \% value was decided using technical considerations and not a technical-economical approach such as optimization. Then, the increasing penetration of renewable energies, the structural changes in consumption patterns, and the energy mix with less and less fossil-based conventional power plants have been modifying the sources of system imbalances. Finally, the short-term criterion does not specify which unit to activate upward or downward and when to face a deficit of available margins versus the required margins.
Therefore, the present work proposes a techno-economic strategy using a stochastic optimization-based approach, including the fixed and variable costs of units, their technical constraints, and the costs of lost load and production (spillage). This strategy involves solving a \textit{probabilistic unit commitment} (UC). 

\begin{figure}[tb]
\centerline{\includegraphics[width=90mm]{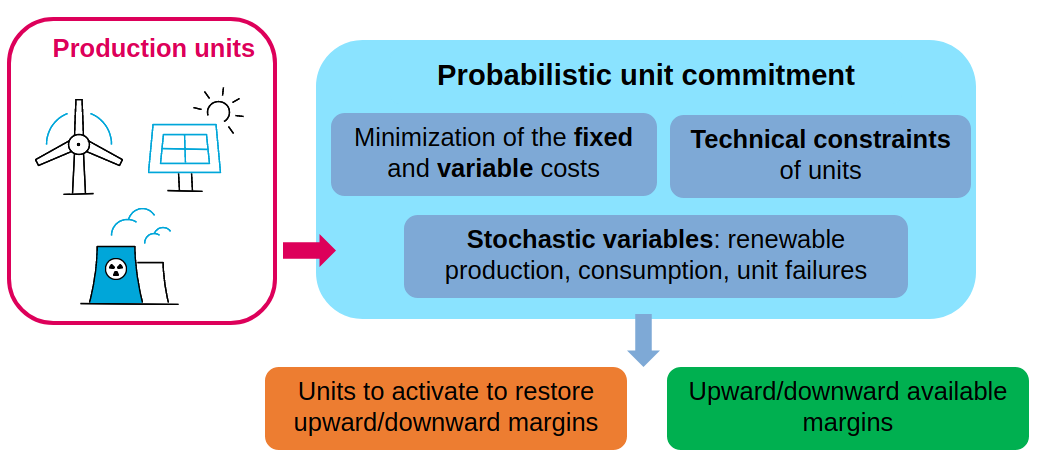}}
\caption{Probabilistic unit commitment approach. This UC tool has two main purposes: 
i) decide which units activate to make upward/downward variations to restore available margins if needed in the case of deficit margins;
ii) provide a techno-economical approach to investigate the value of the short-term criterion by comparing the available upward/downward margins derived from the production plan of the stochastic optimizer to the required margins computed with the one \% criterion.
This study uses this UC tool for the first application.}
\label{fig:graphical-abstract}
\end{figure}

% Probabilistic UC purposes
This probabilistic UC tool, depicted by Figure \ref{fig:graphical-abstract}, has two main purposes:
i) decide which units activate to make upward/downward variations to restore available margins if needed in the case of deficit margins;
ii) provide a techno-economical approach to investigate the value of the short-term criterion. Indeed, by design, the stochastic optimizer builds an optimal dispatch under uncertainty. Then, it is possible to derive the available upward/downward margins from the computed dispatch and to compare them to the required margins computed with the one \% criterion \cite{dumas2023dynamic}. If they exceed the required margins, the short-term criterion's value should be decreased. If smaller, the short-term criterion's value should be increased.
This study aims to design, implement, and test this probabilistic UC on a case study for the first purpose. Further works will be dedicated to using this tool for the second presented purpose.

\subsection{Unit commitment}

%%%%%%%%%%%%%%%%%%%
% UC related works
%%%%%%%%%%%%%%%%%%%%%%%
The unit commitment problem aims to determine the optimal day-ahead or intraday commitment of generators to operate the system at the minimum expected cost. This problem considers the fixed cost of committing generators and the variable cost of dispatching them based on realized uncertainty. 
Stochastic optimization has been employed in UC problems to address the intermittent and stochastic nature of renewable energy and to solve the resulting large-scale computation problems.
Stochastic UC is one of the many applications in the field of uncertainty optimization, including stochastic optimization, chance-constrained optimization, robust optimization, and distributionally robust optimization \cite{ROALD2023108725}.
Two types of stochastic models, the \textit{two-stage} stochastic UC model and the \textit{multi-stage} stochastic UC model, have been studied. 
%
% Approche 1: two-stage
In a two-stage stochastic UC model, the UC decision is the first-stage decision, determined before the uncertainty is realized. Thus, it is not adaptive to specific uncertainty realizations.
% Approche 2: multi-stage
In contrast, a multi-stage stochastic UC model copes with uncertainty dynamically as the UC decision is a function of the realization of load and renewable generation. Therefore, the UC decision adapts to uncertainty realizations, but it may be at the expense of significant computational costs.
We provide a few references for two-stage and multi-stage stochastic UC in the following.

\subsection{Related work on stochastic unit commitment}
 
\begin{figure}[tb]
\centerline{\includegraphics[width=90mm]{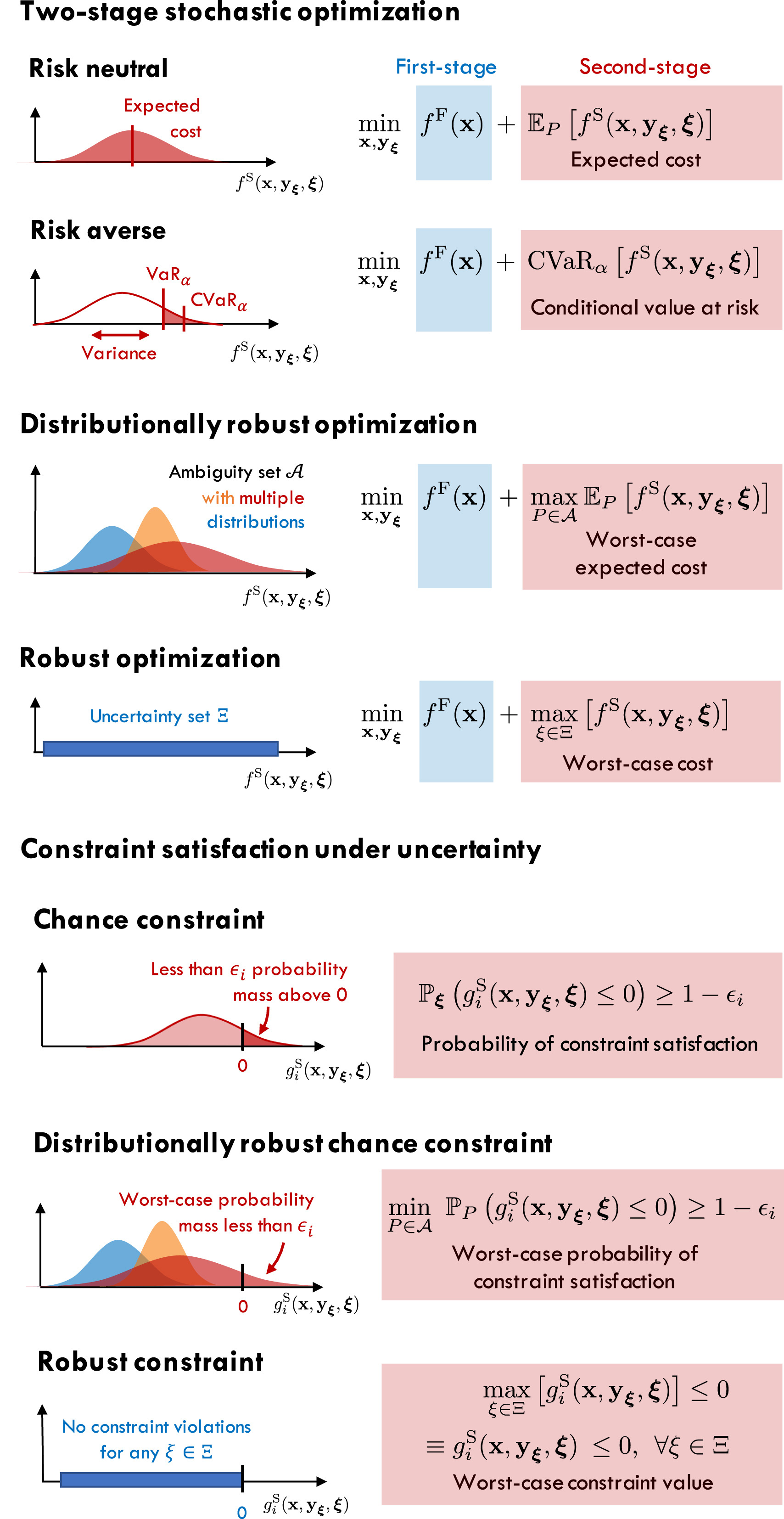}}
\caption{An overview of formulations for power system optimization under uncertainty. This Figure is presented in \cite{ROALD2023108725}.}
\label{fig:sto_stage_strategies}
\end{figure}

%%%%%%%%%%%%%%%%%%%
% Two-stage UC formulations
%%%%%%%%%%%%%%%%%%%%%%%
In the first stage of a two-stage stochastic UC, units are committed, and uncertainty is revealed through realized forecast errors and equipment outages. In the second stage, the system can react by dispatching generators to balance the system while respecting network constraints.
The study \cite{ROALD2023108725} provides an overview of different ways to formulate two-stage optimization problems under uncertainty depicted by Figure \ref{fig:sto_stage_strategies}.
It presents formulations: 
i) that explicitly consider the impact of uncertainty on cost, specifically \textit{risk-neutral} and \textit{risk-averse} versions of two-stage stochastic optimization problems, distributionally robust formulations, and robust min-max formulations;
ii) that focus on providing guarantees of constraint satisfaction despite different realizations of the uncertain parameters, including chance constraints, distributionally robust constraints, and robust constraints.
Various papers (refer to \cite{ROALD2023108725} for relevant references) employ a two-stage stochastic formulation for UC. The choice of the formulation depends on the uncertainty representation, which is often constrained by practical limitations on access to uncertainty data.

%%%%%%%%%%%%%%%%%%%
% Multi-stage UC
%%%%%%%%%%%%%%%%%%%%%%%
In multi-stage stochastic UC, the uncertainty is revealed in multiple steps, allowing us to update the wait-and-see decisions each time new information about the uncertainty is revealed while still accounting for the fact that we do not know what future realizations will look like.
A significant line of research on multi-stage UC has been the development of advanced decomposition algorithms to deal with tractability issues. The study \cite{8532315} reviews several papers mainly using two approaches: \textit{unit decomposition} and \textit{scenario decomposition}.
%
% Unit decomposition = spatial decomposition
In unit decomposition, constraints that couple generation units, such as load balancing, transmission, and spinning reserve constraints, are relaxed so that each sub-problem corresponds to a single generation unit.
% Multi-stage opti with Spatial decomposition technique applied to UC (Edf)
Such a decomposition scheme was first studied in \cite{496196}. In this paper, random disturbances are modeled as scenario trees, and optimization consists of minimizing the average generation cost over this tree-shaped future by implementing an augmented Lagrangian technique.
%
% Scenario decomposition
Alternatively, the scenario decomposition approach attempts to relax the coupling constraints among scenarios, usually named non-anticipativity constraints. The ensuing subproblems correspond to single scenarios (see \cite{8532315} for relevant references to this approach). 

% Limitation des approches de décomposition spatiales et par scénario
However, these solution methods usually apply only to relatively small scenario trees, and tractability issues are rapidly encountered in practical problems with large scenario trees. 
% SDDP 
To deal with large scenario trees, \cite{pereira1991multi} proposes the \textit{Stochastic Dual Dynamic Programming} (SDDP) approach. It is a sampling-based variant of nested Benders decomposition to solve multi-stage stochastic linear programs. It relies on approximating the expected-cost-to-go functions of stochastic dynamic programming by piecewise linear functions. The SDDP approach has been widely applied to multi-stage stochastic hydro-thermal scheduling problems.
The SDDP approach was investigated in a research internship \cite{lucille2023} at RTE R\&D in the context of probabilistic unit commitment. This work allowed us to formulate the problem and propose a first implementation on a simple case study. However, much work is still required to model complex case studies in a reasonable computation time.
%
% Extend SDDP to MPC  2018
SDDP has also been employed in \cite{KUMAR2018493} as a scalable approach to handling complex model predictive control (MPC) applications with uncertainties evolving over long time horizons and with fine time resolutions. 
%
% Extension de SDDP a des MILP: SDDiP
However, the SDDP approach cannot deal with integer variables in the expected cost-to-go. To overcome the intrinsic limitation of Benders-type decomposition for multi-stage stochastic integer programs, \cite{8532315} developed the \textit{Stochastic Dual Dynamic integer Programming} (SDDiP) algorithm.
It uses a new family of valid cuts, termed \textit{Lagrangian cuts}, which can achieve strong duality for mixed integer programs. SDDiP is also a sampling-based algorithm, as SDDP, and provides promising tractability properties in solving large-scale scenario trees.
Finally, a multi-stage stochastic approach is proposed for optimizing the day-ahead unit commitment of power plants and virtual power plants operating in the day-ahead and ancillary services markets by \cite{FUSCO2023120739}. It relies on a novel decomposition method where the main idea is to solve the problem as a sequence of two-stage stochastic programming models. For each two-stage problem, the scenarios are obtained employing a clustering algorithm, a modified k-medoids, choosing the most representative ones from the scenario tree.

\subsection{Contributions}

The main goal of this study is to design, implement, and test a probabilistic UC tool that handles the main uncertainties of system imbalance. Due to the technical constraints of conventional units, the problem is multi-stage by design.
%
% Choice of the formulation
However, we currently do not have enough data to generate a scenario tree for consumption, PV, wind power, and conventional unit failures. Thus, we cannot easily consider a multi-stage stochastic formulation. Indeed, it would require implementing an approach to simulate a relevant scenario tree from available data, which is not straightforward. In addition, due to the technical constraints of the conventional units, the multi-stage stochastic formulation becomes highly complex and intractable when increasing the number of units considered. An SDDP algorithm would be required to tackle this issue. However, as demonstrated in \cite{lucille2023}, this approach is challenging to implement and to obtain interpretable and relevant results over realistic case studies. 
% Approximation proposed
Thus, in the vein of \cite{ALLAWATI2021116882}, this work proposes an approximation of the multi-stage problem with a \textit{two-stage stochastic model predictive control}. It consists of a sequence of two-stage stochastic optimization models where the results from each model feed into each subsequent model, allowing for scenarios to be updated as more information becomes available.
Then, we implement an approximation of this two-stage stochastic MPC with a single two-stage problem.

Overall, the contributions of this study are threefold.
\begin{enumerate}
    \item A formulation of the multi-stage stochastic problem is proposed. The stochastic optimizer computes an optimal dispatch by considering the fixed and variable costs, the technical constraints of power plants, the costs of lost load and production, and the uncertainty of the system imbalance modeled with consumption, PV, and wind power scenarios.
    \item A tractable approximation of the multi-stage stochastic problem is presented with a \textit{two-stage stochastic model predictive control} using a risk-neutral approach with a scenario-based formulation for each two-stage problem.
    \item Finally, an approximated version with a single two-stage stochastic problem is implemented. It is tested on a real case study where the French TSO RTE faced a deficit of downward available margins.
\end{enumerate}

% Structure du papier 
Section \ref{sec:pb-statement} presents the abstract formulation of the probabilistic unit commitment problem considered. 
Section \ref{sec:sto-framework} describes the stochastic framework adopted, and Section \ref{sec:formulations} provides the formulations of the heuristics considered to solve the problem. 
Section \ref{sec:results} presents the results of a case study where the French TSO RTE faced a deficit of downward available margins.
Finally, Section \ref{sec:conclusions} presents the conclusions and future works.
Appendix \ref{appendix:scenario-generation} provides the methodology for generating the PV, wind power, and consumption scenarios.
Appendix \ref{appendix:comparison-with-without-scenario-selection} provides additional results related to section \ref{sec:results}.

\section{Problem statement}\label{sec:pb-statement}

This section presents the abstract formulation of the probabilistic unit commitment problem considered. We adopt a central dispatch approach, and the conventional units' technical constraints make the problem multi-stage by design. Thus, a multi-stage formulation and an approximation with a sequence of two-stage stochastic problems are proposed. Finally, we adopt a risk-neutral approach with a scenario-based formulation where the uncertainties related to consumption and renewable production are modeled with a finite set of scenarios.

\subsection{Central dispatch}

We adopt a \textit{benevolent monopoly} approach where the TSO computes the unit commitment of all units with a central dispatch over several hours. Each day $D$ comprises $T$ periods of duration $\Delta_t$. For instance, $T=24$ with $\Delta_t$ = 1 hour.
% Description des éléments qui composent le central dispatch
The central dispatch comprises: i) conventional power plants composed of nuclear, fossil-based, and hydraulic power plants; ii) renewable production (PV and wind power); iii) consumption. 
Each conventional production unit has technical constraints, such as startup and shutdown delays, minimum flat duration after a power variation, ramping power constraints, minimum time ON/OFF, and minimal and maximal power. These technical constraints make the formulation of the problem multi-stage by design.
In this study, we assume consumption and renewable units are not flexible and result in fatal residual demand that the production of conventional units must meet. Further developments should consider flexible demand, such as electric vehicles and renewable curtailment.

\begin{figure}[tb]
\centerline{\includegraphics[width=90mm]{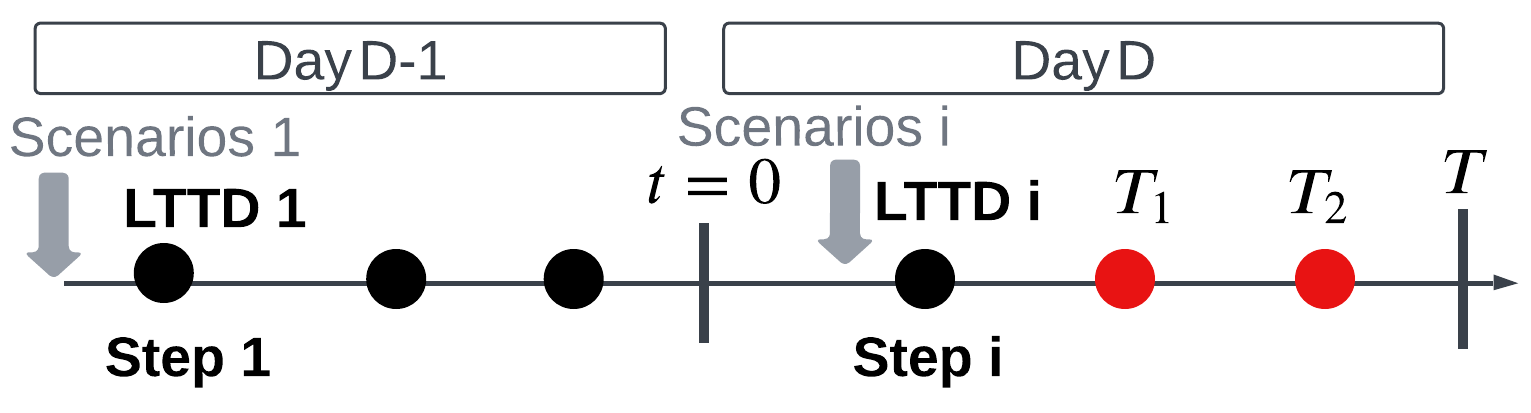}}
\caption{Illustration of the problem timeline with $M$ LTTDs for a dispatch over $[T_1, T_2]$. For instance, the first LTTD at step 1 corresponds to the decision to start or not nuclear power plants based on the uncertainty available at this stage.}
\label{fig:LTTD-time-line}
\end{figure}

% Explication des LTTD
The technical constraints of conventional units imply \textit{last time to decides} (LTTDs) corresponding to the last period where a decision can be taken, such as starting/shutting down a unit or making an upward/downward power variation. For instance, a nuclear power plant has a minimum starting delay of 10 hours, a minimal ON/OFF duration of 24 hours, and a flat duration of two hours between two power variations.
Thus, the central dispatch over $[T_1, T_2]$ is a multi-step stochastic program with $M$ stages corresponding to $M$ LTTDs, as Figure \ref{fig:LTTD-time-line} depicts. The uncertainty results from renewable production and consumption, described as random variables.
At each LTTD, the uncertainty evolves, and conventional units can encounter technical failures.

\subsection{Multi-stage formulation}\label{sec:multi-stage-formulation}

\begin{figure}[tb]
	\begin{subfigure}{.5\textwidth}
		\centering
		\includegraphics[width=\linewidth]{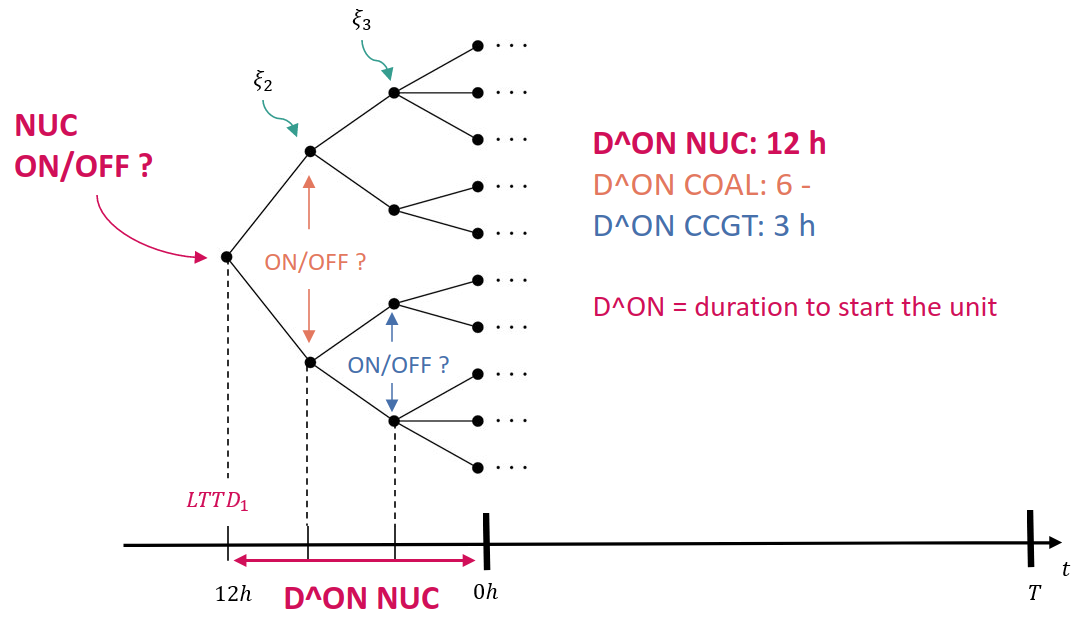}
		\caption{Problem 1: nuclear ON/OFF.}
	\end{subfigure}
 	\begin{subfigure}{.5\textwidth}
		\centering
		\includegraphics[width=\linewidth]{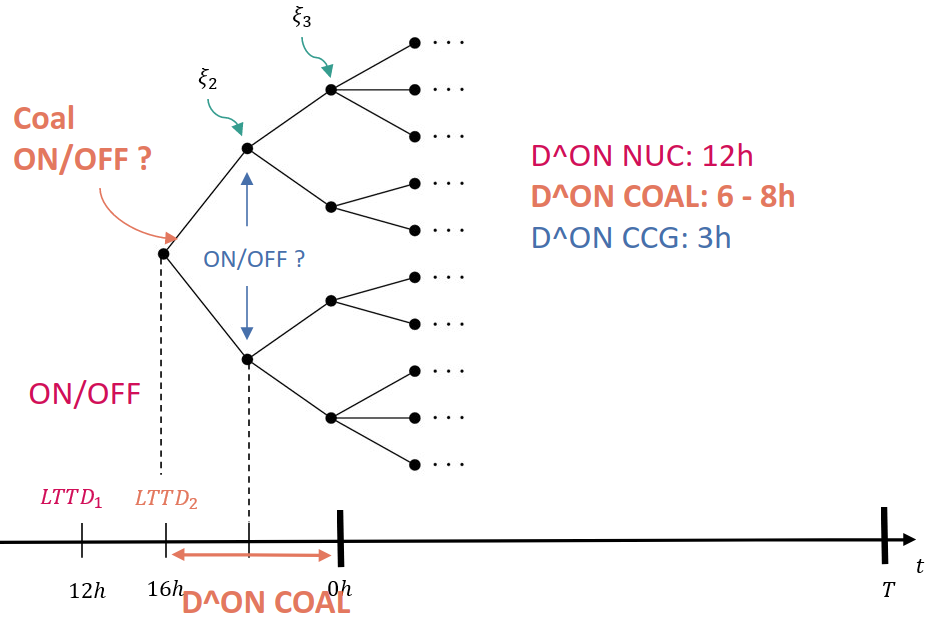}
		\caption{Problem 2: coal ON/OFF.}
	\end{subfigure}
\caption{Illustration of the sequence of multi-stage problems on a simplified example. In this example, we consider only three types of units: nuclear, coal, and CCGT (combined cycle gas turbine) power plants. $D^\text{ON}$ is the minimal duration to start a power plant. Notice that the values provided are examples that are not necessarily realistic.
In the first multi-stage problem, the first-stage variables are nuclear units' ON/OFF status. In the second multi-stage problem, the first-stage variables are coal power plants' ON/OFF status. Then, in a third problem, it would be the ON/OFF of the CCGT power plants. This Figure is adapted from \cite{lucille2023}.}
	\label{fig:multi-stage}
\end{figure}

In the multi-stage setting, the sequence $\boldsymbol \xi_{t \in [0, T]}$, of data vectors is modeled as a stochastic process, \textit{i.e.}, a sequence of random variables with a specified probability distribution. 
$\boldsymbol \xi_{[t, t']}$ is a sequence of random data vectors corresponding to stages $t$ through $t'$, and $\xi_{[t, t']}$ a specific realization of this sequence of random vectors. 
In this setting, the values of the decision vector $\boldsymbol x_t$, chosen at stage $t$, may depend on the information $\boldsymbol \xi_{[1, t]}$ available up to time $t$, but not on the results of future observations. This is the essential requirement of \textit{non-anticipativity}. As $\boldsymbol x_t$ may depend on $\boldsymbol \xi_{[1, t]}$, the sequence of decisions is also a stochastic process. 
In a generic form, a linear $T$-stage stochastic programming problem can be written in the nested formulation  \cite{birge2011introduction,zou2019stochastic,RUSZCZYNSKI20031} as follows
\begin{align}
\label{eq:multi-stage-general}
\min_{\boldsymbol x_1 \in F_1} \bigg \{ & c_1^\intercal \boldsymbol x_1   +\mathbb{E}_{\boldsymbol \xi_{[2, T]}|\xi_{[1, 1]}}\bigg[ \min_{\boldsymbol x_2 \in F_2(\boldsymbol x_1, \xi_2)} \bigg \{ c_2^\intercal \boldsymbol x_2(\xi_2) + \cdots \notag \\  
& + \mathbb{E}_{\boldsymbol \xi_{[T, T]}|\xi_{[1, T-1]}}\bigg[ \min_{\boldsymbol x_T \in F_T(\boldsymbol x_{T-1}, \xi_T)}  \bigg \{ c_T^\intercal \boldsymbol x_T(\xi_T)   \bigg \} \bigg] \bigg \}\bigg] \bigg \},
\end{align}
where $F_T(\boldsymbol x_{t-1}, \xi_t)$ is the feasible set of the stage $t$ problem, which depends on the decision in stage $t-1$ and the information $\xi_t$ available in stage $t$.
$\mathbb{E}_{\boldsymbol \xi_{[t, T]}|\xi_{[1, t-1]}}$ is the expectation in stage $t$ with respect to the
conditional distribution of $\boldsymbol \xi_{[t, T]}$ given realization $ \xi_{[1, t-1]}$ in stage $t-1$. 
The first-stage decisions are $\boldsymbol x_1$, referred to as \textit{here-and-now} decisions, which must be decided before the sequence of random data vectors $\boldsymbol \xi_{[2, T]}$ is known. The multi-stage decision variables $\boldsymbol x_2$, ..., $\boldsymbol x_T$, referred to as \textit{wait-and-see} decisions, are taken in response to the realization of the uncertain random data vectors $\boldsymbol \xi_{[2, T]}$.
Figure \ref{fig:multi-stage} depicts the sequence of multi-stage problems in a simplified example.

\begin{figure}[tb]
    \centering
    \includegraphics[width=\linewidth]{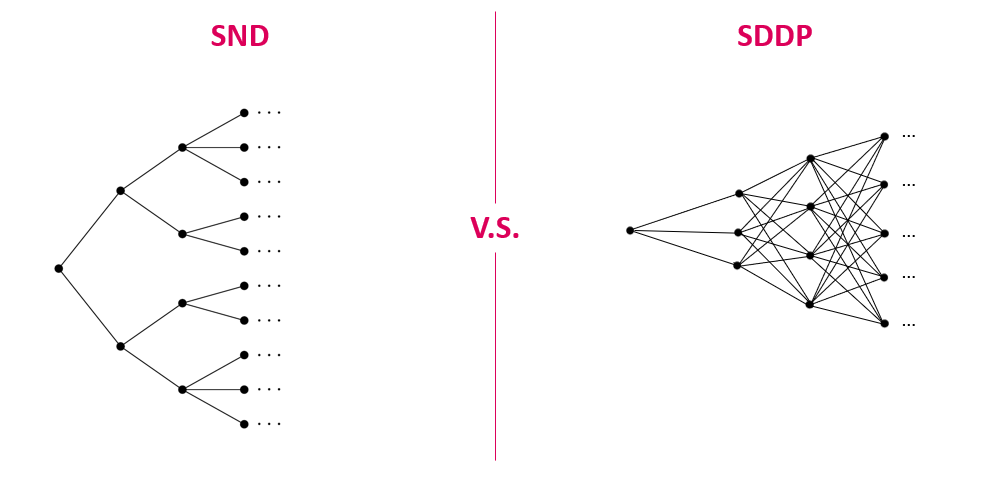}
    \caption{Scenario tree structure for multi-step algorithms Stochastic Nested Decomposition (SND) \cite{sddip} and Stochastic Dual Dynamic Programming \cite{pereira1991multi}, which exploits the Markovian assumption. This Figure is adapted from \cite{lucille2023}.}
    \label{fig:snd-vs-sddp}
\end{figure}

% Markov process
Notice that the random process $(\boldsymbol \xi_2, \ldots, \boldsymbol \xi_T)$ is said to be Markovian, if for each $t =2, \ldots, T-1$ the conditional distribution of $\boldsymbol \xi_{t+1}$ given $\boldsymbol \xi_{[1,t]}$ is the same as the conditional distribution of $\boldsymbol \xi_{t+1}$ given $\boldsymbol \xi_t$. In this case, the stochastic model is simplified considerably. 
Figure \ref{fig:snd-vs-sddp} presents the difference between a classic tree scenario (with no Markovian assumption) and one assuming this framework.

% Scenario approach and scenario tree
In addition, if the number of scenarios is finite, the linear multi-stage stochastic program is formulated as one large (deterministic) linear programming problem. Thus, computational approaches for this problem are based on approximating the stochastic process $(\boldsymbol \xi_2, \ldots, \boldsymbol \xi_T)$ by a process having finitely many realizations in the form of a scenario tree \cite{zou2019stochastic,RUSZCZYNSKI20031}. Monte Carlo methods may construct an approximation as in the sample average approximation approach.

\subsection{Two-stage problems}

\begin{figure}[tb]
\centering
\includegraphics[width=\linewidth]{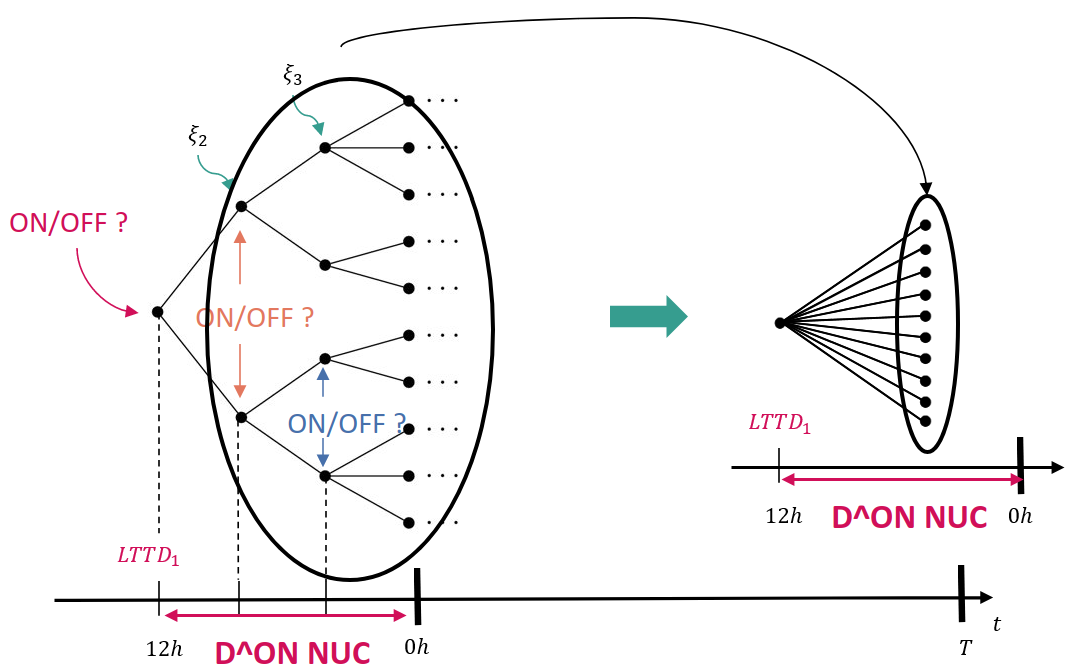}
\caption{Approximation with a two-stage problem on the simplified example of Figure \ref{fig:multi-stage}. This Figure is adapted from \cite{lucille2023}.}
\label{fig:two-stage-approximation}
\end{figure}

In many applications, the sequence of multi-stage problems is too complex to solve because it leads to intractable problems. Thus, a standard approximation is formulating the problem with a sequence of two-stage problems \cite{ROALD2023108725}.
Figure \ref{fig:two-stage-approximation} depicts this approximation in the previous simplified example.
In the two-stage formulation, the multi-stage formulation (\ref{eq:multi-stage-general}) becomes
\begin{subequations}
\begin{align}
        \min_{\boldsymbol x,\boldsymbol y_{\xi}} &f^F(\boldsymbol x) + \mathbb{E}[f^S(\boldsymbol x,\boldsymbol y_{\xi},\boldsymbol \xi)] \\
        \text{s.t. } &h^F(\boldsymbol x) = 0 \\
        &g^F(\boldsymbol x) \leq 0 \\
        &h^S(\boldsymbol x,\boldsymbol y_{\xi},\boldsymbol \xi) = 0 \\
        &g^S(\boldsymbol x,\boldsymbol y_{\xi},\boldsymbol \xi) \leq 0,
\end{align}
\label{eq:two-stage}
\end{subequations}
where $f^F$ and $f^S$ are the first-stage and second-stage costs, $\boldsymbol x$ represents the \textit{first-stage variables}, $\boldsymbol \xi$ the random variables modeling the uncertain parameters, and $\boldsymbol y_{\xi}$ the \textit{second-stage variables} depending on the realization of $\boldsymbol \xi$.
This formulation is \textit{risk-neutral}, as it similarly treats costs above and below the expected value. It is only sometimes desirable, as some unfavorable outcomes may have a disproportionate impact. Problem formulations that focus specifically on limiting the negative impacts of uncertainty realizations are called risk-averse formulations \cite{ROALD2023108725} with variance minimization, for example.

\subsection{Risk-neutral approach}

The study \cite{ROALD2023108725} provides an overview of existing methods for modeling and optimizing uncertainty-related problems. 
We adopt a \textit{scenario-based formulation} \cite{birge2011introduction} for each two-stage stochastic problem as probabilistic forecasters provide scenarios that represent possible realizations of PV, wind power, and consumption. The second-stage variables and constraints are indexed for every timestep of every possible scenario. 
Then, we consider a risk-neutral approach by minimizing the total expected cost, as it equally treats costs above and below the expected value.
These choices are motivated by the complexity of the multi-stage stochastic problem considered. The numerical resolution must be tractable to be computed within a few minutes. In addition, it is attractive to implement a benchmark before investigating more complex formulations, such as risk-averse formulations, by considering the Conditional Value-at-Risk \cite{ROALD2023108725}. 

\section{Stochastic programming framework}\label{sec:sto-framework}

This study considers the sequence of two-stage stochastic problems to approximate the multi-stage formulation presented in Section \ref{sec:multi-stage-formulation}. 
First, the units and their technical constraints considered are presented. These constraints imply LTTDs used to model the sequence of two-stage problems. 
We first present a single-phase approach that provides a baseline. Then, the multi-phase sequential framework is presented.

\subsection{Description of the units considered}

% Description du multi-stage stochastic program
Due to technical constraints, the number of LTTDs increases with the production units. This makes the problem highly complex and intractable. 
Thus, the problem at hand is modeled with a limited number of stages by grouping the technical constraints of similar conventional units, resulting in fewer LTTDs and stages.

\begin{figure}[tb]
\centerline{\includegraphics[width=90mm]{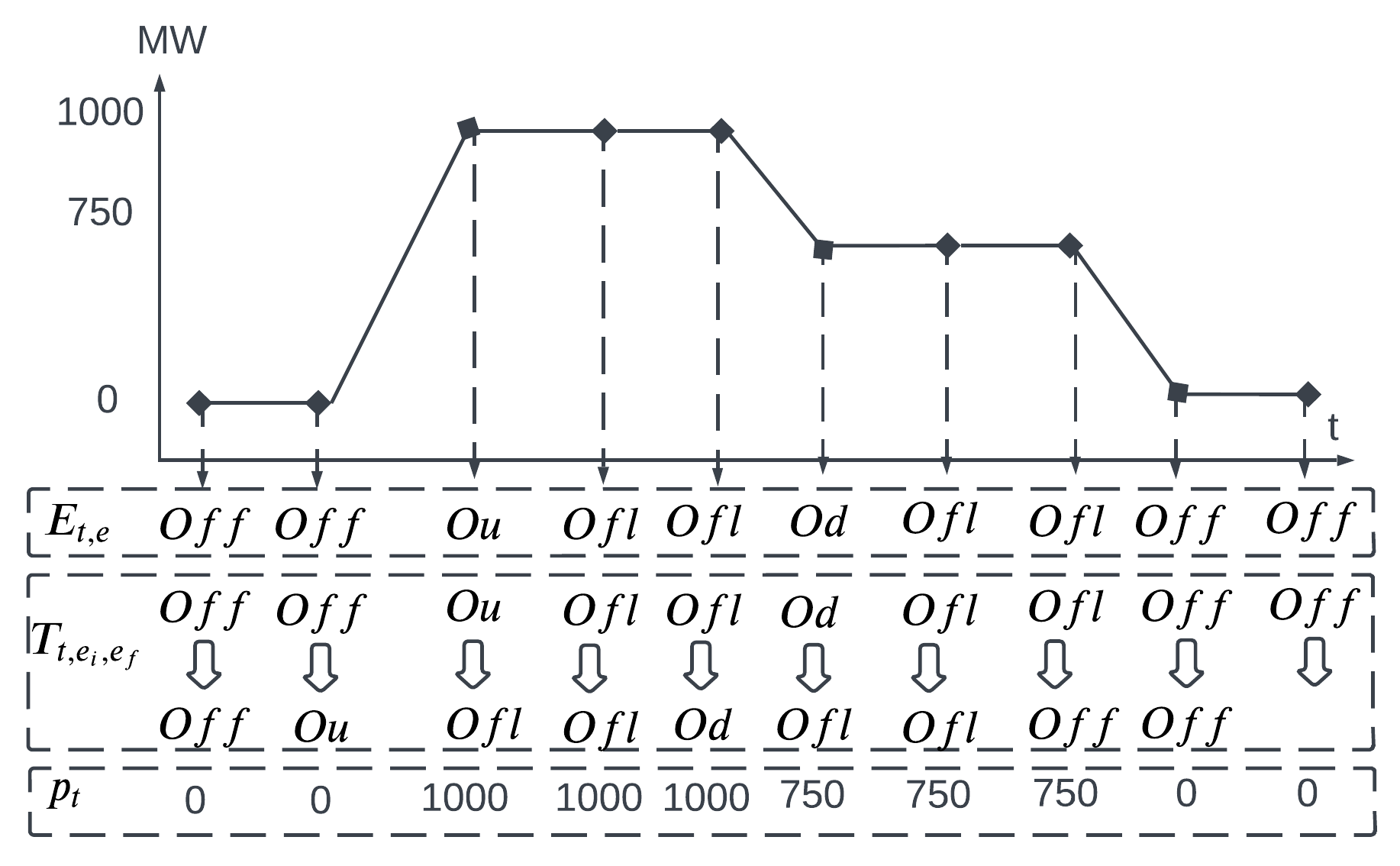}}
\caption{Illustration of the state transitions of a thermal unit.}
\label{fig:state-description}
\end{figure}

% four groups of units
We group the power plants into four categories. The nuclear, coal, CCGT, and OCGT (open-cycle gas turbine). In this study, we include hydro-power plants in the group of OCGT units due to similar time delays for the technical constraints. Notice that the specific case of hydro-power plants is pointed out as one of the other research directions to address in future works. 
The units of each category are assumed to have the same technical constraints and, in particular, the same minimum start-up duration and, thus, the same LTTD.
Then, each unit is described by four states $e \in \Omega_\text{E}$ with: 
\begin{itemize}
    \item $e = $ OU, the unit is ON at $t$ and made an upward variation from $t-1$ to $t$;
    \item $e = $ OD, the unit is ON at $t$ and made a downward variation from $t-1$ to $t$;
    \item $e = $ OFL, the unit is ON at $t$, and its power was constant from $t-1$ to $t$;
    \item $e = $ OFF, the unit is OFF at $t$.
\end{itemize}
Table \ref{tab:transitions-allowed} provides the transitions between states allowed and forbidden, and Figure \ref{fig:state-description} depicts an example of transitions with $E_{t,e}$ as the state $e$ of this unit at time $t$, and $T_{t,e_i, e_f}$  the transition at $t$ from the initial state $e_i$ to the final state $e_f$ at $t+1$.
\begin{table}[tb]
\renewcommand{\arraystretch}{1.25}
\begin{center}
\begin{tabular}{lllll}
\hline \hline
     & OU & OD & OFL & OFF \\
OU &  \checkmark   &   $\times$   & \checkmark & $\times$  \\
OD &  $\times$   &   \checkmark   & \checkmark & \checkmark  \\
OFL &  \checkmark   &   \checkmark  & \checkmark & \checkmark  \\
OFF &  \checkmark   &   $\times$   & $\times$ & \checkmark  \\
\hline \hline
\end{tabular}
\caption{Transitions allowed and forbidden between states.}
\label{tab:transitions-allowed}
\end{center}
\end{table}

\subsection{Technical constraints and sequence of LTTDs}

\begin{figure}[tb]
\begin{subfigure}{0.5\textwidth}
		\centering
		\includegraphics[width=\linewidth]{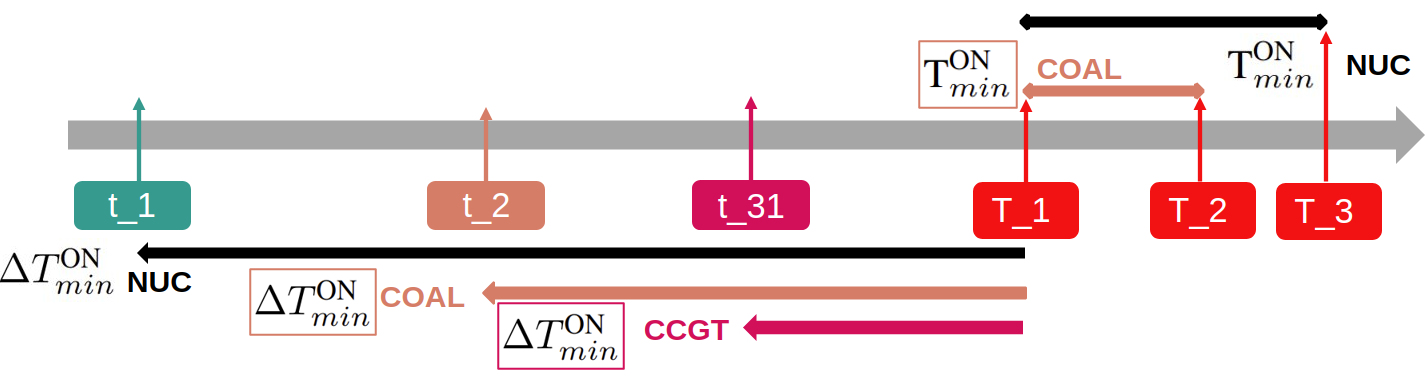}
  \caption{Sequence of the first three LTTDs.}
\label{fig:LTTD-sequence}
	\end{subfigure}
 	\begin{subfigure}{0.5\textwidth}
		\centering
		\includegraphics[width=\linewidth]{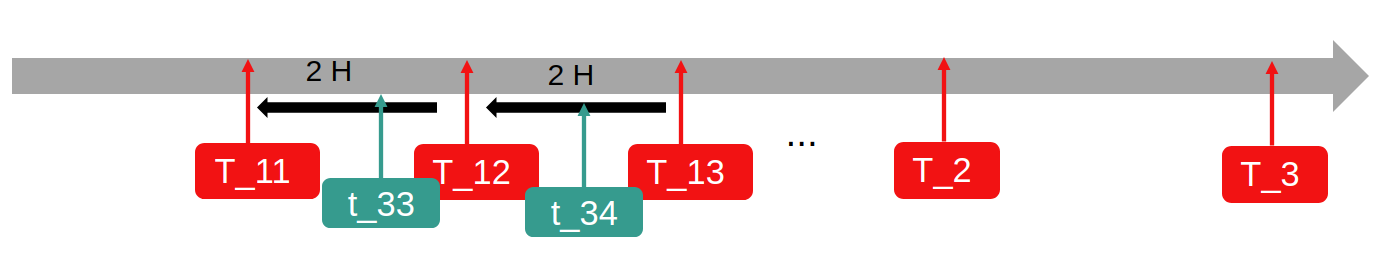}
  \caption{Decomposition of the period of study by 2 hours step.}
    \label{fig:study-period}
	\end{subfigure}
	\caption{Illustration of the timeline with an optimization period $[T_1, T_3]$ and the period of study $[T_1, T_2]$. The period of study is decomposed in a sequence of 2 hours period $[T_1, T_2]= [T_{11}, T_{12}] \cup [T_{12}, T_{13}]\cup \ldots [T_{1n}, T_2] $.}
	\label{fig:time-line-framework}
\end{figure}

% Assumptions
Table \ref{tab:technical-constraints} provides the values of the technical constraints adopted for each group of units.
\begin{table}[tb]
\renewcommand{\arraystretch}{1.25}
\begin{center}
\begin{tabular}{lrrrrrrr}
\hline \hline
Unit & $\Delta T^\text{ON}_{min}$ & $\Delta T^\text{OFF}_{min}$ & $\text{T}_{min}^\text{ON}$ & $\text{T}_{min}^\text{OFF}$ & $\text{T}^\text{FLAT}$ & $\text{T}^\text{ON}_{max}$ & $\text{N}^\text{ON}_{max}$\\
\hline
NUC & 600 & 60 & 1440 & 1440 & 90 & - & - \\
COAL & 480 & 15 & 480 & 480 & 15 & - & - \\
CCGT & 180 & 15 & 180 & 120 & 15 & - & - \\
OCGT & 15 & 15 & 60 & 30 & 15 & 480 & 2 \\
\hline \hline
\end{tabular}
\caption{Values (in minutes) of technical constraints considered for thermal units. The maximum and minimal power duration is $\nabla P^\text{max}$ and is less than 1 hour for all units considered in this study. $\Delta T^\text{ON}_{min}$ is the minimal duration to start the unit. $\Delta T^\text{OFF}_{min}$ is the minimal duration to shut-down the unit. $\text{T}_{min}^\text{ON}$ and $\text{T}_{min}^\text{OFF}$ are the minimal duration when started and when shut-down. $\text{T}^\text{FLAT}$ is the minimal duration between two power variation. $\text{T}^\text{ON}_{max}$ is the maximum ON duration and $\text{N}^\text{ON}_{max}$ is the maximal number of startup within a day.}
\label{tab:technical-constraints}
\end{center}
\end{table}
% Sequence of LTTDs
Due to the minimum start-up duration of the units, we consider three LTTDs related to the start-up of the units: the LTTDs of the nuclear, coal, and CCGT units. Figure \ref{fig:LTTD-sequence} depicts the sequence of the first three LTTDs considered. The decisions to start nuclear, coal, and CCGT units are taken at $t = t_1$, $t = t_2$, and $t = t_{31}$. 
The optimization period is the interval $[T_1, T_3]$, and the period of study is $[T_1, T_2]$. It is decomposed in a sequence of 2 hours period $[T_1, T_2]= [T_{11}, T_{12}] \cup [T_{12}, T_{13}] \cup \ldots [T_{1n}, T_2] $ as depicted by Figure \ref{fig:study-period}. Indeed, it corresponds to the duration of the balancing window where the French TSO can activate units. However, the duration of this window can be changed, and the formulation is not specific to a value of two hours.

\subsection{Single-phase framework}

In the single-phase framework, only one optimization problem is solved, approximating the sequence of two-stage problems. It serves as a baseline, as considered in \cite{ALLAWATI2021116882}, to compare the results with a more complex formulation involving a sequence of two-stage problems in the multi-phase framework.

\begin{figure}[tb]
\includegraphics[width=\linewidth]{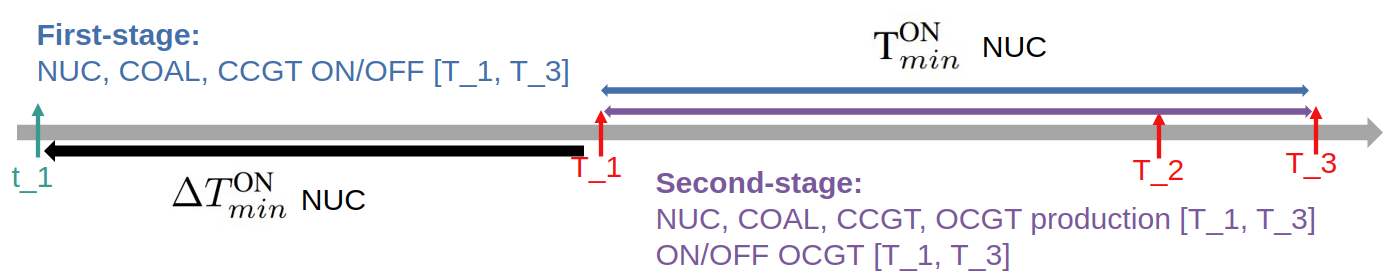}
\caption{First-stage and second-stage variables considered in the single framework. 
The first-stage variables are the ON/OFF status over $[T_1, T_3]$ of nuclear, coal, and CCGT units, represented in blue.
The second-stage variables comprise the ON/OFF status of OCGT units and the production level of all units over $[T_1, T_3]$, represented in purple.
}
\label{fig:single-framework}
\end{figure}

Figure \ref{fig:single-framework} depicts the single framework timeline.
In this formulation, the first-stage variables are the ON/OFF status for the entire optimization period $[T_1, T_3]$ of all units with $\Delta T^\text{ON}_{min} \geq $ 60 minutes, \textit{e.g.}, nuclear, coal, and CCGT units. 
% Approximation motivations
Notice that this approximation is conducted because the nuclear $\text{T}_{min}^\text{ON}$ is approximately 24 hours. Thus, the ON/OFF status variables are decided once and for all over this duration, and $[T_1, T_3]$ is supposed to equal 24 hours in this study. However, this approximation for coal and CCGT units is more questionable as their $\text{T}_{min}^\text{ON}$ is smaller than 24 hours, in particular, CCGT units.
The second-stage variables comprise the ON/OFF status of OCGT units and the production level of all units over $[T_1, T_3]$.
Therefore, resolving this problem provides the ON/OFF status of nuclear, coal, and CCGT units over the period of study $[T_1, T_2]$, which will be used in the evaluation.

\subsection{Multi-phase sequential framework}

In the multi-phase framework, a sequence of two-stage optimization problems is considered to model the LTTDs of the nuclear, coal, and CCGT units.
\begin{figure}[tb]
\centerline{\includegraphics[width=\linewidth]{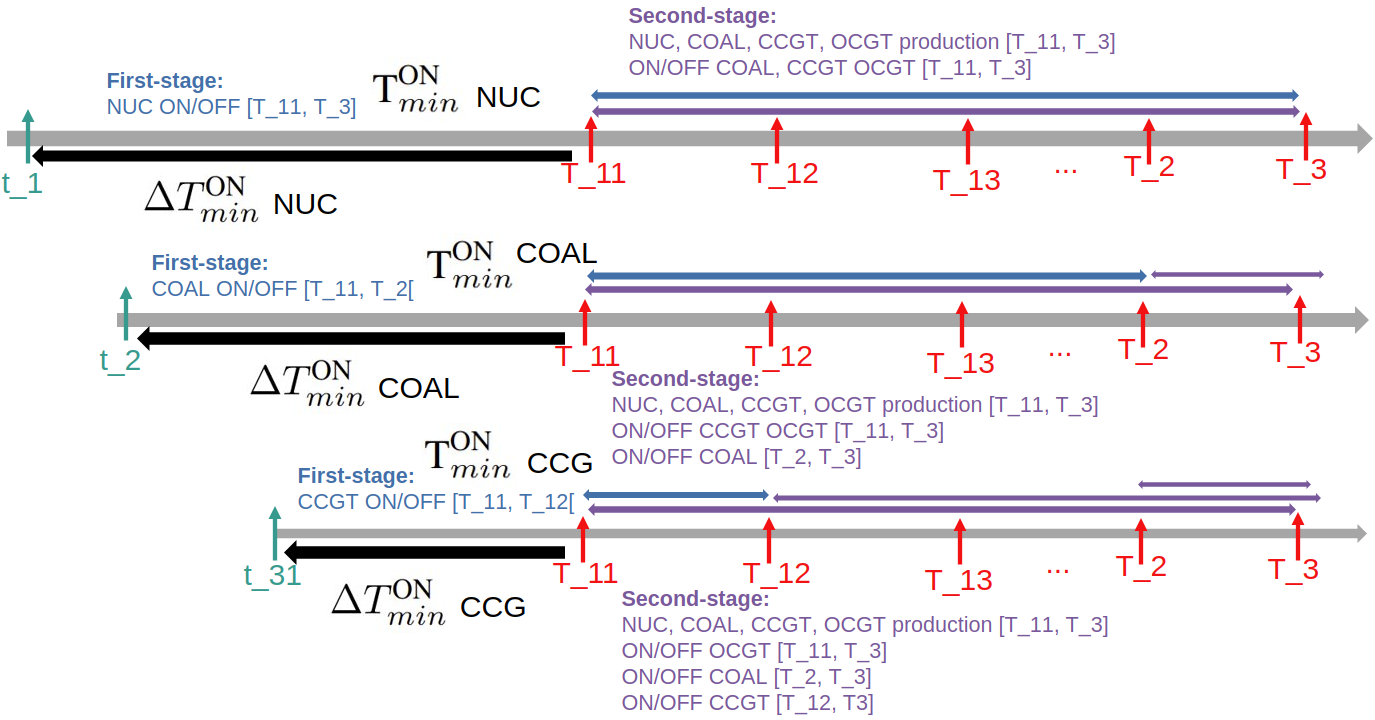}}
\caption{Illustration of the timeline of the multi-phase framework. The blue lines indicate the first-stage variables, and the purple lines are the second-stage variables. For instance, in the second two-stage problem, the first-stage variables are the ON/OFF status of the coal units over $[T_1, T_2]$, and the second-stage variables comprise the ON/OFF status of CCGT and OCGT for $[T_1, T_3]$, the level of production of all units for $[T_1, T_3]$, and the OFF/OFF status of the coal units for $[T_2, T_3]$.}
\label{fig:multi-phase-framework}
\end{figure}
Figure \ref{fig:multi-phase-framework} depicts the multi-phase framework timeline.
In this formulation, three two-stage stochastic problems are considered to compute sequentially the ON/OFF status of nuclear, coal, and CCGT units for the first period of 2 hours $[T_{11}, T_{12}]$.

% NUC LTTD
The first two-stage optimization, at $t=t_1$ corresponding to the LTTD of the nuclear units, provides the ON/OFF status of nuclear units over $[T_1, T_3]$. 
The first-stage variables are the ON/OFF status of the nuclear units over $[T_1, T_3]$.  Indeed, the minimum ON duration of nuclear units $\text{T}_{min}^\text{ON}$ is 1440 minutes, corresponding to this period in this study. 
The second-stage variables comprise the ON/OFF status of COAL, CCGT, and OCGT and the production level of all units over $[T_1, T_3]$.

% COAL LTTD
The second optimization, at $t=t_2$ corresponding to the LTTD of the coal units, provides the ON/OFF status of coal units over $[T_1, T_2[$. 
The first-stage variables are the ON/OFF status of the coal units over $[T_1, T_2[$. Indeed, the minimum ON duration $\text{T}_{min}^\text{ON}$ of coal units is 480 minutes, corresponding to this period. 
The second-stage variables comprise the ON/OFF status of CCGT and OCGT and the production level of all units over $[T_1, T_3]$. It also comprises the ON/OFF status of coal units over $[T_2, T_3]$.

% CCGT LTTD
Finally, the third optimization problem, at $t=t_{31}$ corresponding to the LTTD of the CCGT units, provides the ON/OFF status of CCGT units over $[T_1, T_{12}[$. 
The first-stage variables are the ON/OFF status of the CCGT units over $[T_1, T_{12}[$. Indeed, the minimum ON duration $\text{T}_{min}^\text{ON}$ of CCGT units is 180 minutes, corresponding approximately to this period. 
The second-stage variables comprise the ON/OFF status of OCGT and the level of production of all units over $[T_1, T_3]$. It also comprises the ON/OFF status of coal units over $[T_2, T_3]$ and the ON/OFF status of CCGT units for $[T_{12}, T_3]$.

\begin{figure}[tb]
\includegraphics[width=\linewidth]{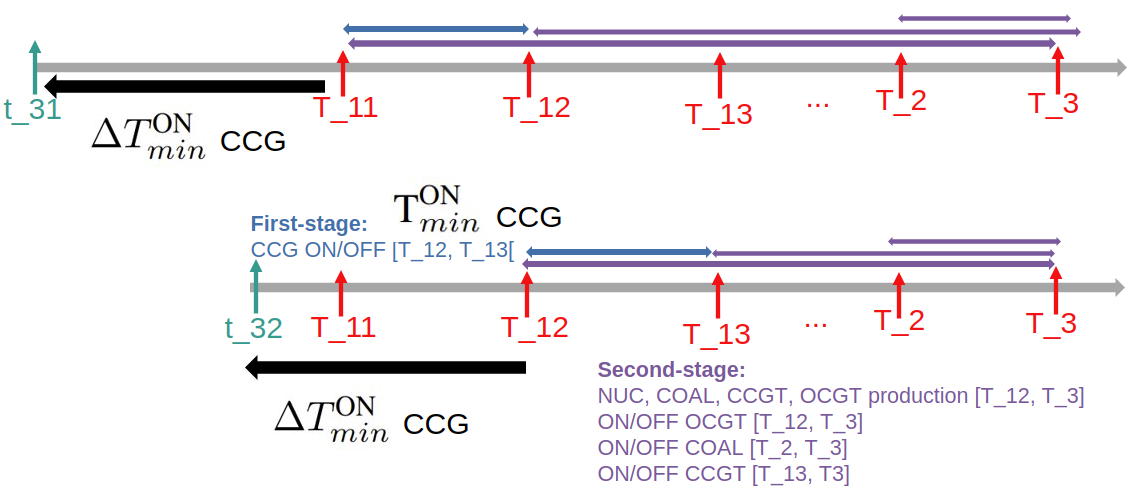}
\caption{Two-stage optimization problem at $t=t_{32}$. The first-stage variables are the ON/OFF status of the CCGT units over $[T_{12}, T_{13}[$. The second-stage variables comprise the ON/OFF status of OCGT and OCGT for $[T_{12}, T_3]$, the level of production of all units for $[T_{12}, T_3]$, the OFF/OFF status of the CCGT units for $[T_{13}, T_3]$, and the OFF/OFF status of the coal units for $[T_2, T_3]$.}
\label{fig:multi-phase-framework-2}
\end{figure}

% result of the sequence
Then, to compute the ON/OFF status for the rest of the period of study, we move forward in time as depicted by Figure \ref{fig:multi-phase-framework-2} for the optimization problem at $t=t_{32}$ which provides the ON/OFF status of nuclear, coal, and CCGT units over $[T_{12}, T_{13}[$.
This process is conducted until the last period $[T_{3n}, T_2]$.
Overall, resolving the sequence of two-stage stochastic problem sequentially allows us to compute the ON/OFF status of the nuclear, coal, and CCGT units over the period of study $[T_1, T_2]$.

\section{Formulations of the proposed approaches}\label{sec:formulations}

This section presents the mathematical formulations of the single-phase and multi-phase sequential approaches.

\subsection{Single-phase framework formulation}\label{sec:sto-simple}

% OBJ 
The mixed-integer linear programming optimization problem modeling the two-stage stochastic problem in the single framework is 
\begin{subequations}
\label{eq:single-obj}	
\begin{align}
\min \  &  J_f^1 + J_f^2 + J_v , \\ 
&  \text{s.t. } \notag \\
& \quad (\ref{eq:lttd-1-constraints-part-1-balance}) - (\ref{eq:lttd-1-constraints-part-6}).
\end{align}
\end{subequations}
Where $J_1^1$ models the fixed cost of nuclear, coal, and CCGT units, $J_f^2$ the fixed cost of OCGT units, and $J_v$ the variable cost of all thermal units (nuclear, fossil-based), the lost load, and the lost production. They are defined as follows
\begin{subequations}
\label{eq:single-obj-details}	
\begin{align}
% first-stage cost
J_f^1  & = \sum_{u \in \Omega_\text{FT}}  \pi_f^u  \bigg[\sum_{t\in \mathcal{T}}  T_{t,\text{OFF},\text{OU}}^u +   T_{\text{OFF},\text{OU}}^{0,u} \bigg] , \\ 
% second-stage fixed cost
J_f^2  & =  \sum_{u \in \Omega_\text{ST}}  \sum_{s \in \Omega_S} p_s \pi_f^u\bigg[\sum_{t\in \mathcal{T}}   T_{t,s,\text{OFF},\text{OU}}^u +   T_{s,\text{OFF},\text{OU}}^{0,u} \bigg],\\
% second-stage variable cost
J_v  & =  \sum_{t\in \mathcal{T}} \sum_{s \in \Omega_S} p_s \bigg[ \sum_{u \in \Omega_\text{T}} \pi_v^u p_{t,s}^u +  \pi^\text{dns} \text{dns}_{t,s}  + \pi^\text{spill} \text{spill}_{t,s} \bigg],
\end{align}
\end{subequations}
where $\mathcal{T} = [T_1, T_3]$ is the set of periods, $\Omega_\text{FT}$ is the set of first-stage units (nuclear, coal, and CCGT in the single-phase framework), $\Omega_\text{ST}$ is the set of second-stage units (OCGT in the single-phase framework), and $\Omega_\text{T}= \Omega_\text{FT} \cup \Omega_\text{ST}$ the set of all thermal units.

% scenarios
The scenarios of consumption and renewable generations of this problem issued at $t_1$ are $\hat{d}_{t,s}^{t^1}$ and $\hat{p}_{t,s}^{t^1, u \in \Omega_\text{Nr}}$ with $t \in \mathcal{T}$.
%
% first-stage
The first-stage variables are $T_{t,e_i,e_f}^{u \in \Omega_\text{FT}}$, $T_{e_i,e_f}^{0,u \in \Omega_\text{FT}}$, and $E_{t,e}^{u \in \Omega_\text{FT}}$ with $ t \in \mathcal{T}$. They are binary variables that allow considering the fixed start-up costs $\pi_f^u$ of the first-stage units. $T_{e_i,e_f}^{0,u}$ allows to model a transition from state $e_i$ at $t<0$ to state $e_f$ at $t=0$. 
%
% second-stage
The second-stage variables are $T_{t,s,e_i,e_f}^{u \in \Omega_\text{ST}}$, $T_{s,e_i,e_f}^{0,u \in \Omega_\text{ST}}$, $E_{t,s,e}^{u \in \Omega_\text{ST}}$, the power of all thermal units $p_{t,s}^{u \in \Omega_\text{T}}$, the lost load $\text{dns}_{t,s}$, and the lost production $\text{spill}_{t,s}$ with $ t \in \mathcal{T}$. 
Overall, resolving this problem provides the ON/OFF status of the nuclear, coal, and CCGT units, $E_{t \in \mathcal{T}, e=Off}^{u \in \Omega_\text{FT}}$.

The constraints of the two-stage stochastic problem in the single-phase framework are valid $\forall t \in \mathcal{T}$ and $\forall s \in \Omega_S$ (when variables are indexed by scenario $s$). They can be decomposed into six parts:
\begin{enumerate}
    \item constraints related to the power balance (\ref{eq:lttd-1-constraints-part-1-balance}) where $\tilde{d}_{t,s}^{t^1}$ is the residual demand (demand minus the total of renewable generation), the thermal unit's maximum and minimum power (\ref{eq:lttd-1-constraints-part-1-pmax})-(\ref{eq:lttd-1-constraints-part-1-pmin}), and instantaneous power variation (\ref{eq:lttd-1-constraints-part-1-gradient-infinity-nuc})-(\ref{eq:lttd-1-constraints-part-1-gradient-infinity-fossil});
    \item constraints defining the state transitions $E_{t,s,e}^u$, $T_{t,s,e_i,e_f}^u$, and $T_{s,e_i,e_f}^{0,u}$ variables (\ref{eq:lttd-1-constraints-part-2-E-def-1})-(\ref{eq:lttd-1-constraints-part-2-E-def-2}), and ensuring there is only one state $e$ per timestep (\ref{eq:lttd-1-constraints-part-2-E-uniq-1}-\ref{eq:lttd-1-constraints-part-2-E-uniq-2});
    \item transition between states which are allowed and forbidden provided by Table \ref{tab:transitions-allowed} with (\ref{eq:lttd-1-constraints-part-3-transition-forbidden-1})-(\ref{eq:lttd-1-constraints-part-3-transition-forbidden-last}), and we assume all nuclear units are ON (OU/OFL/OD) initially (\ref{eq:lttd-1-constraints-part-3-T-ini-nuc});
    \item constraints (\ref{eq:lttd-1-constraints-part-4-OFL})-(\ref{eq:lttd-1-constraints-part-4-OFF}) related to $\text{T}^{u, \text{FLAT}}$, $\text{T}_{min}^{u, \text{ON}}$, $\text{T}_\text{max}^{u, \text{ON}}$, and $\text{T}_\text{min}^{u, \text{OFF}}$ corresponding to minimal FLAT, ON duration, maximal ON duration for OCGT units, and minimal OFF duration of thermal units;
    \item constraints (\ref{eq:lttd-1-constraints-part-5-1})-(\ref{eq:lttd-1-constraints-part-5-2}) related to the state definitions;
    \item constraint (\ref{eq:lttd-1-constraints-part-6}) related to the maximum number of start-up for an OCGT unit per day.
\end{enumerate}
\begin{subequations}
\label{eq:lttd-1-constraints-part-1}	
\begin{align}
% balance equation
0 & = \tilde{d}_{t,s}^{t^1} -\text{dns}_{t,s} + \text{spill}_{t,s} - \sum_{u \in \Omega_\text{T}} p_{t,s}^u \label{eq:lttd-1-constraints-part-1-balance} \\
% Pmax thermal units
p_{t,s}^u & \leq (1-E_{t,\text{OFF}}^u) P_\text{max}^u  \quad   u \in \Omega_\text{FT} \label{eq:lttd-1-constraints-part-1-pmax}\\
p_{t,s}^u & \leq (1-E_{t,s,\text{OFF}}^u) P_\text{max}^u  \quad s \in \Omega_S,  u \in \Omega_\text{ST}\\
% Pmin thermal units 
p_{t,s}^u & \geq (1-E_{t,\text{OFF}}^u) P_\text{min}^u  \quad   u \in \Omega_\text{FT} \\
p_{t,s}^u & \geq (1-E_{t,s,\text{OFF}}^u) P_\text{min}^u  \quad s \in \Omega_S,  u \in \Omega_\text{ST} \label{eq:lttd-1-constraints-part-1-pmin}\\
% Instantaneous power variation for first-stage units from OU to OU or OD to OD
T_{t,e,e}^u & \leq 0  \quad   u \in \Omega_\text{FT}, e \in \{\text{OU}, \text{OD}\} \label{eq:lttd-1-constraints-part-1-gradient-infinity-nuc} \\
% Instantaneous power variation for second-stage units from OU to OU or OD to OD
T_{t,s,e,e}^u & \leq 0  \quad   s \in \Omega_S, u \in \Omega_\text{ST}, e \in \{\text{OU}, \text{OD}\} \label{eq:lttd-1-constraints-part-1-gradient-infinity-fossil}.
\end{align}
\end{subequations}
\begin{subequations}
\label{eq:lttd-1-constraints-part-2}	
\begin{align}
%%%%%%%%%%%%%%%%%%%%%%%%%%%%%%%%%%%%%%%%%%%%%%%%%%%%%%
% E definition from T
% State transition definition for first-stage units t-1 -> t
E_{t,e}^u & = \sum_{e_i \in \Omega_\text{E}} T_{t-1,e_i,e}^u  \quad   u \in \Omega_\text{FT}, t > 0 \label{eq:lttd-1-constraints-part-2-E-def-1}\\ 
% State transition definition for second-stage units t-1 -> t
E_{t,s,e}^u & = \sum_{e_i \in \Omega_\text{E}} T_{t-1,s,e_i,e}^u  \quad   u \in \Omega_\text{ST}, t > 0\\
% State transition definition for first-stage units
E_{t,e}^u & = \sum_{e_f \in \Omega_\text{E}} T_{t,e,e_f}^u  \quad   u \in \Omega_\text{FT}\\
% State transition definition for second-stage units
E_{u,t,s,e}^u & = \sum_{e_f \in \Omega_\text{E}} T_{t,s,e,e_f}^u  \quad   u \in \Omega_\text{ST}\\
% State transition definition for first-stage units ini
E_{t=0,e}^u & = \sum_{e_i \in \Omega_\text{E}} T_{e_i,e}^{0,u}  \quad   u \in \Omega_\text{FT}\\ 
% State transition definition for second-stage units ini
E_{t=0,s,e}^u & = \sum_{e_i \in \Omega_\text{E}} T_{s,e_i,e}^{0,u}  \quad   u \in \Omega_\text{ST} \label{eq:lttd-1-constraints-part-2-E-def-2}\\
% State unicity for first-stage units and second-stage units
\sum_{e \in \Omega_E} E_{t,e}^u & =  1 \quad   u \in \Omega_\text{FT} \label{eq:lttd-1-constraints-part-2-E-uniq-1} \\
\sum_{e \in \Omega_E} E_{t, s, e}^u & =  1 \quad  s \in \Omega_S, u \in \Omega_\text{ST} \label{eq:lttd-1-constraints-part-2-E-uniq-2} 
\end{align}
\end{subequations}
\begin{subequations}
\label{eq:lttd-1-constraints-part-3}	
\begin{align}
%%%%%%%%%%%%%%%%%%%%%%%%%%%%%%%%%%%%%%%%%%%%%%%%%%%%%%
% State transition forbidden
% State transition forbidden NUC OU -> OF or OFF
T_{t,\text{OU},e_f}^u & \leq 0  \quad   u \in \Omega_\text{FT}, e_f \in \{\text{OD}, \text{OFF}\} \label{eq:lttd-1-constraints-part-3-transition-forbidden-1}\\ 
% State transition forbidden fossil OU -> OF or OFF
T_{t,s,\text{OU},e_f}^u & \leq 0  \quad   u \in \Omega_\text{ST}, e_f \in \{\text{OD}, \text{OFF}\} \\ 
% State transition forbidden NUC OD -> OU
T_{t,\text{OD},\text{OU}}^u & \leq 0  \quad   u \in \Omega_\text{FT}\\ 
% State transition forbidden fossil OD -> OU
T_{t,s,\text{OD},\text{OU}}^u & \leq 0  \quad   u \in \Omega_\text{ST} \\ 
% State transition forbidden NUC OFF -> OFL or OD
T_{t,\text{OFF},e_f}^u & \leq 0  \quad   u \in \Omega_\text{FT}, e_f \in \{\text{OD}, \text{OFL}\} \\ 
% State transition forbidden fossil OFF -> OFL or OD
T_{t,s,\text{OFF},e_f}^u & \leq 0  \quad   u \in \Omega_\text{ST}, e_f \in \{\text{OD}, \text{OFL}\} \\
%%%%%%%%%%%%%%%%%%%%%%%%%%%%%%%%%%%%%%%%%%%%%%%%%%%%%%
% State transition forbidden ini
% State transition forbidden NUC OU -> OF or OFF
T_{\text{OU},e_f}^{0,u} & \leq 0  \quad   u \in \Omega_\text{FT}, e_f \in \{\text{OD}, \text{OFF}\} \\ 
% State transition forbidden fossil OU -> OF or OFF
T_{s,\text{OU},e_f}^{0,u} & \leq 0  \quad   u \in \Omega_\text{ST}, e_f \in \{\text{OD}, \text{OFF}\} \\ 
% State transition forbidden NUC OD -> OU
T_{\text{OD},\text{OU}}^{0,u} & \leq 0  \quad   u \in \Omega_\text{FT}\\ 
% State transition forbidden fossil OD -> OU
T_{s,\text{OD},\text{OU}}^{0,u} & \leq 0  \quad   u \in \Omega_\text{ST} \\
%
% Notice that the two following constraints are not required when considering that all units cannot be OFF at t =-1
%
% State transition forbidden NUC OFF -> OFL or OD
T_{\text{OFF},e_f}^{0,u} & \leq 0  \quad   u \in \Omega_\text{FT}, e_f \in \{\text{OD}, \text{OFL}\} \\ 
% State transition forbidden fossil OFF -> OFL or OD
T_{s,\text{OFF},e_f}^{0,u} & \leq 0  \quad   u \in \Omega_\text{ST}, e_f \in \{\text{OD}, \text{OFL}\} \label{eq:lttd-1-constraints-part-3-transition-forbidden-last} \\
%%%%%%%%%%%%%%%%%%%%%%%%%%%%%%%%%%%%%%%%%%%%%%%%%%%%%%
% Initialization: all nuclear units are supposed to be ON at t=-1: ei= OU/OD/OFL, and they can be OFF/OU/OD/OFL at t = 0
T_{\text{OFF},e_f}^{0,u} & \leq 0  \quad   u \in \Omega_\text{N}, e_f \in \Omega_E \label{eq:lttd-1-constraints-part-3-T-ini-nuc} 
%T_{s, Off,e_f}^{0,u} & \leq 0  \quad   u \in \Omega_\text{ST},  e_f \in \Omega_E  \label{eq:lttd-1-constraints-part-3-T-ini-fossil}  
\end{align}
\end{subequations}
%
% constraints for min time FLAT, ON, and OFF durations
%
\begin{subequations}
\label{eq:lttd-1-constraints-part-4}	
\begin{align}
%%%%%%%%%%%%%%%%%%%%%%%%%%%%%%%%%%%%%%%%%%%%%%%%%%%%%%
% Min time Flat
T_{t,\text{OU},\text{OFL}}^u + T_{t,\text{OD},\text{OFL}}^u & \leq  E_{t',\text{OFL}}^u \quad   u \in \Omega_\text{FT},  t' \in \mathcal{T}^{\text{OFL}}\label{eq:lttd-1-constraints-part-4-OFL}\\ 
T_{t,s,\text{OU},\text{OFL}}^u + T_{t,s,\text{OD},\text{OFL}}^u & \leq  E_{t',s,\text{OFL}}^u \quad   u \in \Omega_\text{ST},  t' \in \mathcal{T}^{\text{OFL}} \\ 
% Min time ON
T_{t,\text{OFF},\text{OU}}^u + E_{t',\text{OFF}}^u & \leq  1 \quad   u \in \Omega_\text{FT},  t' \in \mathcal{T}^{On}_\text{min} \label{eq:lttd-1-constraints-part-4-ON}\\
T_{t,s,\text{OFF},\text{OU}}^u + E_{t',s,\text{OFF}}^u & \leq  1 \quad   u \in \Omega_\text{ST},  t' \in \mathcal{T}^{On}_\text{min} \\
% Maximum On duration for OCGT only
E_{t',s,e \neq \text{OFF}}^u + T_{t,s,\text{OFF}, \text{OU}}^u& \leq 1  \quad t' \in \mathcal{T}^{On}_\text{max},  u \in \Omega_\text{ST} \label{eq:lttd-1-constraints-part-4-ON-max}\\
% Min time OFF
T_{t,\text{OFL},\text{OFF}}^u + T_{t,\text{OD},\text{OFF}}^u & \leq  E_{t',\text{OFF}}^u \quad   u \in \Omega_\text{FT},  t' \in \mathcal{T}^{\text{OFF}} \\
T_{t,s,\text{OFL},\text{OFF}}^u + T_{t,s,\text{OD},\text{OFF}}^u & \leq  E_{t',s,\text{OFF}}^u \quad   u \in \Omega_\text{ST},  t' \in \mathcal{T}^{\text{OFF}} \label{eq:lttd-1-constraints-part-4-OFF},
\end{align}
\end{subequations}
with $\mathcal{T}^{\text{OFL}}= [t + \Delta_t, t + \text{T}^{u,\text{Flat}}]$, $\mathcal{T}^{On}_\text{min}= [t + \Delta_t, t + \text{T}^{u,On}_\text{min}]$, $\text{T}_\text{max}^{u,On} = [t + \Delta_t + T^{u,On}_\text{max} + \Delta_t, T]$, and $\mathcal{T}^{\text{OFF}}= [t + \Delta_t, t + \Delta_t + \text{T}^{u,\text{OFF}}_\text{min}]$.
\begin{subequations}
\label{eq:lttd-1-constraints-part-5}	
\begin{align}
%%%%%%%%%%%%%%%%%%%%%%%%%%%%%%%%%%%%%%%%%%%%%%%%%%%%%%
% State variable definition
p^u_{t,s} - p^u_{t-1,s} & \geq \Delta p^{min} E_{t,\text{OU}}^u \nonumber\\ 
& - P^u_\text{max} ( E_{t,\text{OD}}^u +  E_{t,\text{OFF}}^u) \label{eq:lttd-1-constraints-part-5-1} \quad   u \in \Omega_\text{FT}\\ 
p^u_{t,s} - p^u_{t-1,s} & \geq \Delta p^{min} E_{t,s,\text{OU}}^u \nonumber\\ 
& - P^u_\text{max} ( E_{t,s,\text{OD}}^u +  E_{t,s,\text{OFF}}^u)  \quad   u \in \Omega_\text{ST}\\ 
 p^u_{t-1,s} - p^u_{t,s} & \leq -\Delta p^{min} E_{t,\text{OU}}^u \nonumber\\ 
& + P^u_\text{max} ( E_{t,\text{OD}}^u +  E_{t,\text{OFF}}^u)  \quad   u \in \Omega_\text{FT}\\ 
p^u_{t-1,s} - p^u_{t,s} & \leq -\Delta p^{min} E_{t,s,\text{OU}}^u \nonumber\\ 
& + P^u_\text{max} ( E_{t,s,\text{OD}}^u +  E_{t,s,\text{OFF}}^u)  \quad   u \in \Omega_\text{ST}
\label{eq:lttd-1-constraints-part-5-2},
\end{align}
\end{subequations}
with $-\Delta p^{min}$ the minimal power variation allowed.
%
%%%%%%%%%%%%%%%%%%%%%%%%%%%%%%%%%%%%%%%%%%%%%%%%%%%%%%
% OCGT specific constraints
%%%%%%%%%%%%%%%%%%%%%%%%%%%%%%%%%%	
\begin{align} \label{eq:lttd-1-constraints-part-6}
% Maximum number of starts per day for OCGT units
\sum_{t \in \mathcal{T}} T_{t,s,\text{OFF},\text{OU}}^u + T_{s,\text{OFF},\text{OU}}^{u,0} & \leq \text{N}^\text{ON}_{max} \quad   u \in \Omega_\text{ST},
\end{align}
with $\text{N}^\text{ON}_{max}$ the maximal number of startup within a day for OCGT units.

\subsection{Two-stage stochastic MPC}\label{sec:sto-sequence}

Figure \ref{fig:multi-phase-framework} depicts the multi-phase framework timeline. It comprises several phases presented in the following paragraphs. 

\subsubsection{Phase 1: ON/OFF of nuclear units} 

The first two-stage stochastic optimizer at $t=t_1$ consists of computing the ON/OFF status of the nuclear units over $[T_1, T_3]$.
Thus, the first-stage variables are the ON/OFF status of the nuclear units over $[T_1, T_3]$.  
Then, the second-stage variables comprise the ON/OFF status of COAL, CCGT, and OCGT and the production level of all units for $[T_1, T_3]$.
% Objective function
The objective function to optimize is (\ref{eq:single-obj-details}) where  $\Omega_\text{FT}$ is the set of first-stage units (nuclear units), $\Omega_\text{ST}$ is the set of second-stage units (coal, CCGT, and OCGT units), and $\Omega_\text{T}= \Omega_\text{FT} \cup \Omega_\text{ST}$ the set of all thermal units.

\subsubsection{Phase 2: ON/OFF of coal units}

The second two-stage stochastic optimizer at $t=t_2$ consists of computing the ON/OFF status of the coal units over $[T_1, T_2[$ using the ON/OFF status of the nuclear units, from the previous second-stage problem solved at $t_1$, as inputs.
The first-stage variables are the ON/OFF status of the coal units over $[T_1, T_2[$. 
The second-stage variables comprise the ON/OFF status of CCGT and OCGT and the production level of all units over $[T_1, T_3]$. It also comprises the ON/OFF status of coal units for $[T_2, T_3]$.
The optimization problem is
\begin{subequations}
\label{eq:coal-lttd-obj}	
\begin{align}
\min \  &  J_{f,t_2}^1 + J_{f,t_2}^2 + J_{v,t_2} , \\ 
&  \text{s.t. } \notag \\
& \quad (\ref{eq:lttd-1-constraints-part-1-balance}) - (\ref{eq:lttd-1-constraints-part-6}), (\ref{eq:lttd-coal-constraints-nuc-on-off}).
\end{align}
\end{subequations}
$J_{f,t_2}^1$, $J_{f,t_2}^2$, and $J_{v,t_2}$ are defined as follows
\begin{subequations}
\label{eq:coal-obj-details}	
\begin{align}
% first-stage cost
J_{f,t_2}^1  & = \sum_{u \in \Omega_\text{FT}^{t_2}}  \pi_f^u  \bigg[\sum_{t\in [T_1, T_2[}  T_{t,\text{OFF},\text{OU}}^u +   T_{\text{OFF},\text{OU}}^{0,u} \bigg] , \\ 
% second-stage fixed cost
J_{f,t_2}^2  & =  \sum_{u \in \Omega_\text{FT}^{t_2}}  \sum_{s \in \Omega_S} p_s \pi_f^u\bigg[\sum_{t\in [T_2, T_3]}   T_{t,s,\text{OFF},\text{OU}}^u +   T_{s,\text{OFF},\text{OU}}^{0,u} \bigg]\\
& + \sum_{u \in \Omega_\text{ST}^{t_2}}  \sum_{s \in \Omega_S} p_s \pi_f^u\bigg[\sum_{t\in \mathcal{T}}   T_{t,s,\text{OFF},\text{OU}}^u +   T_{s,\text{OFF},\text{OU}}^{0,u} \bigg], \\
% second-stage variable cost
J_{v,t_2}  & =  \sum_{t\in \mathcal{T}} \sum_{s \in \Omega_S} p_s \bigg[ \sum_{u \in \Omega_\text{T}} \pi_v^u p_{t,s}^u +  \pi^\text{dns} \text{dns}_{t,s}  + \pi^\text{spill} \text{spill}_{t,s} \bigg],
\end{align}
\end{subequations}
where $\Omega_\text{FT}^{t_2}$ is the set of first-stage units (coal units), $\Omega_\text{ST}^{t_2}$ is the set of second-stage units (CCGT and OCGT), and $\Omega_\text{T}$ the set of all thermal units (nuclear, coal, CCGT, and OCGT units). 
%
% ON/OFF status of nuclear units fixed
The constraints related to the ON/OFF status of nuclear units computed in the previous two-stage stochastic optimization are
\begin{align} 
\label{eq:lttd-coal-constraints-nuc-on-off}
E_{t, e=\text{OFF}}^{u \in \Omega_\text{N}} & = E_{t, e=\text{OFF}}^{\star, u \in \Omega_\text{N}} \quad \forall t \in \mathcal{T},
\end{align}
with $E_{t, e=\text{OFF}}^{\star, u \in \Omega_\text{N}}$ $\forall t \in \mathcal{T}$ the values of the ON/OFF status of nuclear units computed in the previous two-stage stochastic optimization.
%
% scenarios
%
The scenarios of consumption and renewable generations of this problem issued at $t_2$ for $[T_1, T_3]$ are $\hat{d}_{t,s}^{t^2}$ and $\hat{p}_{t,s}^{t^2, u \in \Omega_\text{Nr}}$, $\forall t \in \mathcal{T}$.
In the constraint related to the power balance (\ref{eq:lttd-1-constraints-part-1-balance}), the residual demand becomes $\tilde{d}_{t,s}^{t^2}$.

\subsubsection{Phase 3: ON/OFF of CCGT units}

The third two-stage stochastic optimizer at $t=t_{31}$ consists of computing the ON/OFF status of the CCGT units over $[T_1, T_{12}[$ using the ON/OFF status of the nuclear and coal units, computed in the two previous two-stage stochastic problems, as inputs. 
The first-stage variables are the ON/OFF status of the CCGT units over $[T_1, T_{12}[$. 
The second-stage variables comprise the ON/OFF status of OCGT and the level of production of all units over $[T_1, T_3]$. It also comprises the ON/OFF status of coal units over $[T_2, T_3]$ and the ON/OFF status of CCGT units for $[T_{12}, T_3]$.
The optimization problem is
\begin{subequations}
\label{eq:CCG-lttd-obj}	
\begin{align}
\min \  &  J_{f,t_{31}}^1 + J_{f,t_{31}}^2 + J_{v,t_{31}} , \\ 
&  \text{s.t. } \notag \\
& \quad (\ref{eq:lttd-1-constraints-part-1-balance}) - (\ref{eq:lttd-1-constraints-part-6}), (\ref{eq:lttd-coal-constraints-nuc-on-off}), (\ref{eq:lttd-CCG-constraints-coal-on-off}).
\end{align}
\end{subequations}
$J_{f,t_{31}}^1$, $J_{f,t_{31}}^2$, and $J_{v,t_{31}}$ are defined as follows
\begin{subequations}
\label{eq:CCG-obj-details}	
\begin{align}
% first-stage cost
J_{f,t_{31}}^1  & = \sum_{u \in \Omega_\text{FT}^{t_{31}}}  \pi_f^u  \bigg[\sum_{t\in [T_1, T_{12}[}  T_{t,\text{OFF},\text{OU}}^u +   T_{\text{OFF},\text{OU}}^{0,u} \bigg] , \\ 
% second-stage fixed cost
J_{f,t_{31}}^2  & =  \sum_{u \in \Omega_\text{FT}^{t_{31}}}  \sum_{s \in \Omega_S} p_s \pi_f^u\bigg[\sum_{t\in [T_{12}, T_3]}   T_{t,s,\text{OFF},\text{OU}}^u +   T_{s,\text{OFF},\text{OU}}^{0,u} \bigg]\\
& + \sum_{u \in \Omega_\text{FT}^{t_2}}  \sum_{s \in \Omega_S} p_s \pi_f^u\bigg[\sum_{t\in [T_2, T_3]}   T_{t,s,\text{OFF},\text{OU}}^u +   T_{s,\text{OFF},\text{OU}}^{0,u} \bigg]\\
& + \sum_{u \in \Omega_\text{ST}^{t_{31}}}  \sum_{s \in \Omega_S} p_s \pi_f^u\bigg[\sum_{t\in \mathcal{T}}   T_{t,s,\text{OFF},\text{OU}}^u +   T_{s,\text{OFF},\text{OU}}^{0,u} \bigg], \\
% second-stage variable cost
J_{v,t_{31}}  & =  \sum_{t\in \mathcal{T}} \sum_{s \in \Omega_S} p_s \bigg[ \sum_{u \in \Omega_\text{T}} \pi_v^u p_{t,s}^u +  \pi^\text{dns} \text{dns}_{t,s}  + \pi^\text{spill} \text{spill}_{t,s} \bigg],
\end{align}
\end{subequations}
where $\Omega_\text{FT}^{t_{31}}$ is the set of first-stage units (CCGT units), $\Omega_\text{ST}^{t_{31}}$ is the set of second-stage units (OCGT units), and $\Omega_\text{T}$ the set of all thermal units (nuclear, coal, CCGT, and OCGT units). 
Notice that in $J_{f,t_{31}}^2$ the set of units $ \Omega_\text{FT}^{t_2}$ is considered. It referred to the coal units with the ON/OFF status as second-stage variables over $[T_2, T_3]$.
%
% ON/OFF status of coal units fixed
The constraint related to the ON/OFF status of coal units is
\begin{align} 
\label{eq:lttd-CCG-constraints-coal-on-off}
E_{t, e=\text{OFF}}^{u \in \Omega_\text{FT}^{t_2}} & = E_{t, e=\text{OFF}}^{\star, u \in \Omega_\text{FT}^{t_2}} \quad \forall t \in [T_1, T_2[
\end{align}
with $E_{t \in [T_1, T_2[, e=\text{OFF}}^{\star, u \in \Omega_\text{FT}^{t_2}}$ the values of the ON/OFF status of coal units computed in the second two-stage stochastic problem.
%
% scenarios
The scenarios of consumption and renewable generations of this problem issued at $t_{31}$ are $\hat{d}_{t,s}^{t_{31}}$ and $\hat{p}_{t,s}^{t_{31}, u \in \Omega_\text{Nr}}$, $\forall t \in \mathcal{T}$.
Then, the constraint related to the power balance (\ref{eq:lttd-1-constraints-part-1-balance}), $\tilde{d}_{t,s}^{t_{31}}$ becomes the residual demand.

\subsubsection{Phase 4: sequence of ON/OFF of CCGT units}

Then, the two-stage stochastic optimizers at $t=t_{3i}$, with $2 \leq i \leq 6$, consist of computing the ON/OFF status of the CCGT units over $[T_{1i}, T_{1(i+1)}[$ using the ON/OFF status of the nuclear and coal units as inputs. 
The first-stage variables are the ON/OFF status of the CCGT units over $[T_{1i}, T_{1(i+1)}[$. 
The second-stage variables comprise the ON/OFF status of OCGT and the production level of all units over $[T_{1i}, T_3]$. It also comprises the ON/OFF status of coal units over $[T_2, T_3]$ and the ON/OFF status of CCGT units for $[T_{1(i+1)}, T_3]$.

\section{Numerical experiments}\label{sec:results}

The numerical experiments are conducted using the single-phase framework presented in section \ref{sec:sto-simple} in a realistic case study where RTE faced a deficit of downward margins.
The stochastic optimizer is implemented using the Python API of FICO Xpress Optimization\footnote{\url{https://www.fico.com/fico-xpress-optimization/docs/latest/overview.html}}, and the problem is solved with the FICO Xpress solver.

This section presents the ex-post \textit{out-of-sample} evaluation approach adopted and the deterministic formulation used to evaluate the ON/OFF plans computed by the stochastic optimizer.
Then, the results are commented on.

\subsection{Case study}

\begin{figure}[tb]
\centerline{\includegraphics[width=\linewidth]{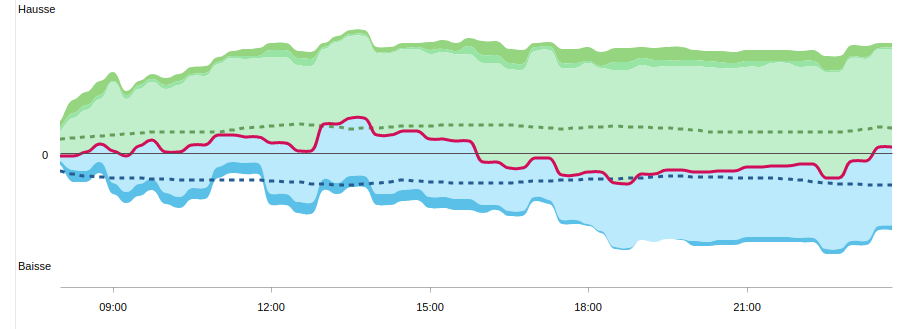}}
\caption{Required (dashed lines) and available (plain lines) upward (green) and downward (blue) margins computed by RTE at 8:30 a.m. on Saturday, 15 July 2023, by steps of 15 minutes until the end of the day. The red line is the forecasted system imbalance.}
\label{fig:marges-15-juillet-2023}
\end{figure}

The case study considers Saturday, 15 July 2023, where RTE faced a deficit of downward margins. Figure \ref{fig:marges-15-juillet-2023} depicts the required (dashed lines) and available (plain lines) upward (green) and downward (blue) margins computed by RTE at 8:30 a.m. by steps of 15 minutes until the end of the day. At 8:30 a.m., the deficit of downward margins is anticipated to be around 10 and 12 a.m. 
Notice that the numerical experiments consider the specific day of Saturday, 15 July 2023. However, it can be extended more generally to other days where RTE faced a deficit of downward margins.
In addition, the approach considered would also be valid for situations with a deficit of upward margins.

% Explanations why this case study is interesting
We selected the downward deficit case as these deficits are increasingly frequent due to the increased penetration of renewable energies. They tend to occur during weekends from spring to autumn, where the consumption tends to be small and renewable production is high. However, these situations also tend to happen in winter as the temperatures get milder.
In these situations, preserving enough available downward margins is sometimes challenging as the system faces excess production. Thus, most fossil-based flexible production, such as CCGT and OCGT, is usually shut down. However, this is not always enough, and it may be necessary to decrease the production of nuclear power plants and shut down some of them. 
Furthermore, the minimum OFF duration of nuclear power plants is at least 24 hours. Thus, the decision to shut down nuclear power plants implies carefully anticipating over at least 24 hours a potential need for upward margin. 
Indeed, in the short term, shutting down several nuclear power plants could help to restore downward available margins. However, if the system needs upward margins a few hours later, these power plants cannot be used, resulting in a deficit of upward margins. This typically could happen during the night when the consumption is small or in the middle of the day. Then, it is followed by the consumption peak in the morning or evening. 
Therefore, fossil-based units, which are more expensive and emit GHG, will be used instead to compensate for the shutdown of nuclear power plants. 
This situation illustrates the need for a tool to help make these decisions under uncertainty. 

% Unit considered
Table \ref{tab:production-units} provides the units considered in this case study for optimization. The other units connected to the French electrical system are assumed to follow their production plan. 
Indeed, most units usually produce at a constant production level, and only a tiny part modulates to adjust the variation of the residual demand.
On Saturday, 15 July 2023, the nuclear units considered in the optimization were the ones that modulated this day. Three were shut down, and the other three decreased production to the minimum. Notice that RTE did not order these units to modulate. However, to face the deficit of downward margin optimally, it is interesting to consider them in the optimization to understand how a probabilistic UC tool could use them. 
The values of the technical parameters related to the technical constraints, such as minimal time ON/OFF, are provided by Table \ref{tab:technical-constraints}.
\begin{table}[tb]
\renewcommand{\arraystretch}{1.25}
\begin{center}
\begin{tabular}{lrrrrrr}
\hline \hline
Unit & $P_\text{min}$ & $P_\text{max}$ & $\pi_f$ & $\pi_v$\\
\hline
% nuclear
NUC 1 & 180 & 915 & 25 & 10  \\
NUC 2 & 180 & 930 & 25 & 10  \\
NUC 3 & 182 & 937 & 25 & 10  \\
NUC 4 & 260 & 1320 & 37 & 10 \\
NUC 5 & 260 & 1297 & 37 & 10 \\
NUC 6 & 375 & 1532 & 42 & 10  \\
% fossil
CCGT 1 &  180 & 440 & 10 & 30  \\
CCGT 2 &  214 & 446 & 10 & 30  \\
CCGT 3 &  196 & 455 & 10 & 30  \\
OCGT 1 &  110 & 181 & 5 & 150 \\
OCGT 2 &  150 & 179 & 5 & 150 \\
OCGT 3 &  106 & 182 & 5 & 150 \\
\hline \hline
\end{tabular}
\caption{Parameters of the thermal units considered. Notice that these start-up and variable cost values are not necessarily those used by the market players.}
\label{tab:production-units}
\end{center}
\end{table}

\subsection{Ex-post \textit{out-of-sample} evaluation}

%%%%%%%%%%%%%%%%%%%%%%% 
% Principles
%%%%%%%%%%%%%%%%%%%%%%%
An ex-post \textit{out-of-sample} simulation is the most common approach \cite{ROALD2023108725} to evaluate the quality of a solution, independent of the approach used for decision-making under uncertainty, such as stochastic programming, robust optimization, chance-constrained optimization, and distributionally robust optimization. 
For given sources of uncertainty $\boldsymbol \xi$ (in our case, PV, wind power, and consumption), consider $N$ samples obtained from a probabilistic forecast or empirical observations, each representing a potential realization of $\xi$. Then, these samples are split arbitrarily into two separate sets of samples, each with $M$ and $K$ samples, such that $M+K = N$ and $K > M$. The set with $M$ samples is used as the training set to conduct the stochastic optimizations. This set solves the problem with probabilistic constraints to obtain the in-sample value for the objective function and a solution for the decision variables.
For instance, the in-sample value $\kappa^{in}$ in the case of a two-stage stochastic optimization (\ref{eq:two-stage}) is $\kappa^{in} = f^F(\boldsymbol x^\star) + \mathbb{E}[f^S(\boldsymbol x^\star, \boldsymbol y_{\xi}^\star, \boldsymbol \xi)]$ where $\boldsymbol x^\star$ and $\boldsymbol y_{\xi}^\star$ are the optimal values obtained for the first-stage and the second-stage variables, according to the $M$ samples used in the training set.
Then, the values of the first-stage variables are fixed to that
achieved in the in-sample simulation, \textit{i.e.}, $\boldsymbol x^\star$, and deterministically solve the second-stage problem $K$ times, each time using a sample that has not been used in the in-sample simulation. Hence, the set of $K$ samples is called the testing set.
For each unseen sample $\xi_i$, the following deterministic problem is solved
\begin{subequations}
\begin{align}
        \min_{\boldsymbol y_{\xi}} & \ \mathbb{E}[f^S(\boldsymbol x^\star, \boldsymbol y_{\xi}, \boldsymbol\xi)] \\
        \text{s.t. } & h^F(\boldsymbol x^\star) = 0 \\
        &g^F(\boldsymbol x^\star) \leq 0 \\
        &h^S(\boldsymbol x^\star, \boldsymbol y_{\xi}, \boldsymbol \xi) = 0 \\
        &g^S(\boldsymbol x^\star, \boldsymbol y_{\xi}, \boldsymbol \xi) \leq 0.
\end{align}
\label{eq:det-eval}
\end{subequations}
Then, the out-of-sample value of the objective function, denoted as $\kappa^{oos}$, is obtained by summing up the first-stage value  $f^F(\boldsymbol x^\star)$ achieved from the in-sample simulation and the average second-stage value $1/K \sum_{i=1}^K  [f^S(\boldsymbol x^\star, \boldsymbol y_{\xi_i},\xi_i)]$ obtained over the $K$ solutions to the deterministic problem (\ref{eq:det-eval})
\begin{align}
\kappa^{oos} & = F(x^\star) + 1/K  \sum_{i=1}^K [f^S(x^\star,y_{\xi_i},\xi_i)].
\end{align}
The difference between $\kappa^{in}$ and $\kappa^{oos}$ is a metric for the solution quality.
A comparatively significant difference indicates a lower-quality solution, implying that the training set with $M$ samples needs to be more sufficient to represent the underlying uncertainty properly. Therefore, the training set should be improved by increasing the number of samples and/or selecting more representative samples.
In our case study, we followed this evaluation strategy. However, we are also interested in other KPIs, such as the lost load and production.

%%%%%%%%%%%%%%%%%%%%%%% 
% In practice
%%%%%%%%%%%%%%%%%%%%%%%

\begin{figure}[tb]
\centerline{\includegraphics[width=\linewidth]{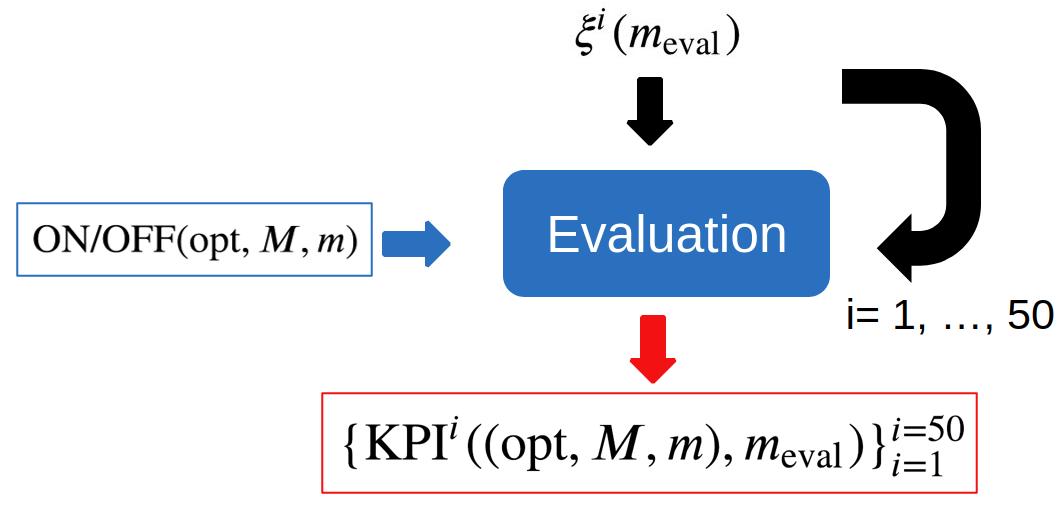}}
\caption{Overview of the ex-post \textit{out-of-sample} evaluation for a given $\text{ON/OFF}(\text{opt}, M, m)$ plan. This ON/OFF plan was computed with $M$ scenarios, randomly selected, with a maximal deviation $m$ \% for a three-sigma confidence interval, and the optimizer opt.}
\label{fig:strategie-evaluation}
\end{figure}

Figure \ref{fig:strategie-evaluation} depicts the overall ex-post \textit{out-of-sample} evaluation process.
In practice, the ON/OFF plan computed is $\text{ON/OFF}(\text{opt}, M, m)$ where $\text{opt}$ is the optimizer, $M$ the number of random scenarios selected, and $m$ the maximal deviation of scenarios for a three-sigma confidence interval.
The greater $m$, the more variable the scenarios generated. Appendix \ref{appendix:scenario-generation} presents the scenario approach method.
The optimizer ($\text{opt}$) can be the stochastic optimizer of the single-phase/multi-phase sequential or any other optimizer.
Then, it computes several $\text{ON/OFF}(\text{opt}, M, m)$ plans for several pairs $(M, m)$.

For a given optimizer ($\text{opt}$), the evaluation is performed for each $\text{ON/OFF}(\text{opt}, M, m)$ with an evaluation scenario $\xi^i(m_\text{eval})$, where $m_\text{eval}$ is the maximal deviation for a three-sigma confidence interval.
It allows to compute the KPI $\text{KPI}^i((\text{opt}, M, m), m_\text{eval})$, such as the lost load/production or the dispatch costs.
Then, by applying this process recursively $\forall i \in [1, K]$ the set of KPIs $\{ \text{KPI}^i((\text{opt}, M, m), m_\text{eval})\}_{i=1}^{K}$ is computed.

\begin{figure}[tb]
\centerline{\includegraphics[width=\linewidth]{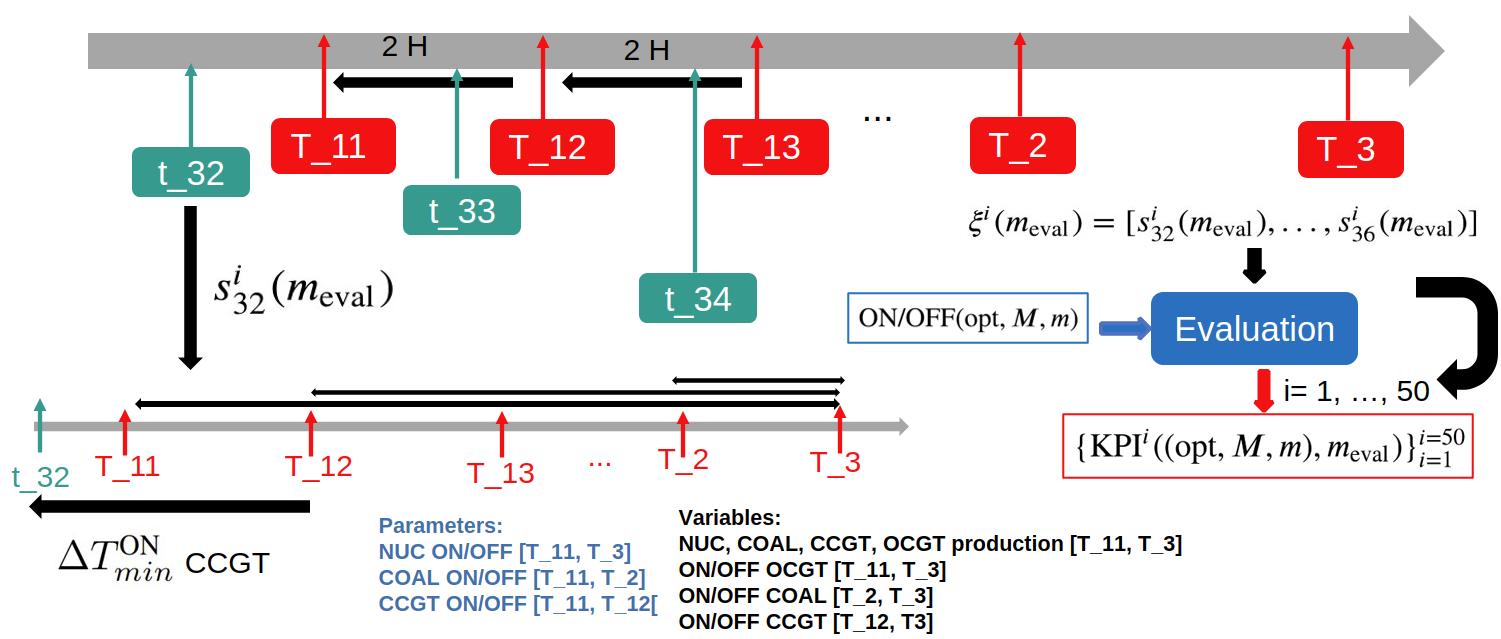}}
\caption{Detailed overview of the ex-post \textit{out-of-sample} evaluation with a focus at $t=t_{32}$. For a given evaluation scenario $\xi^i(m_\text{eval}) = [s_{32}^i(m_\text{eval}), \ldots, s_{36}^i(m_\text{eval})]$, a sequence of deterministic optimizations is solved. At $t=t_{32}$, the scenario $s_{32}^i(m_\text{eval})$ is used as input, and the variables are the ON/OFF status for $[T_{12}, T_3]$ of the CCGT units, ON/OFF status of coal units for $[T_2, T_3]$, ON/OFF status of OCGT units for $[T_{11}, T_3]$, and the level of production of all units for $[T_{11}, T_3]$.}
\label{fig:strategie-evaluation-details}
\end{figure}

In the case study, the period $[T_1, T_2]$ is split into intervals of 2 hours to model the duration of the operational window where the TSO can activate upward and downward offers to face a deficit of margins or system imbalances.
Figure \ref{fig:strategie-evaluation-details} presents the sequence of evaluation for a given scenario $\xi^i(m_\text{eval}) = [s_{32}^i(m_\text{eval}), \ldots, s_{36}^i(m_\text{eval})]$, where $s_{32}^i(m_\text{eval})$ is a scenario of PV, wind power, and consumption computed at $t=t_{32}$ for $[T_{11}, T_3]$, \ldots, and $s_{36}^i(m_\text{eval})$ is a scenario of PV, wind power, and consumption computed at $t=t_{36}$ for $[T_{15}, T_3]$. 
Then, a sequence of deterministic optimizations is performed at $t=t_{32}$, $t=t_{33}$, \ldots, $t=t_{36}$ over $[T_{11}, T_3]$, $[T_{12}, T_3]$, \ldots, $[T_{15}, T_3]$ using $s_{32}^i(m_\text{eval})$, $s_{33}^i(m_\text{eval})$, \ldots, and $s_{36}^i(m_\text{eval})$. 
The production plan for each period of two hours $[T_{11}, T_{12}[$, $[T_{12}, T_{13}[$, \ldots, $[T_{15}, T_2[$ are retrieved sequentially and aggregated to build the production plan over $[T_{11}, T_2[$.
Notice that each optimization at $t=t_{33}$, \ldots, $t=t_{36}$ requires as input the production level at the previous period as the initial conditions. They are computed by conducting a deterministic optimization at $t=t_{32}$, $t=t_{33}$, \ldots, $t=t_{36}$ using the "best forecasts" (mean scenario). 

% example eval t32
For instance, at $t=t_{32}$, the inputs are $s_{32}^i(m_\text{eval})$ and the ON/OFF status of nuclear units for $[T_1, T_3]$, the ON/OFF status of coal units for $[T_1, T_2[$, and the ON/OFF status of CCGT units for $[T_{11}, T_{12}[$. 
The variables of this dispatch problem are the ON/OFF status of the CCGT units over $[T_{12}, T_3]$, the ON/OFF status of coal units for $[T_2, T_3]$, the ON/OFF status of OCGT units for $[T_{11}, T_3]$, and the level of production of all units for $[T_{11}, T_3]$.
This problem allows us to compute the production plan of all units over $[T_1, T_3]$, but only the values of the production level over $[T_{11}, T_{12}[$ are retrieved to compute the KPIs.
% Example eval t33
For instance, at $t=t_{33}$, the inputs are $s_{33}^i(m_\text{eval})$ and the ON/OFF status of nuclear units for $[T_1, T_3]$, the ON/OFF status of coal units for $[T_1, T_2[$, the ON/OFF status of CCGT units for $[T_{12}, T_{13}[$, and the production level of all units over $[T_{11}, T_{12}[$ computed by a deterministic optimization using the "best forecasts" issued at $t=t_{32}$ over $[T_1, T_3]$. 
The variables of this dispatch problem are the ON/OFF status of the CCGT units over $[T_{12}, T_3]$, the ON/OFF status of coal units for $[T_2, T_3]$, the ON/OFF status of OCGT units for $[T_{12}, T_3]$, and the level of production of all units for $[T_{12}, T_3]$.
This problem allows us to compute the production plan of all units over $[T_{12}, T_3]$, but only the values of the production level over $[T_{12}, T_{13}[$ are retrieved to compute the KPIs.

Finally, $\text{KPI}^i((\text{opt}, M, m), m_\text{eval})$ is computed using the results of the sequence of deterministic optimizations for an evaluation scenario $\xi^i(m_\text{eval})$.
Indeed, we sum for each period of 2 hours the volume of lost load and lost production (MWh) and the dispatch costs using the results of the deterministic optimizations.
This process is done recursively $\forall i \in [1, K]$ allowing to compute the set $\{ \text{KPI}^i((\text{opt}, M, m), m_\text{eval})\}_{i=1}^{K}$.

\subsection{Deterministic formulation for the ex-post \textit{out-of-sample} evaluation}\label{sec:det-eval}

In the ex-post \textit{out-of-sample} evaluation, the ON/OFF status of the first-stage variables $E_{t \in \mathcal{T}, e=\text{OFF}}^{u \in \Omega_\text{FT}}$ are parameters and provided by the stochastic optimizer.
Overall, the optimization problem is identical to the two-stage stochastic optimization in the single-phase framework with only one scenario and with the variables $E_{t \in \mathcal{T}, e=\text{OFF}}^{u \in \Omega_\text{FT}}$ being parameters.
% scenarios
The scenarios of consumption and renewable generations of this problem issued at $t_{3(i+1)}$ for $[T_{1i}, T_3]$ (with $1\leq i \leq 5$) are $\hat{d}_{t,s}^{t^{3(i+1)}}$ and $\hat{p}_{t,s}^{t^{3(i+1)}, u \in \Omega_\text{Nr}}$, $\forall t \in \mathcal{T}$.

In the single-phase approach, the ON/OFF status of the nuclear, coal, and CCGT units are computed over $[T_1, T_3]$ in one two-stage optimization problem.
In this case, $E_{t \in \mathcal{T}, e=\text{OFF}}^{u \in \Omega_\text{FT}}$ are parameters with $\Omega_\text{FT}$ the set of nuclear, coal, and CCGT units.
In the multi-phase approach, the first three two-stage stochastic problems allow to compute the ON/OFF status of nuclear units over $[T_1, T_3]$, the ON/OFF status of coal units over $[T_1, T_2[$, and the ON/OFF status of CCGT units over $[T_{11}, T_{12}[$.
Then, a sequence of two-stage stochastic problems computes the ON/OFF status of CCGT units sequentially by intervals of two hours. As a result, at the end of the process, the values of the ON/OFF status of CCGT units are computed over $[T_1, T_2]$.

The evaluation is conducted over periods of two hours.
Figure \ref{fig:strategie-evaluation-details} depicts the ex-post \textit{out-of-sample} evaluations sequence. For each period of 2 hours within the period of study $[T_1, T_2]$, we sample $K$ scenarios at the given $t_{3i}$ and solve $K$ deterministic optimization problems.
For instance, at $t=t_{32}$, the variables are the ON/OFF status for $[T_{12}, T_3]$ of the CCGT units, ON/OFF status of coal units for $[T_2, T_3]$, ON/OFF status of OCGT units for $[T_{11}, T_3]$, and the level of production of all units for $[T_{11}, T_3]$.

%%%%%%%%%%%%%%%%%%%%%%%%%%%%%%%%%%%%
% Evaluation at t_32
%%%%%%%%%%%%%%%%%%%%%%%%%%%%%%%%%%%%
At $t=t_{32}$, the evaluation optimization problem is
\begin{subequations}
\label{eq:eval-eval-t32-obj}	
\begin{align}
\min \  &  J_{f,t_{32}}^\text{eval} + J_{v,t_{32}}^\text{eval} , \\ 
&  \text{s.t. } \notag \\
& \quad (\ref{eq:lttd-1-constraints-part-1-balance}) - (\ref{eq:lttd-1-constraints-part-6}), (\ref{eq:lttd-coal-constraints-nuc-on-off}), (\ref{eq:lttd-CCG-constraints-coal-on-off}), (\ref{eq:lttd-eval-t32-constraints-CCG-on-off}).
\end{align}
\end{subequations}
$J_{f,t_{32}}^\text{eval}$ and $ J_{v,t_{32}}^\text{eval} $ are defined as follows
\begin{subequations}
\label{eq:eval-t32-obj-details}	
\begin{align}
% OCGT
J_{f,t_{32}}^\text{eval}  & = \sum_{u \in \Omega_\text{OCGT}}  \pi_f^u  \bigg[\sum_{t\in \mathcal{T}}  T_{t,\text{OFF},\text{OU}}^u +   T_{\text{OFF},\text{OU}}^{0,u} \bigg] , \\ 
% CCGT
& + \sum_{u \in \Omega_\text{CCGT}}  \pi_f^u  \bigg[\sum_{t \in [T_{12}, T_3]}  T_{t,\text{OFF},\text{OU}}^u +   T_{\text{OFF},\text{OU}}^{0,u} \bigg] , \\ 
% COAL
& + \sum_{u \in \Omega_\text{COAL}}  \pi_f^u  \bigg[\sum_{t \in [T_2, T_3]}  T_{t,\text{OFF},\text{OU}}^u +   T_{\text{OFF},\text{OU}}^{0,u} \bigg] , \\ 
% Variable costs
J_{v,t_{32}}^\text{eval} & =  \sum_{t\in \mathcal{T}} \bigg[ \sum_{u \in \Omega_\text{T}} \pi_v^u p_{t}^u +  \pi^\text{dns} \text{dns}_{t}  + \pi^\text{spill} \text{spill}_{t} \bigg],
\end{align}
\end{subequations}
where $\Omega_\text{OCGT}$ is the set of OCGT units, $\Omega_\text{CCGT}$ is the set of CCGT units, $\Omega_\text{COAL}$ is the set of COAL units, and $\Omega_\text{T}$ is the set of all thermal units (nuclear, coal, CCGT, and OCGT units). The constraint related to the ON/OFF status of CCGT units is
%
% ON/OFF status of CCGT units fixed
\begin{align} 
\label{eq:lttd-eval-t32-constraints-CCG-on-off}
E_{t \in [T_{11}, T_{12}[, e=\text{OFF}}^{u \in \Omega_\text{CCGT}} & = E_{t \in [T_{11}, T_{12}[, e=\text{OFF}}^{\star, u \in \Omega_\text{CCGT}},
\end{align}
with $E_{t \in [T_{11}, T_{12}[, e=\text{OFF}}^{\star, u \in \Omega_\text{CCGT}}$ the result of the third two-stage stochastic optimizer providing the ON/OFF status of CCGT units for $[T_{11}, T_{12}[$.

%%%%%%%%%%%%%%%%%%%%%%%%%%%%%%%%%%%%
% Evaluation at t_3i
%%%%%%%%%%%%%%%%%%%%%%%%%%%%%%%%%%%%
At $t=t_{3i}$ with $ 3\leq i \leq 6$, the evaluation optimization problem is
\begin{subequations}
\label{eq:eval-eval-t3i-obj}	
\begin{align}
\min \  &  J_{f,t_{3i}}^\text{eval} + J_{v,t_{3i}}^\text{eval} , \\ 
&  \text{s.t. } \notag \\
& \quad (\ref{eq:lttd-1-constraints-part-1-balance}) - (\ref{eq:lttd-1-constraints-part-6}), (\ref{eq:lttd-coal-constraints-nuc-on-off}), (\ref{eq:lttd-CCG-constraints-coal-on-off}), (\ref{eq:lttd-eval-t3i-constraints-CCG-on-off}).
\end{align}
\end{subequations}
$J_{f,t_{3i}}^\text{eval}$ and $ J_{v,t_{3i}}^\text{eval} $ are defined as follows
\begin{subequations}
\label{eq:eval-t3i-obj-details}	
\begin{align}
% OCGT
J_{f,t_{3i}}^\text{eval}  & = \sum_{u \in \Omega_\text{OCGT}}  \pi_f^u  \bigg[\sum_{t\in [T_{1(i-1)}, T_3]}  T_{t,\text{OFF},\text{OU}}^u +   T_{\text{OFF},\text{OU}}^{0,u} \bigg] , \\ 
% CCGT
& + \sum_{u \in \Omega_\text{CCGT}}  \pi_f^u  \bigg[\sum_{t \in [T_{1i}, T_3]}  T_{t,\text{OFF},\text{OU}}^u +   T_{\text{OFF},\text{OU}}^{0,u} \bigg] , \\ 
% COAL
& + \sum_{u \in \Omega_\text{COAL}}  \pi_f^u  \bigg[\sum_{t \in [T_2, T_3]}  T_{t,\text{OFF},\text{OU}}^u +   T_{\text{OFF},\text{OU}}^{0,u} \bigg] , \\ 
% Variable costs
J_{v,t_{3i}}^\text{eval} & =  \sum_{t\in [T_{1(i-1)}, T_3]} \bigg[ \sum_{u \in \Omega_\text{T}} \pi_v^u p_{t}^u +  \pi^\text{dns} \text{dns}_{t}  + \pi^\text{spill} \text{spill}_{t} \bigg].
\end{align}
\end{subequations}
The constraint related to the ON/OFF status of CCGT units is
% ON/OFF status of CCG units fixed
\begin{align} 
\label{eq:lttd-eval-t3i-constraints-CCG-on-off}
E_{t \in [T_{1(i-1)}, T_{1i}[, e=\text{OFF}}^{u \in \Omega_\text{CCGT}} & = E_{t \in [T_{1(i-1)}, T_{1i}[, e=\text{OFF}}^{\star, u \in \Omega_\text{CCGT}},
\end{align}
with $E_{t \in [T_{1(i-1)}, T_{1i}[, e=\text{OFF}}^{\star, u \in \Omega_\text{CCGT}}$ the result of the two-stage stochastic optimizer at $t_{3(i-1)}$ providing the ON/OFF status of CCGT units for $[T_{1(i-1)}, T_{1i}[$.

\subsection{Results using the single-phase framework}

Table \ref{tab:on-off-computations} presents the combination of pairs $(M=5, 10, 20, 50, m=2, 5, 10, 25 \%)$ for computing several $\text{ON/OFF}(\text{opt}, M, m)$ plans using the two-stage stochastic optimizer ($\text{opt} = \text{sto}$) of the single-phase framework presented in Section \ref{sec:sto-simple}. 
This optimizer is compared to a version ($\text{opt} = \text{sto}^\star$), where the second-stage binary variables are relaxed to become continuous variables between 0 and 1, and to the deterministic optimizer $(\text{opt}=\text{det})$, which uses the best-forecasts of PV, wind power, and consumption.

%%%%%%%%%%%%%%%%%%%%%%%%%%%%%
% Computation time
%%%%%%%%%%%%%%%%%%%%%%%%%%%%%
Table \ref{tab:on-off-computations} provides the computation time of each of $\text{ON/OFF}(\text{opt}, M, m)$ plan given the number of threads used and the MIPRELSTOP\footnote{\url{https://www.fico.com/fico-xpress-optimization/docs/latest/solver/optimizer/HTML/MIPRELSTOP.html}} criterion fixed. The crosses indicate the optimization did not converge within the time allowed for a given MIPRELSTOP criterion. This value determines when the branch and bound tree search will terminate. Seven computations using the stochastic optimizer, and only one with the relaxed version ($\text{opt} = \text{sto}^\star$) did not reach at least a convergence value of MIPRELSTOP $= 10\%$ within four days.
Then, the ex-post \textit{out-of-sample} evaluation is performed for each $\text{ON/OFF}(\text{opt}, M, m)$ plan with $m_\text{eval}=$ 5, 10, and 25 \% over $K=$ 50 scenarios $\xi^i(m_\text{eval})$ $\forall i \in [1, 50]$. Thus, for each ON/OFF plan the set $\{ \text{KPI}^i((\text{opt}, M, m), m_\text{eval})\}_{i=1}^{i=50}$ is computed.
\begin{table}[tb]
\renewcommand{\arraystretch}{1.25}
\begin{center}
\begin{tabular}{rrlrrr}
\hline \hline
M & m & sto &  threads & MIPRELSTOP & duration \\
\hline
5 & 2 & \checkmark&  5 & 6 & \textbf{3494} \\
5 & 5  &  \checkmark&  10 & 5 & 233 \\
5 & 10  & \checkmark &  5 & 5 & \textbf{1134} \\
5 & 25  &  \checkmark&  5 & \textbf{10} & 133 \\
10 & 2  &  \checkmark &  5 & \textbf{10} & 480 \\
10 & 5  &  $\times$&  5 &  20 & $>$ \textbf{5760} \\
10 & 10  &  $\times$&  5 &  20 & $>$ \textbf{5760} \\
10 & 25  & \checkmark &  20 & \textbf{10} & \textbf{2231} \\
20 & 2 & \checkmark &  20&  \textbf{10} & \textbf{2191} \\
20 & 5  & $\times$ & 5  & 20 & $>$ \textbf{5760} \\
20 & 10  & $\times$&  5  & 20 & $>$ \textbf{5760} \\
20 & 25  & \checkmark & 5  & 5 & \textbf{3196} \\
50 & 2  &  $\times$&  10 & 65 & $>$ \textbf{5760} \\
50 & 5  &  $\times$&  - & - & - \\
50 & 10  &  $\times$&  -& - & -\\
50 & 25  & $\times$&   -& - & -\\
\hline \hline
M & m & $\text{sto}^\star$  & threads & MIPRELSTOP & duration \\ \hline
5 & 2 & \checkmark & 5 & 1 & 235 \\
5 & 5  &   \checkmark & 5 & 1  & 25  \\
5 & 10  &  \checkmark & 5 & 1 & 5 \\
5 & 25  &   \checkmark & 5 & 7 & \textbf{1499} \\
10 & 2  &   \checkmark & 5 & 1.5 & 88  \\
10 & 5  &   \checkmark & 20 & 1 & 119 \\
10 & 10  &   \checkmark & 5 & 5 & 187 \\
10 & 25  &  \checkmark & 20 & \textbf{10} & 351 \\
20 & 2 & \checkmark & 10 & 1.5 & 758 \\
20 & 5  & \checkmark & 20 & 1 & 186 \\
20 & 10  &   \checkmark & 5 & \textbf{10} & \textbf{2125} \\
20 & 25  &  \checkmark & 20 & \textbf{10} & \textbf{2231} \\
50 & 2  &  \checkmark & 5 & 5 & 24 \\
50 & 5  & \checkmark & 20 & 1 & 378 \\
50 & 10  &  $\times$ & 20 & - & -\\
50 & 25  &   \checkmark & 5 & \textbf{10} & 83 \\\hline \hline
\end{tabular}
\caption{$\text{ON/OFF}(\text{opt}, M, m)$ plans computed using the two-stage stochastic optimizer of the single-phase framework. 
$M$ is the number of scenarios randomly selected and used to solve the two-stage stochastic problem, $m$ (\%) is the maximal deviation of scenarios for a three-sigma confidence interval, $\text{opt}=\text{sto}$ is the two-stage stochastic optimizer in the single-framework, $\text{opt}=\text{sto}^\star$ is the version with the binary variables of the second-stage relaxed between 0 and 1, threads is the number of threads used, MIPRELSTOP is the convergence criterion (\%) which determines when the branch and bound tree search will terminate, and duration the total computation time (minutes). 
When the computation time exceeded four days (5760 minutes) and did not reach at least a MIPRELSTOP of 10 \%, no solution was retrieved.}
\label{tab:on-off-computations}
\end{center}
\end{table}
\begin{figure}[tb]
\begin{subfigure}{0.5\textwidth}
		\centering
		\includegraphics[width=\linewidth]{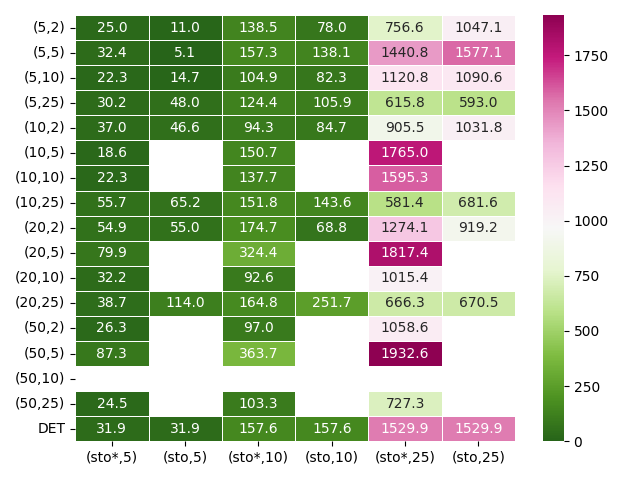}
  \caption{Mean lost load (MWh).}
	\end{subfigure}
 	\begin{subfigure}{0.5\textwidth}
		\centering
		\includegraphics[width=\linewidth]{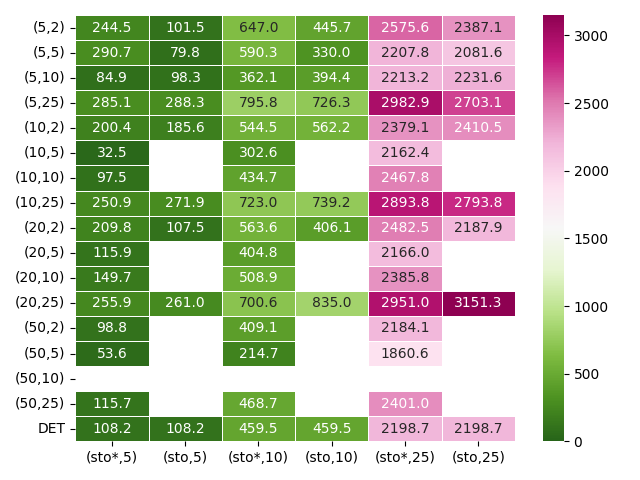}
  \caption{Mean lost production (MWh).}
	\end{subfigure}
	\caption{Average values over the 50 evaluation scenarios of lost load/production $\sum_{i=1}^{50}\frac{1}{50} \text{KPI}^i((\text{opt}, M, m), m_\text{eval})$ using the stochastic optimizer of the single-phase framework. The y-axis provides the pair $(M,m)$ used to compute the $\text{ON/OFF}(\text{opt}, M, m)$ plan, and the x-axis indicates the $m_\text{eval}$ parameter of the evaluation scenario and the optimizer used.}
	\label{fig:single-framework-heatmap-mean-values}
\end{figure}
\begin{figure}[tb]
\begin{subfigure}{0.5\textwidth}
		\centering
		\includegraphics[width=\linewidth]{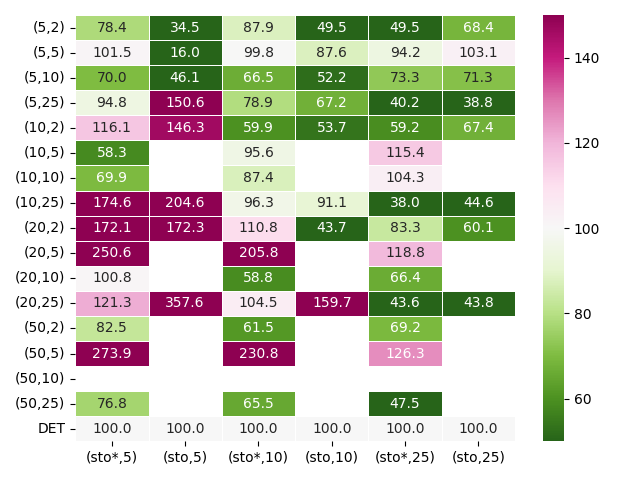}
  \caption{Mean lost load (\%).}
	\end{subfigure}
 	\begin{subfigure}{0.5\textwidth}
		\centering
		\includegraphics[width=\linewidth]{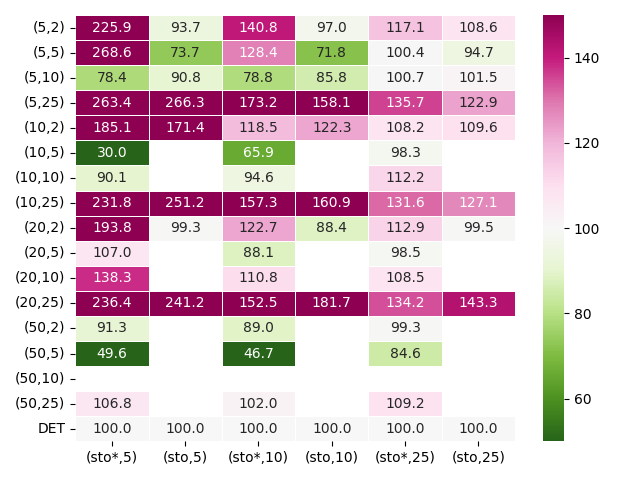}
  \caption{Mean lost production (\%).}
	\end{subfigure}
	\caption{Average values over the 50 evaluation scenarios of lost load/production divided by the deterministic solution $\sum_{i=1}^{50} \frac{1}{50}\frac{\text{KPI}^i((\text{opt}, M, m), m_\text{eval})}{\text{KPI}^i(\text{det}, m_\text{eval})}$ using the stochastic optimizer of the single-phase framework. The y-axis provides the pair $(M,m)$ used to compute the $\text{ON/OFF}(\text{opt}, M, m)$ plan, and the x-axis indicates the $m_\text{eval}$ parameter of the evaluation scenario and the optimizer used. A value higher than 100 \% means an average volume of lost load or lost production higher than the deterministic approach.}
	\label{fig:single-framework-heatmap-mean-values-perc}
\end{figure}
\begin{figure}[tb]
\begin{subfigure}{0.5\textwidth}
		\centering
		\includegraphics[width=\linewidth]{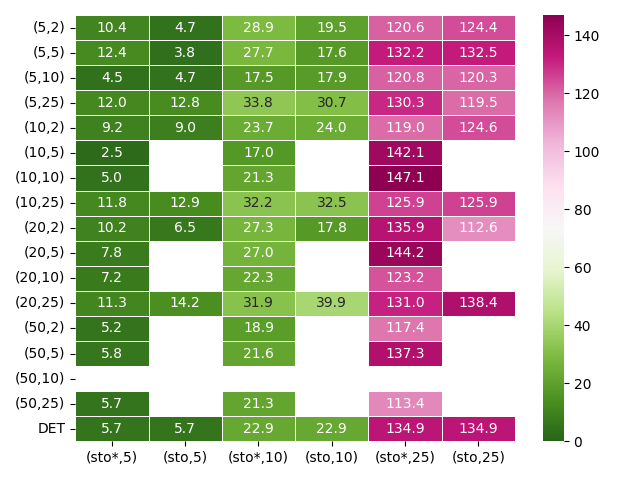}
  \caption{Mean dispatch costs (k\euro).}
	\end{subfigure}
 	\begin{subfigure}{0.5\textwidth}
		\centering
		\includegraphics[width=\linewidth]{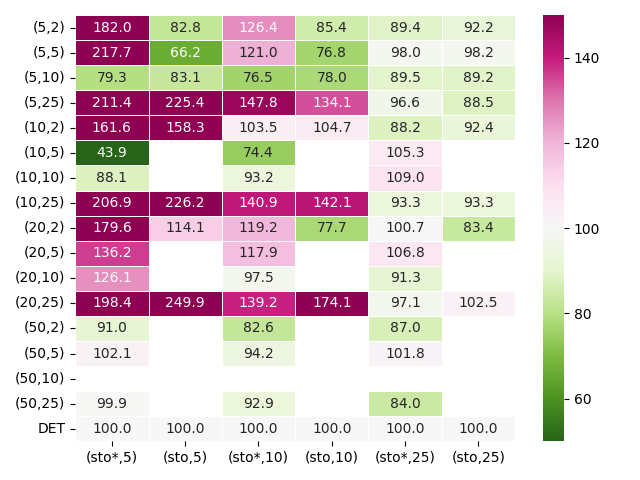}
  \caption{Mean dispatch costs (\%).}
	\end{subfigure}
	\caption{Average dispatch costs (k\euro) and divided by the deterministic solution (\%). The y-axis provides the pair $(M,m)$ used to compute the $\text{ON/OFF}(\text{opt}, M, m)$ plan, and the x-axis indicates the $m_\text{eval}$ parameter of the evaluation scenario and the optimizer used.}
	\label{fig:single-framework-heatmap-mean-values-obj}
\end{figure}

%%%%%%%%%%%%%%%%%%%%%%%%%%%%%%%
% Presentation of the heatmaps 
%%%%%%%%%%%%%%%%%%%%%%%%%%%%%%%%
Figure \ref{fig:single-framework-heatmap-mean-values} depicts heatmaps of the average values over the 50 evaluation scenarios of lost load/production: $\sum_{i=1}^{50}\frac{1}{50} \text{KPI}^i((\text{opt}, M, m), m_\text{eval})$.
Figure \ref{fig:single-framework-heatmap-mean-values-perc} depicts heatmaps of the average values over the 50 evaluation scenarios divided by the deterministic solution: $\sum_{i=1}^{50} \frac{1}{50}\frac{\text{KPI}^i((\text{opt}, M, m), m_\text{eval})}{\text{KPI}^i(\text{det}, m_\text{eval})}$. 
Figure \ref{fig:single-framework-heatmap-mean-values-obj} depicts heatmaps of the mean dispatch costs (k\euro) and divided by the deterministic solution (\%).

%%%%%%%%%%%%%%%%%%%%%%%%%%%%%%%%%%%%%%%%%%%%%%%%%%%%%%%%
% Interpretations lost load/production
%%%%%%%%%%%%%%%%%%%%%%%%%%%%%%%%%%%%%%%%%%%%%%%%%%%%%%%%
Notice that the blanks in the heatmaps refer to ON/OFF plans that did not reach the MIPRELSTOP criterion within four days (see Table \ref{tab:on-off-computations}).
First, the stochastic optimizer (sto) and its version with relaxed second-stage binary variables ($\text{sto}^\star$) tend to compute ON/OFF plans, which generate approximately as much lost load and production during the evaluation. Indeed, for some pairs $(M, m)$, sto produces less than $\text{sto}^\star$; for others, it is the opposite.
Second, the stochastic optimizer (sto) tends to compute ON/OFF plans, which generate less lost load than the deterministic optimizer (det), except for a few plans with $(M > 10, m = 5\%)$. However, the deterministic optimizer tends to compute ON/OFF plans, which generate less lost production than the stochastic optimizers (sto and $\text{sto}^\star$).
Third, the more $m_\text{eval}$ increases, the more the lost load and lost production increase for all optimizers and all pairs $(M, m)$. This was expected because greater values of $m_\text{eval}$ imply more extreme scenarios.
Fourth, with $M$ and $m_\text{eval}$ fixed, the lost load does not tend to decrease when $m$ increases with the stochastic optimizer. And for the lost production, it even manages to increase slightly. This is counter-intuitive as we could have expected that with higher values $m$, the ON/OFF plans should be more robust.
Finally, with $m$ and $m_\text{eval}$ fixed, there is no trend when $M$ increases with the stochastic optimizer for the lost load and production.

%%%%%%%%%%%%%%%%%%%%%%%%%%%%%%%%%%%%%%%%%%%%%%%%%%%%%%%%
% Interpretation of the dispatch costs
%%%%%%%%%%%%%%%%%%%%%%%%%%%%%%%%%%%%%%%%%%%%%%%%%%%%%%%%
Regarding the dispatch costs, the results are similar but slightly different from the lost load and production.
First, the stochastic optimizer tends to compute ON/OFF plans, implying less dispatch costs in the evaluation phase than $\text{sto}^\star$, and the deterministic optimizer, particularly when $m_\text{eval}$ increases. 
Second, the more $m_\text{eval}$ increases, the more the dispatch costs increase for all optimizers and all pairs $(M, m)$. This was expected because greater values of $m_\text{eval}$ imply more extreme scenarios.
Third, with $M$ and $m_\text{eval}$ fixed, the dispatch costs tend to slightly increase when $m$ increases with the stochastic optimizer. 
Finally, similarly, with $m$ and $m_\text{eval}$ fixed, the dispatch costs tend to slightly increase when $M$ increases with the stochastic optimizer.
The two last observations are not what we would like to expect as a more robust ON/OFF plan (by increasing $m$ or $M$) should imply less dispatch costs in the evaluation phase regarding extreme scenarios. 

\subsection{Results using scenario selection}

\begin{figure}[tb]
\centerline{\includegraphics[width=\linewidth]{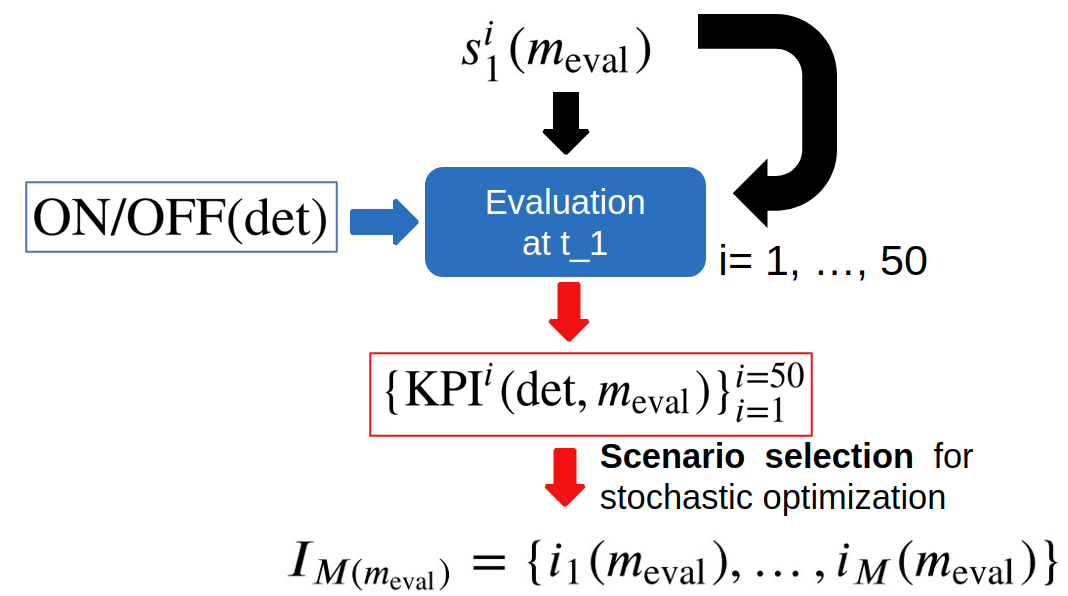}}
\caption{Overview of the importance sampling strategy to select the set $I_{m_\text{eval}} = \{i_1(m_\text{eval}), ..., i_{M(m_\text{eval})(m_\text{eval})}\}$ composed of $M(m_\text{eval})$ worst scenarios to be used as input of the stochastic optimizer.}
\label{fig:importance-sampling-strategy}
\end{figure}

%%%%%%
% Overview of the scenario selection approach
%%%%
In contrast to the previous section, we compute ON/OFF plans by selecting specific worst scenarios instead of random ones. 
Figure \ref{fig:importance-sampling-strategy} depicts the overall importance sampling process to select $M$ scenarios to compute the ON/OFF plan using the stochastic optimizer.
The ON/OFF plan of the deterministic optimizer is evaluated at $t_1$ using the ex-post \textit{out-of-sample} evaluation strategy with scenarios $s_1^i(m_\text{eval})$.
It allows the computation $\text{KPI}^i(\text{det}, m_\text{eval})$, with KPI the volume (MWh) of lost load and production over the period of study $[T_1, T_2]$.
Indeed, we solve a deterministic optimization using as inputs $\text{ON/OFF}(\text{det})$ and $s_1^i(m_\text{eval})$ over $[T_1, T_3]$ and retrieve the lost load and production over $[T_1, T_2]$.
By applying this process recursively $\forall i \in [1, K]$ the set of KPIs $\{ \text{KPI}^i(\text{det}, m_\text{eval})\}_{i=1}^{K}$ is computed.

This set is sorted out by increasing volumes of lost load and production, providing the scenarios $i \in [1, K]$ that cause these highest values. It allows to build two sets:
i) $I^\text{lost load} = \{i_1(m_\text{eval}), \ldots, i_K(m_\text{eval})\}$, with $i_1(m_\text{eval})$ the evaluation scenario which induces the highest value of lost load over all scenarios $i \in [1, K]$ with $m_\text{eval}$ parameter;
ii) $I^\text{lost production} = \{j_1(m_\text{eval}), \ldots, j_K(m_\text{eval})\}$, with $j_1(m_\text{eval})$ the evaluation scenario which induces the highest value of lost load over all scenarios $j \in [1, K]$ with $m_\text{eval}$ parameter.
Then, we select the first $M(m_\text{eval}) \in I^\text{lost load} = I^\text{lost load}_{M(m_\text{eval})} = \{i_1(m_\text{eval}), \ldots, i_{M(m_\text{eval})(m_\text{eval})}\}$ and $M(m_\text{eval}) \in I^\text{lost production} = I^\text{lost production}_{M(m_\text{eval})} = \{i_1(m_\text{eval}), \ldots, i_{M(m_\text{eval})(m_\text{eval})}\}$ worst scenarios which are used for the stochastic optimization.
Finally, the set of selected scenarios is $I_{M(m_\text{eval})} = I^\text{lost load}_{M(m_\text{eval})} \cup I^\text{lost production}_{M(m_\text{eval})}$.
Thus, it is possible to select several values of $M(m_\text{eval})$ for a fixed $m_\text{eval}$. The higher the value $M(m_\text{eval})$, the more we anticipate that the stochastic optimizer should be robust to extreme scenarios.

In practice, we select $M(m_\text{eval})=2, 5, 10$, and $m_\text{eval}=5, 10, 25 \%$. It results in 9 possible values of $M(m_\text{eval})$ and in 9 sets $I_{M(m_\text{eval})}$ of worst scenarios leading to 9 $\text{ON/OFF}(\text{opt}, I_{M(m_\text{eval})}, m)$ plans for a given stochastic optimizer.
\begin{figure}[tb]
\centerline{\includegraphics[width=\linewidth]{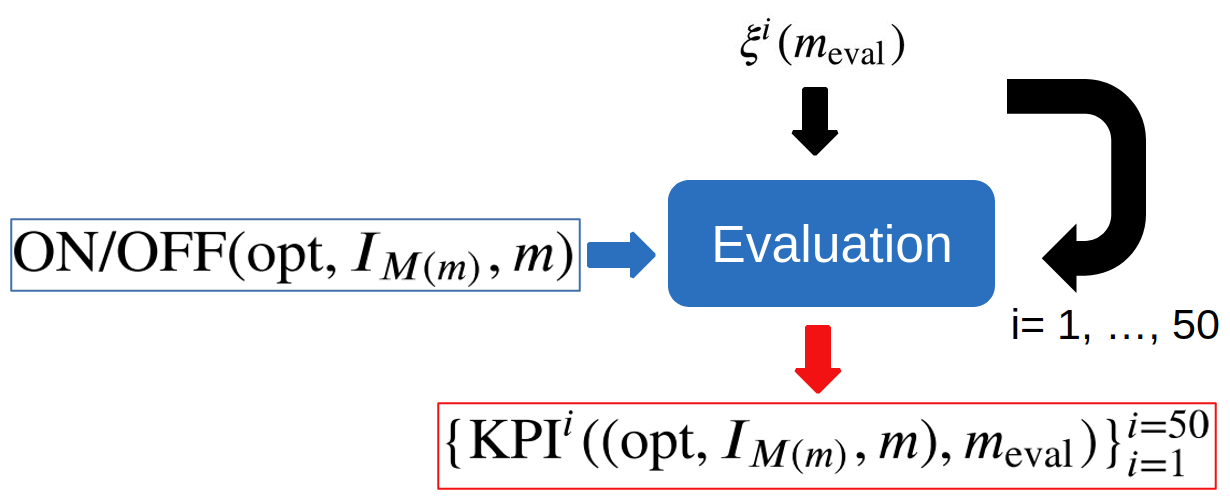}}
\caption{Overview of the ex-post \textit{out-of-sample} evaluation for a given $\text{ON/OFF}(\text{opt}, I_{M(m)}, m)$ plan with $I_{M(m)}$ the set of $M(m)$ worst selected scenarios.}
\label{fig:strategie-evaluation-with-importance}
\end{figure}

Figure \ref{fig:strategie-evaluation-with-importance} depicts the ex-post \textit{out-of-sample} evaluation for a given $\text{ON/OFF}(\text{opt}, I_{M(m)}, m)$ plan with $I_{M(m)}$ the set of $M(m)$ worst selected scenarios.
%
% Interpretation computation and comparison without scenario selection
%
Table \ref{tab:on-off-with-importance-computations} presents the computation time for each $\text{ON/OFF}(\text{opt}, I_{M(m)}, m)$ plan given the number of threads used and the MIPRELSTOP criterion fixed. The crosses indicate the optimization did not converge within the time allowed for a given MIPRELSTOP criterion.
The scenario selection increases the number of computations that reach the convergence criterion compared to Table \ref{tab:on-off-computations}. Only one computation $M(m)=10, m = 10\%$ did not reach at least MIPRELSTOP $= 20\%$ within four days. In addition, as the computation time decreases, smaller values of MIPRELSTOP can be reached. 
\begin{table}[tb]
\renewcommand{\arraystretch}{1.25}
\begin{center}
\begin{tabular}{rrlrrr}
\hline \hline
$M(m)$ & m & sto &  threads & MIPRELSTOP & duration \\
\hline
2 & 5  &  \checkmark   &  20 & 7 & 368 \\
5 & 10  & \checkmark   &  20 & 5 & \textbf{1029} \\
10 & 25  &  \checkmark &  10 & 1 & 358 \\
2 & 5  &  \checkmark   &  20 &  \textbf{11} &  \textbf{3728} \\
5 & 10  &  \checkmark  &  10 &  \textbf{10} &  \textbf{2780} \\
10 & 25  & \checkmark  &  10 & \textbf{10} & 211 \\
2 & 5  & \checkmark    &  20 & \textbf{10} & \textbf{4398} \\
5 & 10  & $\times$     &  10 & \textbf{20} & $>$ \textbf{5760} \\
10 & 25  & \checkmark  &  10 & \textbf{20} & \textbf{3322} \\\hline \hline
$M(m)$ & m & $\text{sto}^\star$  & threads & MIPRELSTOP & duration \\ \hline
2 & 5  &   \checkmark  & 10 & 1  & \textbf{2147}  \\
5 & 10  &  \checkmark  & 10 & 1 & 8 \\
10 & 25  &  \checkmark & 10 & 1 & 23 \\
2 & 5  &   \checkmark  & 10 & 1 & \textbf{2713} \\
5 & 10  &   \checkmark & 10 & 5 & 798 \\
10 & 25  &  \checkmark & 10 & 2 & \textbf{1919} \\
2 & 5  & \checkmark    & 10 & 1 & \textbf{4177} \\
5 & 10  &   \checkmark & 15 & 3 & \textbf{5237} \\
10 & 25  &  \checkmark & 20 & \textbf{10} & \textbf{2239} \\\hline \hline
\end{tabular}
\caption{$\text{ON/OFF}(\text{opt}, I_{M(m)}, m)$ plans computed using the two-stage stochastic optimizer of the single-phase framework. 
$I_{M(m)}$ is the set of worst scenarios selected in a previous phase used to solve the two-stage stochastic problem, $m$ (\%) is the maximal deviation of scenarios for a three-sigma confidence interval, $\text{opt}=\text{sto}$ is the two-stage stochastic optimizer in the single-framework, $\text{opt}=\text{sto}^\star$ is the version with the binary variables of the second-stage relaxed between 0 and 1, threads is the number of threads used, MIPRELSTOP is the convergence criterion (\%) which determines when the branch and bound tree search will terminate, and duration the total computation time (minutes). 
When the computation time exceeded four days (5760 minutes) and did not reach at least a MIPRELSTOP of 20 \%, no solution was retrieved.}
\label{tab:on-off-with-importance-computations}
\end{center}
\end{table}
\begin{figure}[tb]
\begin{subfigure}{0.5\textwidth}
		\centering
		\includegraphics[width=\linewidth]{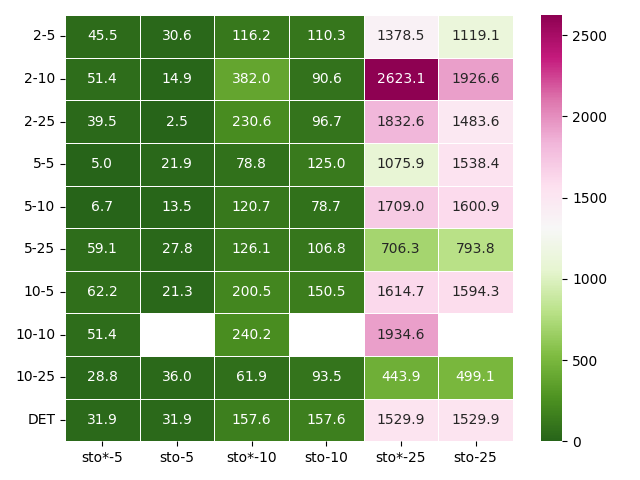}
  \caption{Mean lost load (MWh).}
	\end{subfigure}
 	\begin{subfigure}{0.5\textwidth}
		\centering
		\includegraphics[width=\linewidth]{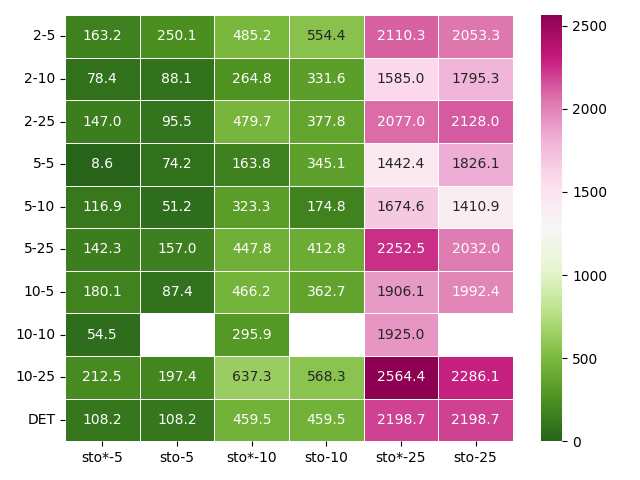}
  \caption{Mean lost production (MWh).}
	\end{subfigure}
	\caption{Average values over the 50 evaluation scenarios of lost load/production using the stochastic optimizer with the set of worst scenarios $I_{M(m)}$. The y-axis provides the pair $(M(m),m)$ used to compute the $\text{ON/OFF}(\text{opt}, I_{M(m)}, m)$ plan, and the x-axis indicates the $m_\text{eval}$ parameter of the evaluation scenario and the optimizer used.}
	\label{fig:sto-sampling-heatmap-mean-values}
\end{figure}
\begin{figure}[tb]
\begin{subfigure}{0.5\textwidth}
		\centering
		\includegraphics[width=\linewidth]{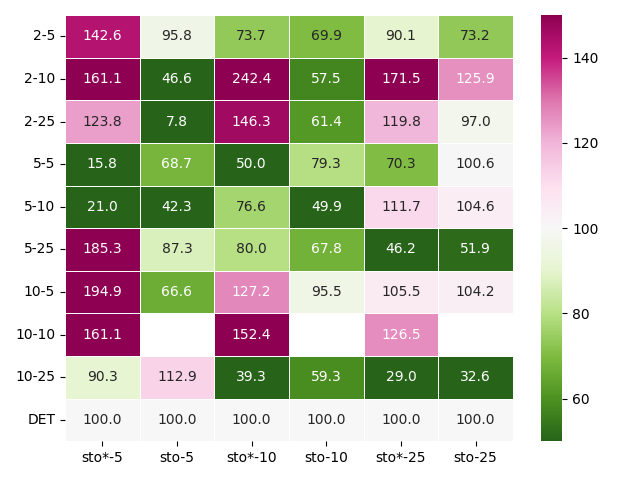}
  \caption{Mean lost load (\%).}
	\end{subfigure}
 	\begin{subfigure}{0.5\textwidth}
		\centering
		\includegraphics[width=\linewidth]{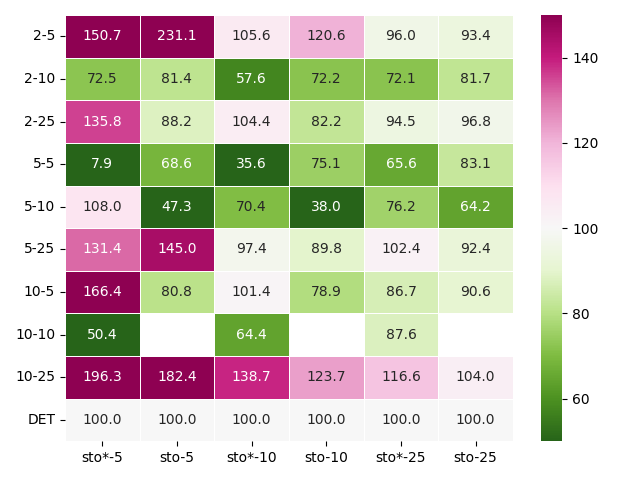}
  \caption{Mean lost production (\%).}
	\end{subfigure}
	\caption{Average values over the 50 evaluation scenarios of lost load/production divided by the deterministic solution using the stochastic optimizer with the set of worst scenarios $I_{M(m)}$. The y-axis provides the pair $(M(m),m)$ used to compute the $\text{ON/OFF}(\text{opt}, I_{M(m)}, m)$ plan, and the x-axis indicates the $m_\text{eval}$ parameter of the evaluation scenario and the optimizer used. }
	\label{fig:sto-sampling-heatmap-mean-perc}
\end{figure}
\begin{figure}[tb]
\begin{subfigure}{0.5\textwidth}
		\centering
		\includegraphics[width=\linewidth]{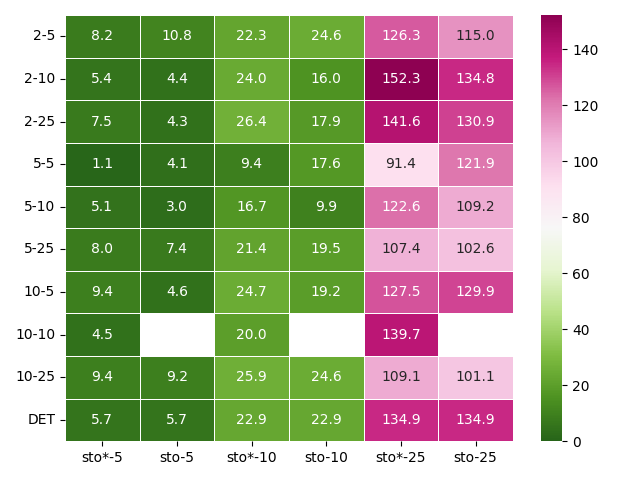}
  \caption{Mean dispatch costs (k\euro).}
	\end{subfigure}
 	\begin{subfigure}{0.5\textwidth}
		\centering
		\includegraphics[width=\linewidth]{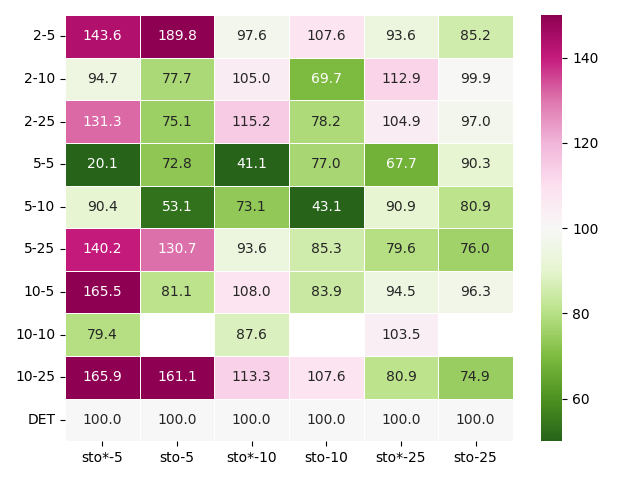}
  \caption{Mean dispatch costs (\%).}
	\end{subfigure}
	\caption{Average dispatch costs (k\euro) using the stochastic optimizer with the set of worst scenarios $I_{M(m)}$ and divided by the deterministic solution (\%). The y-axis provides the pair $(M(m),m)$ used to compute the $\text{ON/OFF}(\text{opt}, I_{M(m)}, m)$ plan, and the x-axis indicates the $m_\text{eval}$ parameter of the evaluation scenario and the optimizer used. }
	\label{fig:sto-sampling-heatmap-mean-obj}
\end{figure}

%%%%%%%%%%%%%%%%%%%%%%%%%%%%%%%
% Presentation of the heatmaps 
%%%%%%%%%%%%%%%%%%%%%%%%%%%%%%%%
Figure \ref{fig:sto-sampling-heatmap-mean-values} depicts heatmaps of the average values over the 50 evaluation scenarios of lost load/production $\sum_{i=1}^{50}\frac{1}{50} \text{KPI}^i((\text{opt}, I_{M(m)}, m), m_\text{eval})$ using the stochastic optimizer with the set of worst scenarios $I_{M(m)}$.
Figure \ref{fig:sto-sampling-heatmap-mean-perc} depicts heatmaps of the average values, using the stochastic optimizer with the set of worst scenarios $I_{M(m)}$, over the 50 evaluation scenarios divided by the deterministic solution $\sum_{i=1}^{50} \frac{1}{50}\frac{\text{KPI}^i((\text{opt}, I_{M(m)}, m), m_\text{eval})}{\text{KPI}^i(\text{det}, m_\text{eval})}$. 
Figure \ref{fig:sto-sampling-heatmap-mean-obj} depicts heatmaps of the mean dispatch costs (k\euro), using the stochastic optimizer with the set of worst scenarios $I_{M(m)}$, which are divided by the deterministic solution (\%).

%%%%%%%%%%%%%%%%%%%%%%%%%%%%%%%%%%%%%%%%%%%%%%%%%%%%%%%%
% Interpretations
%%%%%%%%%%%%%%%%%%%%%%%%%%%%%%%%%%%%%%%%%%%%%%%%%%%%%%%%
Notice that the blanks in the heatmaps refer to ON/OFF plans that did not reach the MIPRELSTOP criterion within four days (see Table \ref{tab:on-off-with-importance-computations}).
First, the stochastic optimizer (sto) generates less lost load than $\text{sto}^\star$. However, there is no clear trend for the lost production.
Second, the stochastic optimizer (sto) tends to compute ON/OFF plans, which generate less lost load than the deterministic optimizer (det), except for the ON/OFF plans with $(N(m)=2, m=10\%)$, $(N(m)=5, m=10\%)$, and $(N(m)=10, m=5\%)$. This trend is similar when considering the lost production with  $m_\text{eval} > 5 \%$.
Third, the more $m_\text{eval}$ increases, the more the lost load and lost production increase for all optimizers and all pairs $(M, m)$. This was expected because greater values of $m_\text{eval}$ imply more extreme scenarios.
Fourth, with $M(m)$ and $m_\text{eval}$ fixed, the lost load tends to decrease when $m$ increases with the stochastic optimizer, particularly with $m_\text{eval}= 25 \%$ except with $M(m)=2, m=5, 10, 25 \%$. This result is consistent as the higher $m$, the more robust the ON/OFF plan computed by the stochastic optimizer to extreme evaluation scenarios. However, when considering the lost production, there is no clear trend.
Finally, with $m=25\%$ and $m_\text{eval}$ fixed, the lost load decreases when $M(m)$ increases with the stochastic optimizer. This result is also consistent as the higher $M(m)$, the more robust the ON/OFF plan computed by the stochastic optimizer to extreme evaluation scenarios. However, no such observation exists for $m=5, 10 \%$. Similarly, when considering the lost production, there is no clear trend.

%%%%%%%%%%%%%%%%%%%%%%%%%%%%%%%%%%%%%%%%%%%%%%%%%%%%%%%%
% Interpretations dispatch costs
%%%%%%%%%%%%%%%%%%%%%%%%%%%%%%%%%%%%%%%%%%%%%%%%%%%%%%%%
Regarding the dispatch costs, the results are fourfold. 
First, the stochastic optimizer tends to compute ON/OFF plans, implying less dispatch costs in the evaluation phase than $\text{sto}^\star$, and the deterministic optimizer, particularly when $m_\text{eval}$ increases. 
Second, the more $m_\text{eval}$ increases, the more the dispatch costs increase for all optimizers and all pairs $(M(m), m)$. This was expected because greater values of $m_\text{eval}$ imply more extreme scenarios.
Third, with $M(m)$ and $m_\text{eval}$ fixed, the dispatch costs tend to decrease slightly when $m$ increases with the stochastic optimizer. 
Finally, similarly, with $m$ and $m_\text{eval}$ fixed, there is no clear trend when $M$ increases with the stochastic optimizer.

Appendix \ref{appendix:comparison-with-without-scenario-selection} compares the mean values of lost load/production and dispatch costs of the stochastic optimizer without and with scenario selection. 
Figures \ref{fig:sto-comparison-val-mean} and \ref{fig:sto-comparison-mean-perc} indicate that overall scenario selection provides better results with these KPIs. 
Figure \ref{fig:on-off-comparison} depicts ON/OFF plans with and without scenario selection and deterministic optimizer. 
The differences between the deterministic and stochastic optimizers are noticeable. Stochastic optimizers tend to maintain more nuclear power plants ON than deterministic ones. In addition, they tend to start CCGT units during the morning and evening consumption peaks. Indeed, they consider scenarios where the residual demand is higher than the mean scenario used by the deterministic optimizer.
However, the differences between the several versions of the stochastic optimizers are more challenging to interpret.
Thus, it illustrates the need for KPIs to evaluate these ON/OFF plans.

\section{Conclusions and perspectives}\label{sec:conclusions}

% reminder of the proposed approach
This study presents a probabilistic unit commitment strategy with a stochastic optimization-based approach, including unit fixed and variable costs and lost load and production costs. 
The problem is by design formulated with a multi-stage stochastic program due to the technical constraints of the conventional units. However, this problem is intractable.
Therefore, it is approximated with a multi-phase framework using a \textit{two-stage stochastic model predictive control}, which consists of solving a sequence of two-stage stochastic programs.
In this study, we decided to use a single two-stage stochastic program and investigate the multi-phase framework in future works.

% Case study and main results
The case study comprises the production units likely to make upward/downward variations on a specific day where the French TSO RTE faced a deficit of downward margins. 
The results indicate that the ON/OFF plans computed by the stochastic optimizer (with or without scenario selection) tend to be more robust to the uncertainty regarding lost load and production than the ON/OFF plan of the deterministic approach.
These results are improved by using a scenario selection. In addition, selecting the worst scenarios allows for a decrease in the time computation and an increase in the number of computations that reach the convergence criterion.
In real-time, in anticipation of a deficit of available margins a few hours ahead, this UC probabilistic could be used to determine which units to activate or not. 
However, there is still much work to achieve to use such a strategy to help ensure supply security. 

% Limitation 1: computation time limitation
The computational burden is the first limitation. Indeed, the stochastic optimizer must deliver results with a computation time in line with the operational requirements, \textit{i.e.}, a few minutes. However, even with the scenario selection implemented, the computation time increases drastically with the number of scenarios, limiting the computation to approximately 10 to 20 scenarios, which takes a few minutes to reach several hours and even a few days.
Therefore, several extensions could be investigated to tackle this issue. 

% research direction to tackle computation time issue
First, use more advanced importance sampling techniques to identify relevant scenarios for stochastic optimization. They could decrease the number of scenarios required to model the uncertainty.
Second, use advanced decomposition techniques such as \textit{accelerated Benders decomposition} approach \cite{zhao2023benefits} and \cite{ramirez2023benders} allowing to speed the computation. Such approaches could be combined with specific techniques to reformulate constraints responsible for the computational burden. 
Third, investigate a hybrid approach combining machine learning and optimization. Such an approach could help to compute quickly a feasible solution. Indeed, optimization proxies \cite{chen2023endtoend, CHEN2022108566} allow replacing a time-consuming optimization model with a machine-learning proxy that can be used in real-time and/or in computationally demanding applications. 

% Limitation 2: pb formulation
The second limitation is the problem formulation, as the single-phase approach provides the results of this study. However, the problem is multi-stage by design due to the technical constraints. This limitation could be addressed with two research directions.
First, benefit from time reduction computation from the previous extensions to implement and test the two-stage stochastic MPC approach and compare it with the single-phase approach investigated in this study. 
Second, better modeling the uncertainty with a scenario tree. Then, use a multi-stage formulation solved with a decomposition technique such as SDDP \cite{pereira1991multi}. This research direction was investigated in a research internship \cite{lucille2023}. However, much work is still required to generate a relevant scenario tree and implement decomposition techniques to tackle the tractability issues.

% Limitation 3: risk aversion
The third limitation lies in the risk aversion that must be well-calibrated to the TSO risk policy. Indeed, other risk-aversions strategies could be implemented than the neutral one, such as with a chance-constraint approach \cite{ROALD2023108725} or by considering the Conditional Value-at-Risk \cite{ROCKAFELLAR20021443}.

% Limitation 4: a too simple case study
Finally, the last limitation is the actual case study, which encloses too many approximations. We should design more realistic case studies by including hydraulic units and refining the temporal resolution from one hour to 15 minutes. In addition, the operational window where the French TSO can activate upward/downward units is expected to decrease from two to one hour. This regulation modification should be considered. 

% Other case studies
In terms of application, this work could serve as a baseline for TSO needs other than the deficit of available margins. 
First, it could help TSOs to formulate optimal bidding on the European reserves markets (TERRE\footnote{\url{https://www.entsoe.eu/network_codes/eb/terre/}}, MARI\footnote{\url{https://www.entsoe.eu/network_codes/eb/mari/}}).
Second, a stochastic formulation could be derived from the proposed one in this study to size operational reserves.

\section*{Acknowledgment}
The authors would like to acknowledge the help of Nathalie Grisey and Sébastien Finet in challenging the formulation and results. 
We also acknowledge Patrick Ponciatici, Lucas Saludjian, and Manuel Ruiz for the fruitful discussions, comments, and support. 
Finally, We acknowledge Quentin Louveaux and Bertrand Cornélusse, assistant professors at Liège University, for giving relevant insights on the formulation.

\bibliographystyle{IEEEtran}
\bibliography{biblio}

%\newpage
\section{Appendix: notation for the optimization problems}

\subsection*{\textbf{Set and indices}}

\begin{supertabular}{lp{0.7\columnwidth}}
Name & Description \\
\hline
$t$ & Time period index. \\
$s$ & Scenario index. \\
$u$ & Unit index. \\
$e$ & Unit state index. \\
$T$ & Number of periods per day. \\
$\#\Omega_S$ & Number of scenarios. \\
$\mathcal{T}$ & Set of time periods, $\mathcal{T}= \{1,2, \ldots, T\} = [T_1, T_3]$. \\
$\mathcal{T}^{\text{OFL}}$ & Truncated set of time periods, 
$\mathcal{T}^{\text{OFL}}= [t + \Delta_t, t + \text{T}^{u,\text{Flat}}]$. \\
$\mathcal{T}^{On}_\text{min}$ & Truncated set of time periods, 
$\mathcal{T}^{On}= [t + \Delta_t, t + \text{T}^{u,On}_\text{min}]$. \\
$\mathcal{T}^{On}_\text{max}$ & Truncated set of time periods, 
$\mathcal{T}^{On}_\text{max}= [t + \Delta_t + T^{u,On}_\text{max} + \Delta_t, T]$. \\
$\mathcal{T}^{\text{OFF}}$ & Truncated set of time periods, 
$\mathcal{T}^{\text{OFF}}= [t + \Delta_t, t + \Delta_t + \text{T}^{u,\text{OFF}}_\text{min}]$. \\
$\Omega_S$ & Set of scenarios, $\Omega_S= \{1,2, \ldots, \#\Omega_S\}$. \\
$\Omega_E$ & Set of states of power plants, $\Omega_E= \{\text{OU}, \text{OD}, \text{OFL}, \text{OFF}\}$. \\
$\Omega_\text{N}$ & Set of nuclear power plants. \\
$\Omega_\text{FT}$ & Set of first-stage units. \\
$\Omega_\text{ST}$ & Set of second-stage units. \\
$\Omega_\text{Nr}$ & Set of renewable units. \\
$\Omega_\text{T}$ & Set of thermal units, $\Omega_\text{T} = \Omega_\text{FT} \cup\Omega_\text{ST}$. \\
\end{supertabular}

\subsection*{\textbf{Parameters}}

\begin{supertabular}{l p{0.65\columnwidth} }
Name & Description  \\
\hline
$\pi^\text{dns}$ & Cost of lost load [k\euro/MWh].   \\
$\pi^\text{spill}$ & Cost of lost-production (spillage) [k\euro/MWh].   \\
$\pi_f^u$ & Fixed start-up cost of unit $u$ [k\euro].   \\
$\pi_v^u$ & Variable cost of unit $u$ [\euro/MWh].   \\
$\Delta_t$ & Duration of a period [hour].   \\
$\hat{d}_{t,s}^{t'}$ & Demand at $t$ in scenario $s$ issued for $\text{LTTD}_{t'}$ [MW].   \\
$\hat{p}_{t,s}^{t',u}$ & Renewable production  at $t$ for unit $u$ in scenario $s$ issued for $\text{LTTD}_{t'}$ [MW].   \\
$\tilde{d}_{t,s}^{t'}$ & Residual demand at $t$ (demand minus the total renewable generation) in scenario $s$ for $\text{LTTD}_{t'}$ [MW].   \\
$P_\text{max}^u$ & Maximal production of thermal unit $u$ [MW].   \\
$P_\text{min}^u$ & Minimal production of thermal unit $u$ [MW].   \\
$\text{T}_{min}^{u,\text{ON}}$ & Minimal time ON of unit $u$ when started [hours].   \\
$\text{T}_{max}^{u,\text{ON}}$ & Maximal time ON of unit $u$ when started [hours].   \\
$\text{T}_{min}^{u,\text{OFF}}$ & Minimal time OFF of unit $u$ when shut-down [hours].   \\
$\text{T}^{u,\text{FLAT}}$ & Minimal time FLAT of unit $u$ after an OU or OD states [hours].   \\
$\Delta T^{\text{ON}}_{u,min}$ & Minimal time delay before starting a unit $u$ [hours].   \\
$\Delta T^{\text{OFF}}_{u,min}$ & Minimal time delay before shutting down a unit $u$ [hours].   \\
\end{supertabular}

\subsection*{\textbf{Variables}}
\noindent Omit index $s$ when a variable is defined as a first-stage variable. 
\begin{supertabular}{l p{0.8\columnwidth}}
Name & Description \\
\hline
$p_{u, t, s}$ & Production of unit $u$ in scenario $s$ [MW]. \\
$\text{dns}_{t,s}$ & Demand not-served (lost load) in scenario $s$ [MW]. \\
$\text{spill}_{t,s}$ & Lost production (spillage) in scenario $s$ [MW]. \\
$E_{t, s, e}^u$ & Unit $u$ in scenario $s$ is in state $e$ [-]. \\
$T_{t,s, e_i,e_f}^u $  & State transition of unit $u$ in scenario $s$ at $t$ from $e_i$ at $t$ to $e_f$ at $t+1$ [-].  \\
$T_{s, e_i,e_f}^{0,u} $  & State transition of unit $u$ in scenario $s$ at $t_0$ from $e_i$ to $e_f$ at $t$ [-].  \\
\end{supertabular}

\section{Appendix: scenario generation methodology}\label{appendix:scenario-generation}

This appendix presents the methodology to generate the set of unbiased wind power, PV, and consumption scenarios for several time periods $t_1$, $t_2$, $t_{3i}$ used to compute the two-stage stochastic programs, and for the evaluation phase.
The goal is to define an ideal unbiased predictor with a fixed variance overall lead times. 
In this appendix, let $t$ be the current time index ($t_1$, $t_2$, $t_{3i}$), $k$ be the lead time of the prediction, $K$ be the maximum lead time of the prediction, $y_{t+k}$ be the true value of the signal $y$ at time $t + k$, and $\widehat{y}_{t+k|t}$ be the value of $y_{t+k}$ predicted at time $t$. 

\begin{figure}[tb]
\begin{subfigure}{0.25\textwidth}
		\centering
		\includegraphics[width=\linewidth]{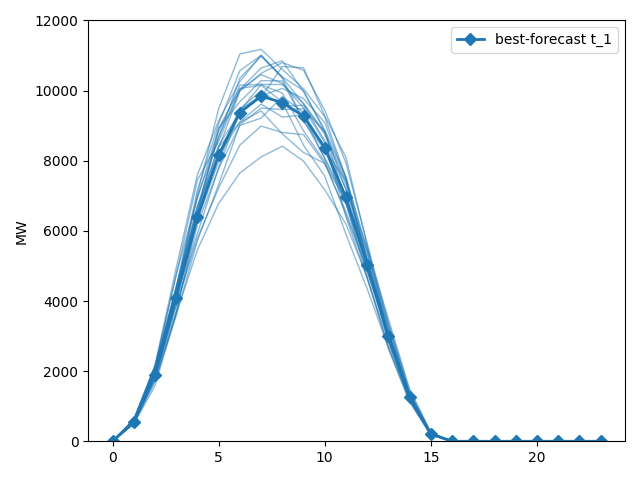}
  \caption{PV scenarios $t_1$.}
	\end{subfigure}%
 	\begin{subfigure}{0.25\textwidth}
		\centering
		\includegraphics[width=\linewidth]{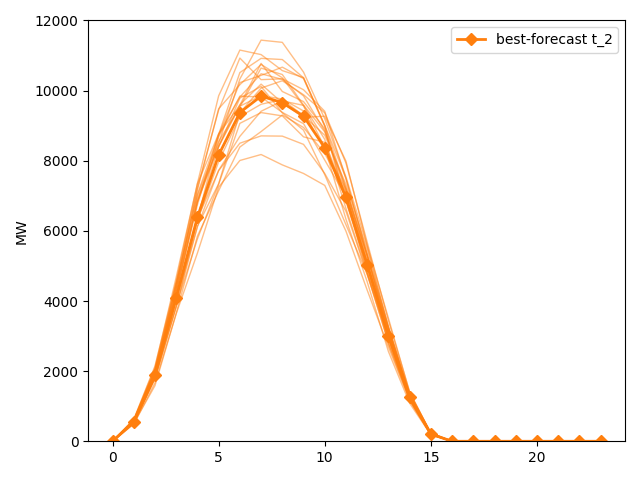}
  \caption{PV scenarios $t_2$.}
	\end{subfigure}
 \begin{subfigure}{0.25\textwidth}
		\centering
		\includegraphics[width=\linewidth]{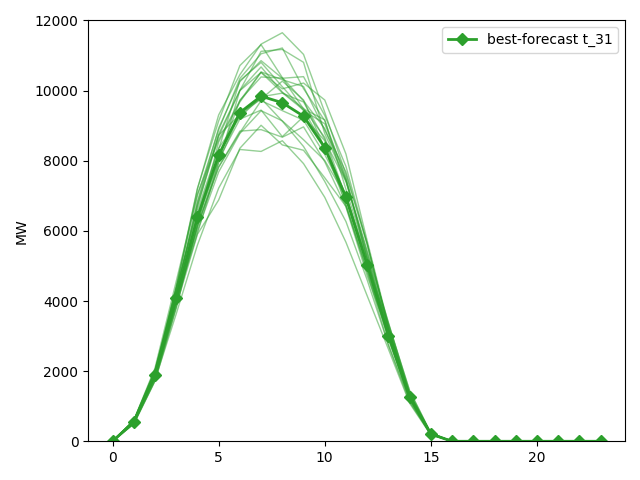}
  \caption{PV scenarios $t_{31}$.}
	\end{subfigure}%
 	\begin{subfigure}{0.25\textwidth}
		\centering
		\includegraphics[width=\linewidth]{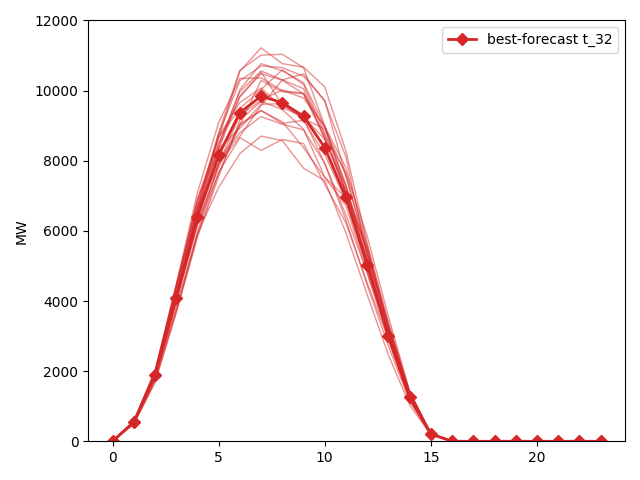}
  \caption{PV scenarios $t_{32}$.}
	\end{subfigure}
 \begin{subfigure}{0.25\textwidth}
		\centering
		\includegraphics[width=\linewidth]{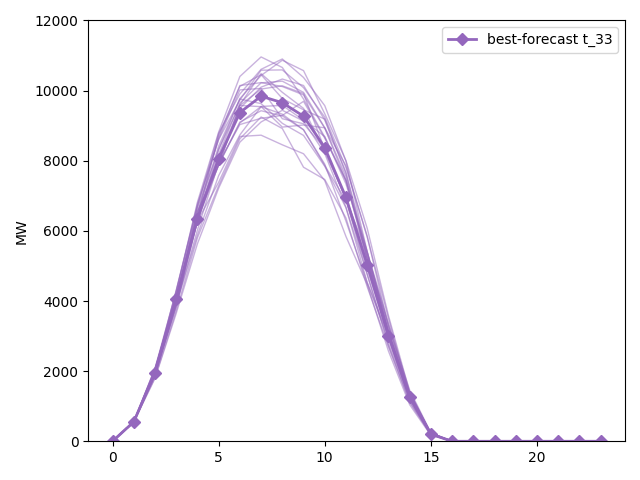}
  \caption{PV scenarios $t_{33}$.}
	\end{subfigure}%
 	\begin{subfigure}{0.25\textwidth}
		\centering
		\includegraphics[width=\linewidth]{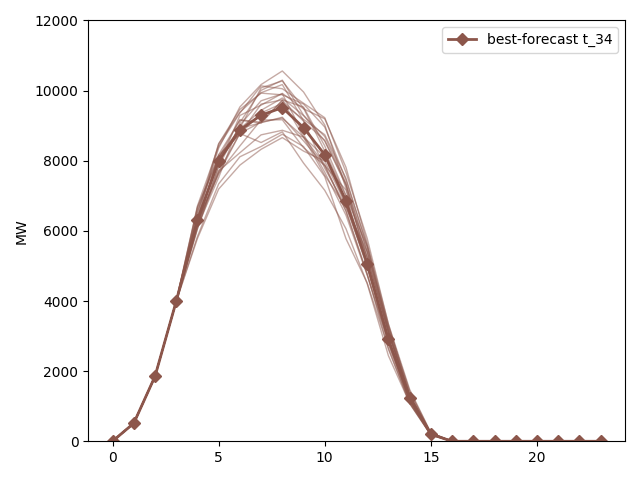}
  \caption{PV scenarios $t_{34}$.}
	\end{subfigure}
 \begin{subfigure}{0.25\textwidth}
		\centering
		\includegraphics[width=\linewidth]{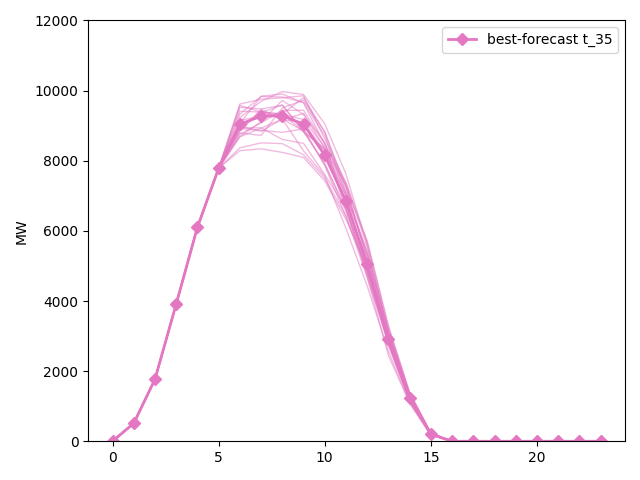}
  \caption{PV scenarios $t_{35}$.}
	\end{subfigure}%
 	\begin{subfigure}{0.25\textwidth}
		\centering
		\includegraphics[width=\linewidth]{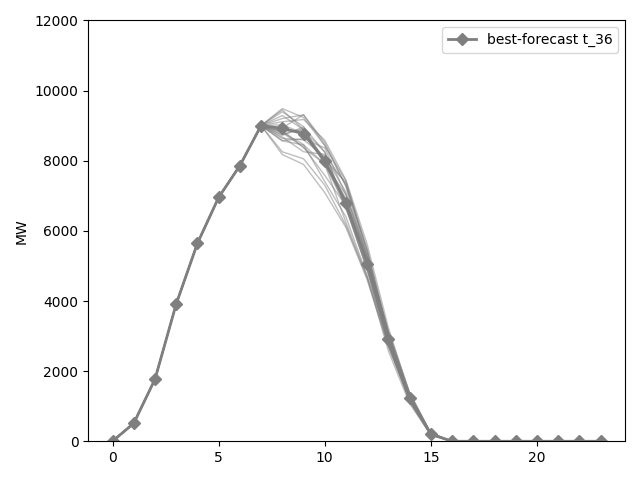}
  \caption{PV scenarios $t_{36}$.}
	\end{subfigure}
	\caption{PV scenarios issued at different periods: $t_1=$ 10 p.m., $t_2=$ 0 a.m., $t_{31}=$ 3 a.m., $t_{32}=$ 5 a.m., $t_{33}=$ 7 a.m., $t_{34}=$ 9 a.m., $t_{35}=$ 11 a.m., $t_{36}=$ 1 p.m. The x-axis starts at 0 a.m., meaning 6 a.m. on day $D$, and ends at 23, meaning 5 a.m. on day $D+1$. Thus, for periods from $t_{33}$ to $t_{36}$, the scenarios cover fewer periods. }
	\label{fig:pv-scenarios}
\end{figure}

For instance, at $t_1 = $ 10 p.m. at day $D-1$, forecasts are computed for $[T_1, T_3]$ with $T_1$ = 6 a.m. at day $D$ and $T_3$  = 6 a.m. at day $D+1$. Then, with $\Delta_t = $ 1 hour, $K = 10 + 24 = 34$. The PV, wind power, and consumption forecasts are needed for lead times from $k=10$ ($T_1$ = 6 a.m. at day $D$) to $k=K=34$ ($T_3$  = 6 a.m. at day $D+1$).
For forecasts computed at $t_{31} = $ 3 a.m. at day $D$ for $[T_1, T_3]$, forecasts are needed for lead times from $k=3$  to $k=K=24+3$.
Then, $\widehat{y}_{t+k|t}$ and $y_{t+k}$ are assumed to be related by
\begin{align}\label{eq:prediction_error_definition}	
\widehat{y}_{t+k|t}  & =  y_{t+k} (1+\epsilon_k).
\end{align}
The error term $\epsilon_k$ is generated by the moving-average model defined in Chapter 3 of \cite{box2015time}
\begin{subequations}
	\begin{align}\label{eq:epsilon_definition}	
	\epsilon_1 & =  \eta_1 \\ 
	\epsilon_k & =  \eta_k + \sum_{i=1}^{k-1} \alpha_i \eta_{k-i} \quad \forall k \in \{2, ..., K\},
	\end{align}
\end{subequations}
with $\{ \alpha_i \}_{i=1}^{K-1} $ scalar coefficients, $\{\eta_k \}_{k=1}^K $ independent and identically distributed sequences of random variables from a normal distribution $\mathcal{N}(0, \sigma)$. Thus, the variance of the error term is
\begin{subequations}
	\begin{align}\label{eq:epsilon_var_definition}	
	& \mathrm{Var} (\epsilon_1)  =  \sigma^2 \\ 
	& \mathrm{Var} (\epsilon_k)  =  \big( 1 + \sum_{i=1}^{k-1} \alpha_i^2 \big) \sigma^2 \quad \forall k \in \{2, ..., K\}.
	\end{align}
\end{subequations}
It is possible to simulate with this model an increase of the prediction error variance with the lead time $k$ by choosing
\begin{align}\label{eq:alpha_definition}	
\alpha_i & =  p^i \quad \forall i \in \{1, ..., K-1\}.
\end{align}
(\ref{eq:epsilon_var_definition}) becomes, $\forall k \in \{1, ..., K\}$
\begin{align}\label{eq:epsilon_definition_2}	
\mathrm{Var} (\epsilon_k) &  = \sigma^2 A_{\epsilon_k},
\end{align}
with $A_{\epsilon_k}$ defined $\forall k \in \{1, ..., K\}$ by
\begin{align}\label{eq:Ak_definition}	
A_{\epsilon_k} &  =  \sum_{i=0}^{k-1} (p^2)^i = \frac{1-(p^2)^{k}}{1-p^2}.
\end{align}
Then, with $0 \leq p < 1 $, it is possible to make the prediction error variance independent of the lead time as it increases. Indeed
\begin{align}\label{eq:Ak_limit}	
\lim_{k \to\infty}  A_{\epsilon_k} & = A_\infty = \frac{1}{1-p^2}.
\end{align}
Finally, the $\sigma$ value to set a maximum $m$ with a high probability of 0.997, corresponding to a three standard deviation confidence interval from a normal distribution, is found by imposing $m = 3 \sqrt{\mathrm{Var} (\epsilon_K)}$
\begin{align}\label{eq:epsilon_max_value}	
\sigma & \approx \frac{m}{3\sqrt{A_\infty}}.
\end{align}
Figure \ref{fig:pv-scenarios} depicts PV scenarios issued at different periods using this approach with $p=0.9$ and a maximal deviation $m=$ 25 \%.

\section{Appendix: comparison with/without scenario selection}\label{appendix:comparison-with-without-scenario-selection}

\begin{figure}[tb]
\begin{subfigure}{0.5\textwidth}
		\centering
		\includegraphics[width=\linewidth]{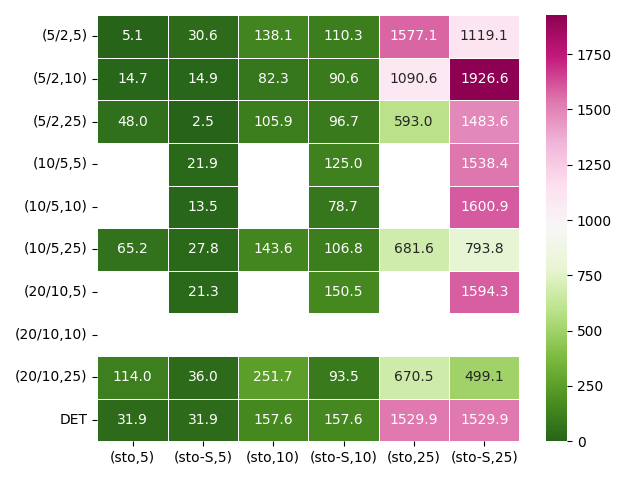}
  \caption{Mean lost load (MWh).}
	\end{subfigure}
 	\begin{subfigure}{0.5\textwidth}
		\centering
		\includegraphics[width=\linewidth]{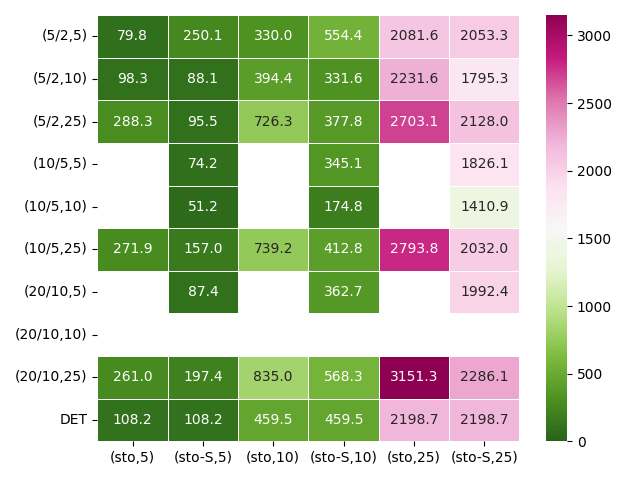}
  \caption{Mean lost production (MWh).}
	\end{subfigure}
 	\begin{subfigure}{0.5\textwidth}
		\centering
		\includegraphics[width=\linewidth]{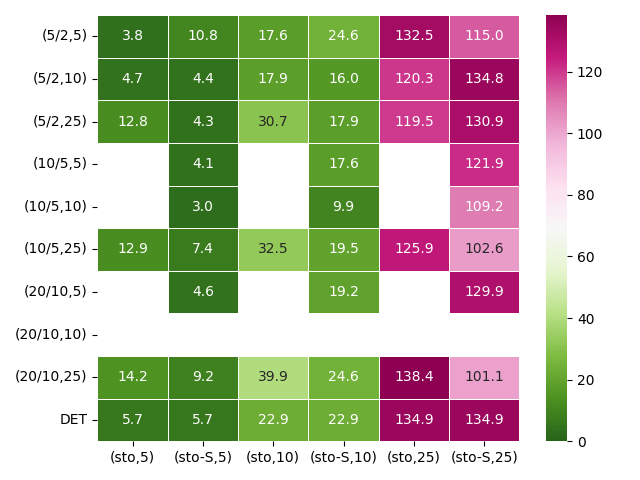}
  \caption{Mean dispatch costs (\euro).}
	\end{subfigure}
	\caption{Comparison of the average values over the 50 evaluation scenarios of lost load/production and dispatch costs using the stochastic optimizer with (sto-S)/without (sto) the set of worst scenarios $I_{M(m)}$. The y-axis provides the pair $(M/M(m),m)$ used to compute the $\text{ON/OFF}(\text{opt}, M/I_{M(m)}, m)$ plan, and the x-axis indicates the $m_\text{eval}$ parameter of the evaluation scenario and the optimizer used.}
	\label{fig:sto-comparison-val-mean}
\end{figure}

\begin{figure}[tb]
\begin{subfigure}{0.5\textwidth}
		\centering
		\includegraphics[width=\linewidth]{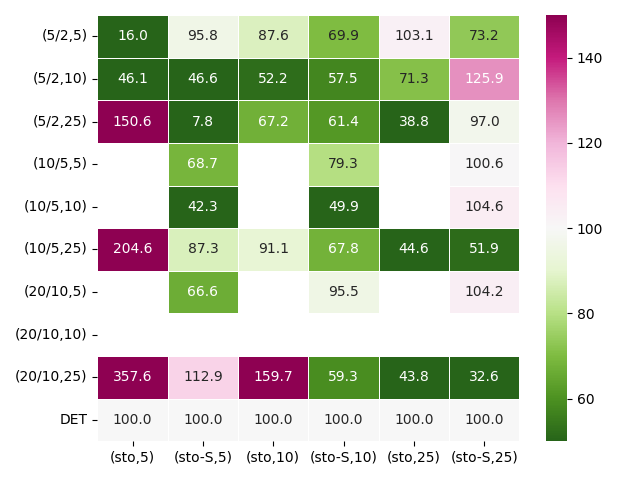}
  \caption{Mean lost load (\%).}
	\end{subfigure}
 	\begin{subfigure}{0.5\textwidth}
		\centering
		\includegraphics[width=\linewidth]{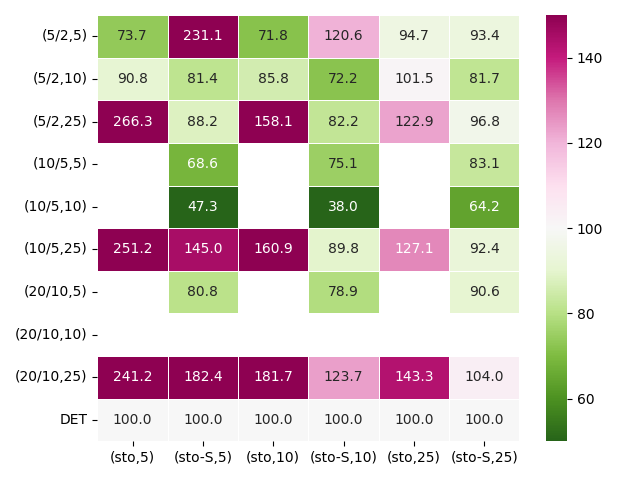}
  \caption{Mean lost production (\%).}
	\end{subfigure}
 	\begin{subfigure}{0.5\textwidth}
		\centering
		\includegraphics[width=\linewidth]{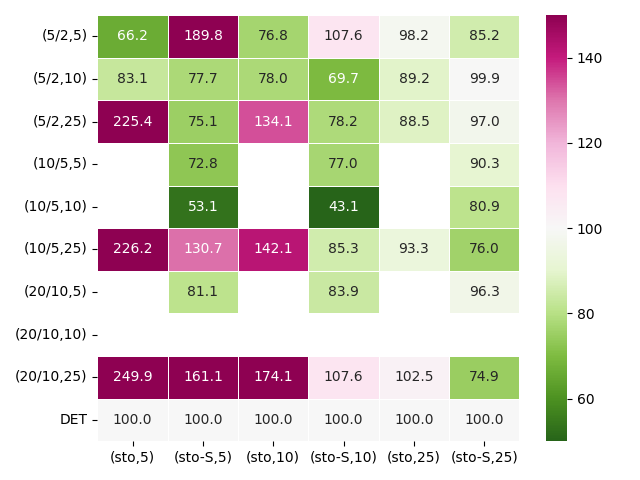}
  \caption{Mean lost production (\%).}
	\end{subfigure}
	\caption{Average values over the 50 evaluation scenarios of lost load/production and dispatch costs normalized by the deterministic solution using the stochastic optimizer with (sto-S)/without (sto) the set of worst scenarios $I_{M(m)}$. The y-axis provides the pair $(M/M(m),m)$ used to compute the $\text{ON/OFF}(\text{opt}, I_{M(m)}, m)$ plan, and the x-axis indicates the $m_\text{eval}$ parameter of the evaluation scenario and the optimizer used.}
	\label{fig:sto-comparison-mean-perc}
\end{figure}

\begin{figure}[tb]
\begin{subfigure}{0.25\textwidth}
		\centering
		\includegraphics[width=\linewidth]{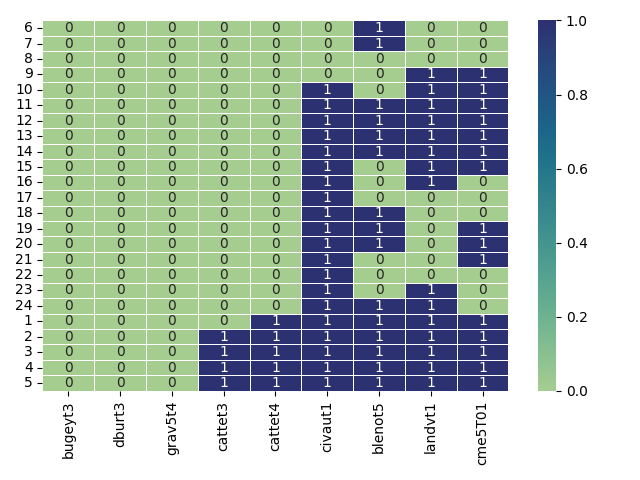}
  \caption{$\text{ON/OFF}(\text{sto}^\star, 10, 25)$.}
	\end{subfigure}%
 	\begin{subfigure}{0.25\textwidth}
		\centering
		\includegraphics[width=\linewidth]{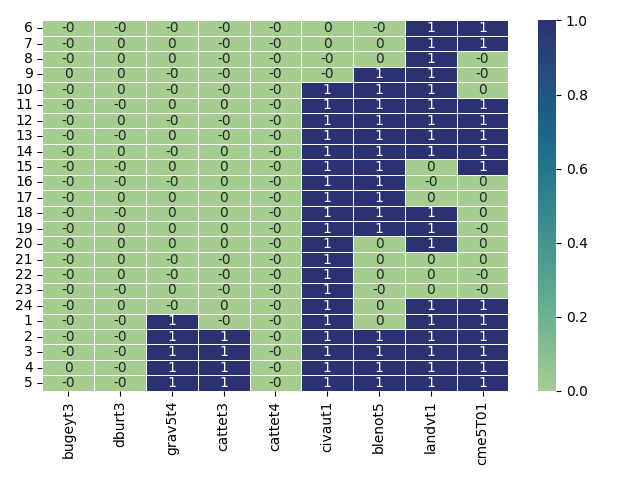}
  \caption{$\text{ON/OFF}(\text{sto}, 10, 25)$.}
	\end{subfigure}
\begin{subfigure}{0.25\textwidth}
		\centering
		\includegraphics[width=\linewidth]{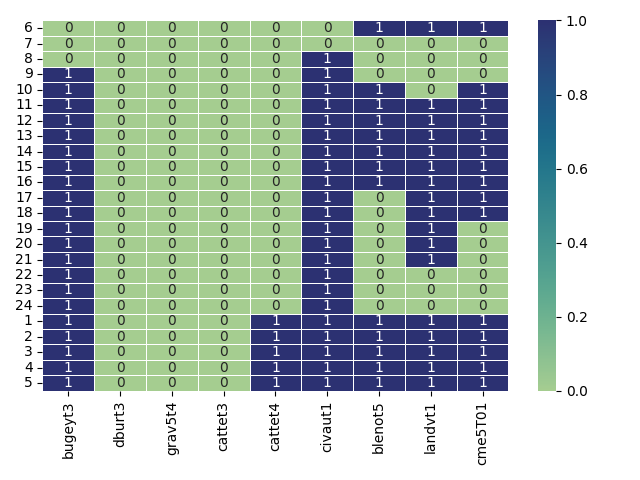}
  \caption{$\text{ON/OFF}(\text{sto}^\star, I_{5}, 25)$.}
	\end{subfigure}%
 	\begin{subfigure}{0.25\textwidth}
		\centering
		\includegraphics[width=\linewidth]{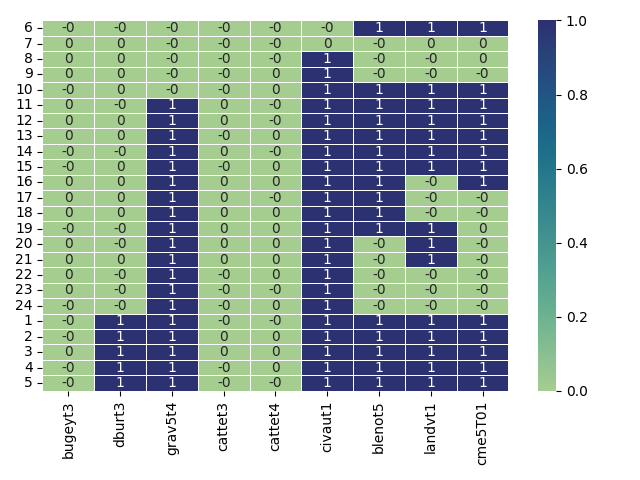}
  \caption{$\text{ON/OFF}(\text{sto}, I_{5}, 25)$.}
	\end{subfigure}
 \centering
  	\begin{subfigure}{0.25\textwidth}
		\centering
		\includegraphics[width=\linewidth]{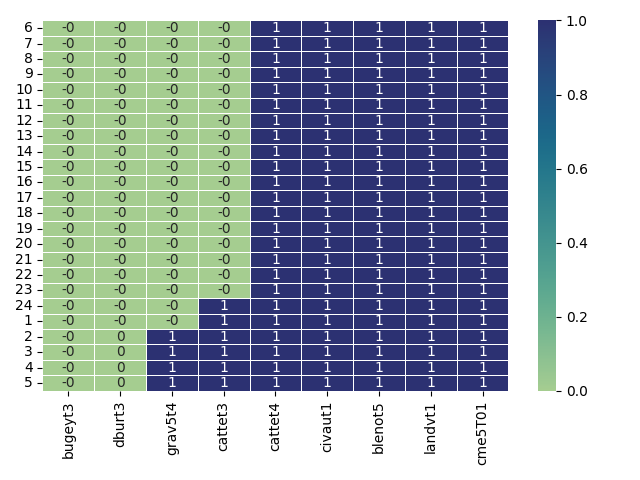}
  \caption{$\text{ON/OFF}(\text{det})$.}
	\end{subfigure}
	\caption{ON/OFF plans comparison with (upper left and right) and without (middle left and right) scenario selection, and deterministic optimizer (lower).
 The production units are listed on the x-axis, with the first six being the nuclear power plants and the last three being the CCGT units. A zero value (green) means the OFF status is 0, and the unit is ON. One value (blue) means the OFF status is 1, and the unit is OFF. 
The y-axis lists the 24 hours considered from 6 a.m. on day D to 5 a.m. on day D+1.}
	\label{fig:on-off-comparison}
\end{figure}

Figures \ref{fig:sto-comparison-val-mean} and \ref{fig:sto-comparison-mean-perc} provide heatmaps with comparison over the 50 evaluation scenarios of the mean values (normalized by the deterministic KPI for the \%) of lost load/production and dispatch costs using the stochastic optimizer with (sto-S)/without (sto) the set of worst scenarios $I_{M(m)}$.
%
% Explications pour faire la correspondence entre sto without scenario selection and with.
%
The $\text{ON/OFF}(\text{opt}, M, m)$ plans of the stochastic optimizer without scenario selection are compared to $\text{ON/OFF}(\text{opt}, I_{M(m)}, m)$ with this equivalence: i) $(M=5, m)$ vs. $I_{M(m)=2}$; ii) $(M=10, m)$ vs. $I_{M(m)=5}$; iii) $(M=20, m)$ vs. $I_{M(m)=10}$. Indeed, $I_{M(m)}$ is the union of two sets composed of the worst scenarios for lost load and lost production. Thus, $I_{M(m)}$ comprises $2 \times M(m)$ scenarios.
%
% Interpretations
These Figures demonstrate that the scenario selection allows, overall, to provide better results in terms of lost load/production and dispatch costs. However, there are some values of $M, m, m_\text{eval}$ where the scenario selection did not improve these KPIs.
%
%%%%%%%%%%%%%%%%%%%%%%%%%%%%%%%%%%%%%%%%%%%%%%%%%%%%%
% ON/OFF examples
%%%%%%%%%%%%%%%%%%%%%%%%%%%%%%%%%%%%%%%%%%%%%%%%%%%%
%
Figure \ref{fig:on-off-comparison} depicts ON/OFF plans with and without scenario selection and deterministic optimizer. The production units are listed on the x-axis, with the first six being the nuclear power plants and the last three being the CCGT units. A zero value (green) means the OFF status is 0, and the unit is ON. One value (blue) means the OFF status is 1, and the unit is OFF. 
The y-axis lists the 24 hours considered from 6 a.m. on day D to 5 a.m. on day D+1.

Interestingly, the deterministic optimizer only started four nuclear power plants and no CCGT units.
Overall, stochastic optimizers maintain more nuclear power plants ON than deterministic ones.
In addition, they tend to start CCGT units during the morning and evening consumption peaks. They anticipate a production shortage for extreme scenarios during consumption peaks occurring during these periods.  Indeed, they consider scenarios where the residual demand is higher than the mean scenario used by the deterministic optimizer.
However, the differences between the several versions of the stochastic optimizers are more challenging to interpret.
Thus, it illustrates the need for KPIs to evaluate these ON/OFF plans.

\end{document}